\documentclass[aos,MSNbibl,nameyear,dvips]{arximspdf}

\doi{10.1214/13-AOS1156} %kopijuoti is PTS
\volume{41}
\issue{5}
\pubyear{2013}
\firstpage{2462}
\lastpage{2504}

\makeatletter

\newcommand{\llangle}{\langle\!\langle}
\newcommand{\rrangle}{\rangle\!\rangle}

\newcommand{\rrvert}{\vert}
\newcommand{\llvert}{\vert}

\newcommand{\cal}{\mathcal}

\newproclaim{remark}{Remark}

\newtheorem{lem}{Lemma}
\newtheorem{theorem}{Theorem}
\newtheorem{cor}{Corollary}
\newtheorem{prop}{Proposition}[section]

\newcommand{\bA}{{\mathbf A}}
\newcommand{\bB}{{\mathbf B}}
\newcommand{\bI}{{\mathbf I}}
\newcommand{\bM}{{\mathbf M}}
\newcommand{\bQ}{{\mathbf Q}}
\newcommand{\bU}{{\mathbf U}}
\newcommand{\bV}{{\mathbf V}}
\newcommand{\bX}{{\mathbf X}}
\newcommand{\be}{{\mathbf e }}
\newcommand{\bg}{{\mathbf g }}

\newcommand{\cA}{{\cal A}}
\newcommand{\cB}{{\cal B}}
\newcommand{\cF}{{\cal F}}
\newcommand{\cG}{{\cal G}}
\newcommand{\cI}{{\cal I}}
\newcommand{\cM}{{\cal M}}
\newcommand{\cP}{{\cal P}}
\newcommand{\cQ}{{\cal Q}}
\newcommand{\cX}{{\cal X}}

\newcommand{\bsigma}{\bolds\sigma}
\newcommand{\brho}{\bolds\rho}

\makeatother

\begin{document}
\begin{frontmatter}

\title{Asymptotic equivalence of quantum state tomography and noisy matrix completion}
\runtitle{Quantum state tomography and noisy matrix completion}

\begin{aug}
\author[A]{\fnms{Yazhen} \snm{Wang}\corref{}\thanksref{t1}\ead[label=e1]{yzwang@stat.wisc.edu}}
\runauthor{Y. Wang}
\affiliation{University of Wisconsin--Madison}
\address[A]{Department of Statistics \\
University of Wisconsin--Madison\\
1300 University Avenue \\
Madison, Wisconsin 53706\\
USA \\
\printead{e1}} %adresu isvedimo komanda gale!
\end{aug}

\thankstext{t1}{Supported in part by the NSF Grants DMS-10-5635 and DMS-12-65203.}

% HISTORY:
\received{\smonth{6} \syear{2012}}
\revised{\smonth{1} \syear{2013}}

% ABSTRACT
%
\begin{abstract}
Matrix completion and quantum tomography are two unrelated research
areas with great current interest %topics
in many modern scientific studies. This paper investigates the
statistical relationship between trace regression in matrix completion
and quantum state tomography in quantum physics and quantum information
science. As quantum state tomography and trace regression share the
common goal of recovering an unknown matrix, it is nature to put them
in the Le Cam paradigm for statistical comparison. Regarding the two
types of matrix inference problems as two statistical experiments, we
establish their asymptotic equivalence in terms of deficiency distance.
The equivalence study motivates us to introduce a new trace regression
model. The asymptotic equivalence provides a sound statistical
foundation for applying matrix
completion methods to quantum state tomography. % under appropriate
%circumstances.
We investigate the asymptotic equivalence for sparse density matrices
and low rank density matrices and demonstrate that sparsity and low
rank are not necessarily helpful for achieving the asymptotic
equivalence of quantum state tomography and trace regression. In
particular, we show that popular Pauli measurements are bad for
establishing the asymptotic equivalence for sparse density matrices and
low rank density matrices.
%%As Pauli measurements are widely used in quantum physics and quantum
%computation,
%The equivalence study inspires us to apply a compressed sensing
%estimation methodology for trace regression to
%quantum state tomography and develop a density matrix estimator for
%quantum state tomography with
%Pauli measurements. Surprisingly the developed density matrix
%estimator works for quantum state tomography
%with Pauli measurements, even when trace regression and quantum state
%tomography are not asymptotically equivalent.
%derive asymptotic equivalence of quantum state tomography and trace
%regression at both scales.
\end{abstract}

% KEYWORDS
\begin{keyword}[class=AMS]
\kwd[Primary ]{62B15}
\kwd[; secondary ]{62P35}
\kwd{62J99}
\kwd{65F10}
\kwd{65J20}
\kwd{81P45}
\kwd{81P50}
\end{keyword}
\begin{keyword}
\kwd{Compressed sensing}
\kwd{deficiency distance}
\kwd{density matrix}
\kwd{observable}
\kwd{Pauli matrices}
\kwd{quantum measurement}
\kwd{quantum probability}
\kwd{quantum statistics}
\kwd{trace regression}
\kwd{fine scale trace regression}
\kwd{low rank matrix}
\kwd{sparse matrix}
\end{keyword}

\end{frontmatter}

%s1 #&#
\section{Introduction}\label{sec1}

Compressed sensing and quantum tomography are two disparate scientific fields.
The fast developing field of compressed sensing provides innovative
data acquisition techniques and supplies efficient
accurate reconstruction methods for recovering sparse signals and
images from highly undersampled observations %with or without noise
[see \citet{Don06}]. Its wide range of applications include signal
processing, medical imaging
and seismology. The problems to solve in compressed sensing often
involve large data sets with complex structures such as
data on many variables or features observed over a much smaller number
of subjects. As a result, the developed theory of compressed sensing
can shed crucial insights on high-dimensional statistics. Matrix
completion, a current research focus point in compressed sensing,
%which grows naturally from traditionally well-established vector
%reconstruction.
is to reconstruct a low rank matrix based on under-sampled
observations. % with or without noise.
Trace regression is often employed in noisy matrix completion for low
rank matrix estimation.
%The recent developments in noisy matrix completion are on the low rank
%matrix estimation based on trace regression.
Recently several methods were proposed
to estimate a low rank matrix by minimizing the squared residual sum
plus some penalty. The penalties used include
nuclear-norm penalty [Cand\'es and Plan (\citeyear{CanPla09,CanPla11}),
\citet{KolLouTsy11} and \citet{NegWai11}],
rank penalty [\citet{BunSheWeg11} and \citet{Klo11}], the
von Neumann entropy penalty [\citet{Kol11}], and
the Schatten-p quasi-norm penalty [\citet{RohTsy11}].

Contemporary scientific studies often rely on understanding and
manipulating quantum systems.
Examples include quantum computation, quantum information and quantum
simulation
%for the purpose of computations and communications, quantum simulation
%in biochemistry and nano-technology
%for the design of bio-molecules and nano-materials, and nuclear fusion
%in astronomy for star analysis
[\citet{NieChu00} and Wang (\citeyear{Wan11}, \citeyear{Wan12})]. The studies
particularly frontier research in quantum
computation and quantum information stimulate great interest in and
urgent demand on quantum tomography.
A~quantum system is described by its state, and the state is often
characterized by a complex matrix %of complex entries
on some Hilbert space. The matrix is called density matrix. % or
%density operator.
A density matrix used to characterize a quantum state usually grows
exponentially
with the size of the quantum system.
%As a result, it takes an exponential number of bits of memory
%on a classical computer to store the state of a quantum system, and
%simulations of quantum systems
%via classic computers face great computational challenge. On the other
%hand, as quantum systems
%are able to store and keep track an exponential number of complex
%numbers
%and perform data manipulations and calculations as the systems evolve,
%quantum
%computation and quantum information are to grapple with understanding
%how to
%take advantage of the enormous information hidden in the quantum
%systems and to
%harness the immense potential computational power of atoms and
%molecules for
%the purpose of information processing and computation.
For the study of a quantum system,
it is important but very difficult to know its state. If we do not know
in advance the
state of the quantum system, we may deduce the quantum state by
performing measurements on
the quantum system. In statistical terminology, we want to estimate the
density matrix based on
measurements performed on a large number of quantum systems which are
identically prepared in
the same quantum state. In the quantum literature, quantum state
tomography refers to the reconstruction
of the quantum state based on measurements obtained %through probing
from measuring identically prepared quantum systems.

In this paper, we investigate statistical relationship between quantum
state tomography and noisy matrix
completion based on trace regression. Trace regression is used to
recover an unknown matrix from noisy
observations on the trace of the products of the unknown matrix and
matrix input variables. Its connection with
quantum state tomography is through quantum probability on quantum
measurements. Consider a finite-dimensional
quantum system with a density matrix. According to the
theory of quantum physics, when we measure the quantum system by
performing measurements on observables which are
Hermitian (or self-adjoint) matrices, the measurement outcomes for each
observable are real eigenvalues of the observable,
and the probability of observing a particular eigenvalue is equal to
the trace of the product of the density matrix and
the projection matrix onto the eigen-space corresponding to the
eigenvalue, with the expected measurement outcome equal to the
trace of the product of the density matrix and the observable. Taking
advantage of the connection \citet{Groetal10} has
applied matrix completion methods with nuclear norm penalization to
quantum state tomography for reconstructing
low rank density matrices. % (actually nuclear norm penalization is not
%suitable for density matrix estimation, since all density matrices
%have nuclear norm equal to one).
As trace regression and quantum state tomography
share the common goal of recovering the same matrix parameter, we
naturally treat them as two statistical models in
the Le Cam paradigm and study their asymptotic equivalence via Le Cam's
deficiency distance. Here equivalence means
that each statistical procedure for one model has a corresponding
equal-performance statistical procedure for another model.
The equivalence study motivates us to introduce a new fine scale trace
regression model. We derive bounds on the deficiency
distances between trace regression and quantum state tomography with
summarized measurement data and between fine scale
trace regression and quantum state tomography with individual
measurement data, and then under suitable conditions we establish
asymptotic equivalence of trace regression and quantum state tomography
for both cases.
The established asymptotic equivalence provides a sound statistical
foundation for applying matrix completion procedures %techniques
to quantum state tomography under appropriate circumstances.
%The obtained results show
%that under some conditions on the sample sizes, trace regression and
%quantum state tomography are asymptotically equivalent
%regardless of the size of the unknown matrix. In other words, the
%asymptotic equivalence is free of the dimensionality
%of the matrix parameter to be estimated, and even if the size of the
%unknown matrix %exceeds sample sizes or even
%grows exponentially faster than simple sizes, we can still obtain the
%asymptotic equivalence of trace regression and quantum
%state tomography. The results stand in contrast to traditional
%asymptotic theory which suffers from the curse of dimensionality.
%The existing asymptotic equivalence
%usually holds for low dimensional statistics such as parametric models
%with fixed finite dimensional parameters and
%nonparametric models with functions in one and two dimensions. This is
%the first asymptotic equivalence result in the ultra-high
%dimension setting.
We further analyze the asymptotic equivalence of trace regression and
quantum state tomography for
sparse matrices and low rank matrices. The detailed analyses indicate
that the asymptotic equivalence does not require sparsity
nor low rank on matrix parameters, and depending on the density matrix
class as well as the set of observables used for performing
measurements, sparsity and low rank may or may not make the asymptotic
equivalence easier to achieve. In particular, we show that the
Pauli matrices as observables are bad for establishing the asymptotic
equivalence for sparse matrices and low rank matrices; and
for certain class of sparse or low rank density matrices, we can obtain
the asymptotic equivalence of quantum state tomography and trace
regression in the ultra high dimension setting where the matrix size of
the density matrices is comparable to or even exceeds the number of the
quantum measurements on the observables.

%Since Pauli matrices are widely used in quantum physics and quantum
%computation to obtain Pauli measurements, inspiring from
%the asymptotic equivalence study we develop a density matrix estimator
%for Pauli measurements
%by applying a compressed sensing estimation methodology for trace
%regression to quantum state tomography. %a method called PET
%Surprisingly the developed density matrix estimator works for quantum
%state tomography with Pauli measurements,
%even when trace regression and quantum state tomography are not
%asymptotically equivalent.
%%The method is based Pauli measurement estimation by thresholding (PET)

The rest of paper proceeds as follows. Section~\ref{sec2} reviews trace
regression and quantum state tomography and states statistical models
and data structures.
We consider only finite square matrices, since trace regression handles
finite matrices, and density matrices %describing quantum states
are square matrices. Section~\ref{sec3} frames trace regression and quantum
state tomography with summarized measurements
as two statistical experiments in Le Cam paradigm and
studies their asymptotic equivalence. Section~\ref{sec4} introduces a
fine scale
trace regression model to match quantum state tomography with
individual measurements
and investigates their asymptotic equivalence. We illustrate the
asymptotic equivalence for sparse density matrix class and low rank
density matrix class in Sections~\ref{sec5} and \ref{sec6}, respectively.
%In particular we show that Pauli measurements are bad for establishing
%the asymptotic equivalence for sparse matrices
%and low rank matrices.
%Section \ref{sec7} applies a compressed sensing estimation methodology
%for
%trace regression to quantum state tomography
%and proposes a density matrix estimator for Pauli measurements. The
%proposed density matrix estimator works for quantum state %tomography
%with Pauli measurements, even when trace regression and quantum state
%tomography are not asymptotically equivalent.
We collect technical proofs in Section~\ref{proofs}, with additional
proofs of technical lemmas in the \hyperref[app]{Appendix}.

%s2 #&#
\section{Statistical models and data structures}\label{sec2}
%%Statistical
%equivalence of trace regression and quantum state tomography}

%s2.1 #&#
\subsection{Trace regression in matrix completion}\label{sec2.1}

Suppose that we have $n$ independent random pairs $(\bX_1, Y_1),\ldots,
(\bX_n, Y_n)$ from the model
%
%e1 #&#
%
\begin{equation}
\label{trace-regression} Y_k = %\langle\langle\bM_k, \brho\rangle\rangle
\operatorname{tr}\bigl(
\bX_k^\dagger\brho\bigr) + \varepsilon_k,\qquad k=1,\ldots, n,
\end{equation}
where $\mathrm{tr}$ is matrix trace, $\dagger$ denotes conjugate transpose,
$\brho$ is an unknown $d$ by $d$ matrix,
$\varepsilon_k$ are zero mean random errors, and $\bX_k$ are matrix
input variables of size $d$ by $d$. %and $\bX_k$ are independent of $
We consider both fixed and random designs. For the random design case,
each $\bX_k$ is randomly sampled from a set of matrices. %which will
%be specified below.
In the fixed design case, $\bX_1,\ldots, \bX_n$ are fixed matrices.
Model (\ref{trace-regression}) is called trace regression and employed
in matrix completion.
Matrix input variables $\bX_k$ are often sparse in a sense that each
$\bX_k$ has a relatively small number of nonzero entries. Trace
regression masks the entries of $\brho$ through $\bX_k^\dagger\brho
$, and each observation $Y_k$ is the trace of the masked
$\brho$ corrupted by noise $\varepsilon_k$. The statistical problem
%%of estimating the unknown matrix $\brho$
is to estimate all the entries of $\brho$ based on observations
$(\bX_k, Y_k)$, $k=1,\ldots, n$, which is often referred to as noisy
matrix completion.
%When $d=1$,
Model (\ref{trace-regression}) and matrix completion are matrix
generalizations of a linear model and sparse signal estimation in
compressed sensing. See Cand\'es and Plan
(\citeyear{CanPla09,CanPla11}), \citet{CanRec09},
\citet{CanTao10}, \citet{KesMonOh10},
\citet{KolLouTsy11}, and \citet{NegWai11},
\citet{Kol11} and \citet{RohTsy11}.

Matrix input variables $\bX_k$ are selected from a matrix set $\cB=\{
\bB_1,\ldots, \bB_p\}$, where $\bB_j$ are $d$ by $d$ matrices.
Below we list some examples of such matrix sets used in matrix completion.
%We use the following labeling convention to link $j$ with $(\ell_1,
%
\begin{longlist}[(iii)]
\item[(i)] Let
%
%e2 #&#
%
\begin{eqnarray}
\label{BasisI}
\cB &=& \bigl\{\bB_{j}=\be_{\ell_1}
\be_{\ell_2}^\prime,
j = (\ell_1-1) d +
\ell_2,\nonumber\\[-8pt]\\[-8pt]
&&\hspace*{5.1pt} j=1,\ldots, p=d^2, \ell_1, \ell
_2=1,\ldots, d \bigr\},\nonumber
\end{eqnarray}
where $\be_{\ell}$ % $=(0,\ldots, 1, \ldots,0)^\prime$
is the canonical basis in Euclid space $\mathbb{R}^d$.
In this case, if $\brho=(\rho_{ab})$, then $\operatorname{tr}(\bB_{j} \brho) =
\rho_{\ell_1 \ell_2}$, and the observation $Y_k$ is equal to
some entry of $\brho$ plus noise $\varepsilon_k$.
More generally, instead of using single $\be_{\ell_1} \be_{\ell
_2}^\prime$,
we may define $\bB_j$ as the sum of several $\be_{\ell_1} \be_{\ell
_2}^\prime$, and then
$\operatorname{tr}(\bB_{j} \brho)$ is equal to the sum of some entries of $\brho$.
\item[(ii)] Set
%
%e3 #&#
%
\begin{equation}
\label{Hermitian-basis} \cB=\bigl\{\bB_{j}, j = 1,\ldots,
p=d^2\bigr\},
\end{equation}
where we identify $j$ with $(\ell_1, \ell_2)$, $j =1,\ldots, p$,
$\ell_1,\ell_2 = 1,\ldots, d$,
$\bB_{j} = \be_{\ell_1} \be_{\ell_2}^\prime$ for $\ell_1=\ell_2$,
\[
\bB_{j} = \frac{1}{\sqrt{2}} \bigl( \be_{\ell_1}
\be_{\ell_2}^\prime+ \be_{\ell_2} \be_{\ell
_1}^\prime
\bigr) \qquad\mbox{for } \ell_1 <\ell_2
\]
and
\[
\bB_{j} = \frac{\sqrt{-1}}{\sqrt{2}} \bigl( \be_{\ell_1} \be
_{\ell_2}^\prime- \be_{\ell_2} \be_{\ell_1}^\prime
\bigr) \qquad\mbox{for } \ell_1 > \ell_2.
\]
\item[(iii)] For $d=2$ define
\begin{eqnarray*}
\bsigma_0 &=& \pmatrix{ 1 & 0
\cr
0 & 1},\qquad \bsigma_1 =
\pmatrix{ 0 & 1
\cr
1 & 0 },\\
\bsigma_2 &=& \pmatrix{ 0 & -\sqrt{-1}
\cr
\sqrt{-1} & 0 },\qquad \bsigma_3 = \pmatrix{ 1 & 0
\cr
0 & -1 },
\end{eqnarray*}
%
% \bsigma_1=\left( \begin{array}{cc} 0 & 1 \\ 1 & 0 \end{array}
% \bsigma_2=\left( \begin{array}{cc} 0 & \sqrt{-1} \\ -\sqrt{-1} & 0
% \bsigma_3=\left( \begin{array}{cc} 1 & 0 \\ 0 & -1 \end{array}
%The four Pauli matrices form a basis in the space $\mathbb{C}^{2\times
%2}$ of all $2$ by $2$ complex matrices.
where $\bsigma_1$, $\bsigma_2$ and $\bsigma_3$ are called the Pauli matrices.
For $d=2^b$ with integer $b$, we may use $b$-fold tensor products of
$\bsigma_0$, $\bsigma_1$, $\bsigma_2$ and $\bsigma_3$ to define
general Pauli matrices and obtain the Pauli matrix set
%$\cB$ consists of all the $b$-fold tensor products of $\sigma_\ell$, $
%Pauli basis for the space $\mathbb{C}^{d\times d}$ of all $d$ by $d$
%complex matrices and let $\cB$ consisting of all Pauli basis matrices,
%
%e4 #&#
%
\begin{equation}
\label{Pauli-basis} \cB= \bigl\{ %2^{-b/2}
\bsigma_{\ell_1} \otimes
\bsigma_{\ell_2} \otimes\cdots\otimes\bsigma_{\ell_b}, (
\ell_1, \ell_2,\ldots, \ell_b) \in\{0, 1, 2,
3\}^b \bigr\},
\end{equation}
where $\otimes$ denotes tensor product. The Pauli matrices %in (
%to represent complex Hermitian matrices
are widely used in quantum physics and quantum information science.
\end{longlist}

Matrices in (\ref{BasisI}) are of rank 1 and have eigenvalues $1$ and $0$.
For matrices in~(\ref{Hermitian-basis}), the diagonal matrices are of
rank $1$ and have eigenvalues $1$ and $0$,
and the nondiagonal matrices are of rank $2$ and have eigenvalues $\pm
1$ and $0$. Pauli matrices in (\ref{Pauli-basis}) are of full rank,
and except for the identity matrix all have eigenvalues $\pm1$.
Denote by $\mathbb{C}^{d\times d}$ the space of all $d$ by $d$ complex
matrices and define an inner product %equipped with the trace norm
$\llangle\bA_1, \bA_2\rrangle=\operatorname{tr}( \bA^\dagger_2
\bA_1)$ for $\bA_1, \bA_2 \in\mathbb{C}^{d\times d}$. Then both
(\ref{Hermitian-basis}) and (\ref{Pauli-basis}) form orthogonal bases
for all complex Hermitian matrices, and the real matrices
in (\ref{Hermitian-basis}) or (\ref{Pauli-basis}) form orthogonal
bases for all real symmetric matrices.

For the random design case, with $\cB=\{\bB_j, j=1,\ldots, p\}$, we
assume that matrix input variables $\bX_k$ are independent and
%identically distributed as matrix random variable $\bX$ which is
%sampled from $\cB$ according to
sampled from $\cB$ according to a distribution $\Pi(j)$ on $\{1,\ldots,
p\}$,
%
%e5 #&#
%
\begin{equation}
\label{design-Pi} P(\bX_k =\bB_{j_k}) =
\Pi(j_k),\qquad k=1,\ldots, n, j_k \in\{ 1,\ldots, p\}.
\end{equation}
The observations from (\ref{trace-regression}) are $(\bX_k, Y_k)$,
$k=1,\ldots, n$, with $\bX_k$ sampled from $\cB$ according to the
distribution $\Pi(\cdot)$.
For the fixed design case, matrix input variables $\bX_1,\ldots, \bX
_n$ form a fixed set of matrices, and we assume $n=p$ and
$\cB=\{ \bX_1,\ldots, \bX_n\}=\{\bB_1,\ldots, \bB_p\}$. The
observations from (\ref{trace-regression}) are $(\bX_k, Y_k)$,
$k=1,\ldots, n$, with deterministic $\bX_k$.

%s2.2 #&#
\subsection{Quantum state and measurements}\label{sec2.2}
%Quantum mechanics}

%Quantum mechanics describes microscopic phenomena such as position and
%momentum of an individual particle like an atom or electron, spin of an
%electron, detection of light photons, and the emission and absorption
%of light by atoms.
%A quantum system is characterized by its state and the time evolution
%of the state.
%The theory of quantum mechanics is mathematically described by a
%Hilbert space and self-adjoint
%(or Hermitian) operators. In this paper we consider the
%finite-dimensional situation.
%A quantum system is characterized by its state,
For a finite-dimensional quantum system, we describe its quantum state
by a density matrix $\brho$ on $d$-dimensional complex
space $\mathbb{C}^d$, where density matrix $\brho$ is a $d$ by $d$
complex matrix satisfying (1) Hermitian, that is, $\brho$ is equal to
its conjugate transpose;
(2) semi-positive definite; (3) unit trace, that is, $\operatorname{tr}(\brho)=1$.
%Following the convention in quantum information science we reserve
%notation $\brho$ for state or density matrix.
%See Griffiths (2004), \citet{SakNap10} and \citet{Sha94}.
%Denote by superscripts $*$, $\prime$ and $\dagger$ the conjugate of a
%complex number, the transpose
%of a vector or matrix, and conjugate transpose operation, respectively.

Experiments are conducted to perform measurements on the quantum system
and obtain data for studying the quantum system. %and test the theory
%of quantum mechanics.
Common quantum measurements are on some observable %such as position,
%momentum and spin, where an {\it observable} $\bM$
$\bM$, which is defined as a Hermitian matrix on $\mathbb{C}^d$.
Assume that the observable $\bM$ has the following spectral
decomposition:
%a discrete spectrum with the following diagonal form,
%
%e6 #&#
%
\begin{equation}
\label{diagonal} \bM= \sum_{a=1}^r
\lambda_a \bQ_a,
\end{equation}
where $\lambda_a$ %\in\mathbb{R}$
are $r$ different real eigenvalues of $\bM$, and $\bQ_a$ are
projections onto the eigen-spaces %of $\bM$
corresponding to $\lambda_a$.
For the quantum system prepared in a state $\brho$, we need a
probability space $(\Omega, {\cal F}, P)$ to describe
measurement outcomes when performing measurements on the observable
$\bM$.
Denote by $R$ the measurement outcome of $\bM$. According to the
theory of quantum
mechanics, $R$ is a random variable on $(\Omega, {\cal F}, P)$ taking
values in $\{\lambda_1, \lambda_2,\ldots, \lambda_r\}$,
with probability distribution given by
%
%e7 #&#
%
\begin{equation}
\label{measurement} P(R = \lambda_a) =\operatorname{tr}(\bQ_a
\brho),\qquad
a=1, 2,\ldots, r,\qquad E(R)= \operatorname{tr}(\bM\brho).
\end{equation}
%
%With the probability we derive the expectation under pure state $|
%We conduct experiments of performing measurements on quantum systems.
%%in the laboratory.
See \citet{Hol82}, \citet{SakNap10}, \citet{Sha94} and
\citet{Wan12}.

Suppose that an experiment is conducted to perform measurements on $\bM
$ independently for $m$
quantum systems which are identically prepared in the same quantum
state $\brho$. From the experiment we obtain individual measurements
$R_1,\ldots, R_m$,
which are i.i.d. according to distribution (\ref{measurement}), and
denote their average by $N=(R_1 + \cdots+ R_m)/m$.

The following proposition provides a simple multinomial
characterization for the distributions of $(R_1,\ldots, R_m)$ and $N$.
%
%pr2.1 #&#
%
\begin{prop} \label{prop1}
As random variables $R_1,\ldots, R_m$ take eigenvalues $\lambda
_{1},\ldots, \lambda_r$, we count
the number of $R_1,\ldots, R_m$ taking $\lambda_a$ and define the
counts by $U_a= \sum_{\ell=1}^m 1(R_\ell=\lambda_a)$, $a=1,\ldots, r$.
Then %(\ref{measurement}) indicates that
the counts $U_1,\ldots, U_r$ jointly follow the following multinomial
distribution:
%with $m$ trials and event probabilities $\operatorname{tr}(\bQ_{a} \brho)$, $a=1,
%
%e8 #&#
%
\begin{eqnarray}
\label{measurement-U} P(U_1 = u_1,\ldots,
U_r = u_r) &=& \pmatrix{ m
\cr
u_1,\ldots,
u_r} \bigl[\operatorname{tr}(\bQ_1 \brho)\bigr]^{u_1} \cdots
\bigl[\operatorname{tr}(\bQ_r \brho)\bigr]^{u_r},\nonumber\\[-8pt]\\[-8pt]
\sum_{a=1}^r u_a &=& m\nonumber
\end{eqnarray}
and
%
%e9 #&#
%
\begin{equation}
\label{measurement-RU} N= (R_1 + \cdots+ R_m)/m = (
\lambda_1 U_1 + \cdots+ \lambda_a
U_a)/m.
\end{equation}
\end{prop}

We note the difference between the observable $\bM$ which is a
Hermitian matrix and its measurement result $R$ which is a real-valued
random variable. To illustrate the connection between density matrix
$\brho$ and the measurements of $\bM$, % and emphasize the difference
%from usual statistics,
we assume that $\bM$ has $d$ different eigenvalues. As in
\citet{ArtGilGut05}, we use the normalized eigenvectors of
$\bM$ to form an orthonormal basis,
%in the class of complex Hermitian matrices.
represent $\brho$ under the basis and denote the resulting matrix by
$(\rho_{\ell_1 \ell_2})$. Then from (\ref{measurement}) we obtain
%$\brho= (\be_j) (\rho_{ij})$. With this special representation we have
%$\bQ_a (\be_a)$ is equal to a matrix whose elements are all zero
%except for $1$ at the $(a,a)$, and thus
%we have following simple expression for $P_{\brho}$:
%
\[
P(R = \lambda_a) =\operatorname{tr}(\bQ_a \brho) %= \operatorname{tr}( \bQ_a (\be_j) (\rho_{ij}))
=
\rho_{aa},\qquad a=1, 2,\ldots, d.
\]
%
%For example, if $\bM$ itself is a projection operator??,
That is, with the representation under the eigen basis of $\bM$,
measurements on single observable $\bM$ contain only information about
the diagonal elements of $(\rho_{\ell_1 \ell_2})$. No matter how
many measurements we perform on $\bM$, we cannot draw any inference
about the off-diagonal elements of $(\rho_{\ell_1 \ell_2})$ based on
the measurements on $\bM$. We usually need to perform measurements on
enough different observables %(at least $d+1$)
in order to estimate the whole density matrix $(\rho_{\ell_1 \ell_2})$.
%With the unit trace constraint on density matrices, the class of
%density matrices have dimension $d^2-1$, and the number of obserables
%need to be at least $d+1$.
See \citet{ArtGilGut05}, \citet{BarGilJup03} and
\citet{ButGutArt07}.

%As dim(Sd) = d2 - 1 = (d - 1)(d + 1), to have a one-to-one map between
%states and probability distributions of results we must measure, on
%many identical systems, at least d + 1 different observables. The
%probing of identically prepared quantum systems from different
%'angles' to reconstruct their state is named quantum state tomography
%in the physics literature.

%s2.3 #&#
\subsection{Quantum state tomography}\label{sec2.3}
In physics literature quantum state tomography refers to the
reconstruction of a quantum state based on measurements
obtained from quantum systems that are identically prepared under the
state. %an experiment,
Statistically it is the problem of estimating the density matrix from
the measurements.
%It is very important in quantum applications such as quantum
%computation and quantum information.
Suppose that quantum systems are identically prepared in a state $\rho
$, $\cB=\{\bB_1,\ldots, \bB_p\}$
is a set of observables available to perform measurements, and each
$\bB_j$ has a spectral decomposition
%
%e10 #&#
%
\begin{equation}
\label{basis-diagonal} \bB_j= \sum_{a=1}^{r_j}
\lambda_{ja} \bQ_{ja},
\end{equation}
where $\lambda_{ja}$ are $r_j$ different real eigenvalues of $\bB_j$,
and $\bQ_{ja}$ are projections onto the eigen-spaces
corresponding to $\lambda_{ja}$. We select an observable, say $\bB_j
\in\cB$, and perform measurements on $\bB_j$ for
the quantum systems. According to the observable selection we classify
the quantum state tomography experiment as either a fixed design or a
random design.
%We may consider both random experiment and fixed experiment for
%performing measurements.
In a random design, we choose an observable at random from $\cB$ to
perform measurements for the quantum systems,
while a fixed design is to perform measurements on every observable in
$\cB$ for the quantum systems.

Consider the random design case. We sample an observable $\bM_k$ from
$\cB$ to perform measurements independently
for $m$ quantum systems, $k=1,\ldots, n$, %by sampling an index $J_k$
%from $\{1, 2,\ldots, p\}$
where observables $\bM_1,\ldots, \bM_n$ are independent and
%%identically distributed as random matrix $\bM$ with distribution
sampled from $\cB$ according to a distribution $\Xi(j)$ on $\{1,\ldots,
p\}$,
%
%e11 #&#
%
\begin{equation}
\label{design-Xi} P(\bM_k = \bB_{j_k}) =
\Xi(j_k),\qquad k=1,\ldots, n, j_k \in\{ 1,\ldots, p\}.
\end{equation}
%
%according to a probability distribution $\Xi(j)$ over $\{1, 2,\ldots,
%p\}$ to
Specifically we perform measurements on each observable $\bM_k$
independently for $m$ quantum systems that are identically
prepared under the state $\brho$, and denote by $R_{k1},\ldots,
R_{km}$ the $m$ measurement outcomes and $N_k$ the average
of the $m$ measurement outcomes. The resulting individual measurements
are the data $(\bM_k, R_{k1},\ldots, R_{km})$,
$k=1,\ldots, n$, and the summarized measurements are the pairs $(\bM
_k, N_k)$, $k=1,\ldots, n$, where
%$J_1,\ldots, J_n$ are independent and identically distributed as
%random variable $J$ with probability distribution
%
%e12 #&#
%
\begin{equation}
\label{tomography0} N_k=(R_{k1} + \cdots+
R_{km})/m,
\end{equation}
%
%$R_{k1},\ldots, R_{km}$ are the results of measuring $\bM_{k}$
%independently $m$ times
$R_{k\ell}$, $k=1,\ldots, n$, $\ell=1,\ldots, m$, are
independent, and given
$\bM_k=\bB_{j_k}$ for some $j_k \in\{1,\ldots, p\}$,
the conditional distributions of $R_{k1},\ldots, R_{km}$ are given by
%
%e13 #&#
%e14 #&#
%
\begin{eqnarray}
\label{tomography1} P(R_{k\ell}=\lambda_{j_k a}|
\bM_k=\bB_{j_k}) &=& \operatorname{tr}(\bQ_{j_k a} \brho),\nonumber\\[-8pt]\\[-8pt]
&&\eqntext{a=1,\ldots, r_{j_k}, \ell
=1,\ldots, m, j_k \in\{1,\ldots, p\},}
\\
\label{tomography2}
E(R_{k\ell}|\bM_k=\bB_{j_k}) &=& \operatorname{tr} (
\bB_{j_k} \brho), \nonumber\\[-8pt]\\[-8pt]
\operatorname{Var}(R_{k\ell}|\bM_k=
\bB_{j_k}) &=& %\left\{
\operatorname{tr} \bigl(\bB_{j_k}^2 \brho
\bigr) - \bigl[\operatorname{tr} (\bB_{j_k} \brho) \bigr]^2.\nonumber
% \right\} /m. &
\end{eqnarray}
The statistical problem is to estimate $\brho$ from the individual
measurements $(\bM_k, R_{k1},\ldots, R_{km})$, $k=1,\ldots, n$,
or from the summarized measurements $(\bM_1, N_1),\ldots, (\bM_n, N_n)$.

For the fixed design case, we take $p=n$ and $\cB=\{\bB_1,\ldots,
\bB_n\}$.
%we fix a subset $\cB_0$ of $\cB$ such that the subset $\cB_0$ has $n$
%observables.
%For example we may assume that the observables in $\cB_0$ are $\bM_1,
%,\ldots, \bM_n$.
We perform measurements on every observable $\bM_k=\bB_k \in\cB$
independently for $m$ quantum systems that are identically
prepared under the state $\brho$, and denote by $R_{k1},\ldots,
R_{km}$ the $m$ measurement outcomes and $N_k$ the average
of the $m$ measurement outcomes. The resulting individual measurements
are the data $(\bM_k, R_{k1},\ldots, R_{km})$,
$k=1,\ldots, n$, and the summarized measurements are the pairs $(\bM
_k, N_k)$, $k=1,\ldots, n$, where $N_k$ is the same
as in (\ref{tomography0}),
% N_k=(R_{k1} + \cdots+ R_{km})/m,
$R_{k\ell}$, $k=1,\ldots, n$, $\ell=1,\ldots, m$, are
independent, and the distributions of $R_{k1},\ldots, R_{km}$ are
given by
%
%e15 #&#
%e16 #&#
%
\begin{eqnarray}
\label{tomography4} P(R_{k\ell}=\lambda_{ka}) &=& \operatorname{tr}(
\bQ_{ka} \brho),\qquad a=1,\ldots, r_k, \ell=1,\ldots, m,
\\
\label{tomography5}
E(R_{k\ell}) &=& \operatorname{tr} (\bM_k \brho),\qquad \operatorname{Var}(R_{k\ell
}) =
\operatorname{tr} \bigl(\bM_k^2 \brho\bigr) - \bigl[\operatorname{tr} (
\bM_k \brho) \bigr]^2.% \right\}/m.
\end{eqnarray}
The statistical problem is to estimate $\brho$ from the individual
measurements $(\bM_k, R_{k1},\ldots, R_{km})$, $k=1,\ldots, n$,
or from the summarized measurements $(\bM_1, N_1),\ldots, (\bM_n, N_n)$.

Because of convenient statistical procedures and fast implementation algorithms,
the summarized measurements instead of the individual measurements are
often employed in quantum state tomography
[\citet{Groetal10}, \citet{Kol11},
\citet{NieChu00}]. However, in Section~\ref{section-fine} we
will show that quantum state tomography based
on the summary measurements may suffer from substantial loss of
information, and we can develop more efficient
statistical inference procedures by the individual measurements than by
the summary measurements.

In order to estimate all $d^2-1$ free entries of $\brho$, we need the
quantum state tomography model identifiable. Suppose that
all $\bB_j$ have exact $r$ distinct eigenvalues. The identifiability
may require $n \geq(d^2 -1)/(r-1)$ (which is at least $d+1$) and $m
\geq r-1$
for the individual measurements and $n \geq d^2 -1$ for the summarized
measurements.
There is a trade-off between $r$ and $m$ in the individual measurement
case. For large $r$, we need less observables but more measurements on
each observable, while for small $r$, we require more observables but
less measurements on each observable.
In terms of the total number, $m n$, of measurement data, the
requirement becomes $m n \geq d^2-1$. %, regardless of $r$.
%how many different eigenvalues the observables have.

%If the models are not identifiable, we may reparametrize $\brho$
%through $\brho\bQ_{ja}$ and consider estimable linear functionals
%$\operatorname{tr}(\brho\bQ_{ja})$, which is similar to a linear model where the
%design matrix is not full rank and we consider
%estimable linear functions of parameters.

%s3 #&#
\section{Asymptotic equivalence}\label{sec3}
\label{section-P-equivalence}
Quantum state tomography and trace regression share the common goal of
estimating the same unknown matrix
$\brho$, and it is nature to put them in the Le Cam paradigm for
statistical comparison.
%based on data observed on random matrices and random variables.
We compare trace regression and quantum state tomography in either the
fixed design case or the random design case.

First, we consider the fixed design case.
Trace regression (\ref{trace-regression}) generates data on dependent
variables $Y_k$ with deterministic
matrix input variables $\bX_k$, and we denote by $\mathbb
{P}_{1,n,\brho}$ the joint distribution of
$Y_k$, $k=1,\ldots, n$. Quantum state tomography performs
measurements on a fixed set of observables
$\bM_k$ and obtains average measurements $N_k$ on $\bM_k$ whose
distributions are specified by
(\ref{tomography0}) and (\ref{tomography4})--(\ref{tomography5}), and
we denote by $\mathbb{P}_{2,n,\brho}$ the joint distribution of
$N_k$, $k=1,\ldots, n$.
Both $\mathbb{P}_{1,n,\brho}$ and $\mathbb{P}_{2,n,\brho}$ are probability
distributions on measurable space $(\mathbb{R}^n, \cF_\mathbb
{R}^n)$, where $\cF_\mathbb{R}$ is the Borel $\sigma$-field
on $\mathbb{R}$.

Second we consider the random design case. Trace regression (\ref
{trace-regression}) generates data on the
pairs $(\bX_k, Y_k)$, $k=1,\ldots, n$, where matrix input variables
$\bX_k$ are sampled from $\cB$ according to the
distribution $\Pi(j)$ given by (\ref{design-Pi}). We denote by
$\mathbb{P}_{1,n,\brho}$ the joint distribution of
$(\bX_k, Y_k)$, $k=1,\ldots, n$, for the trace regression model.
Quantum state tomography yields observations in the
form of observables $\bM_{k}$ and average measurement results $N_k$ on
$\bM_{k}$, $k=1,\ldots, n$, where the distributions of
$(\bM_{k}, N_k)$ are specified by (\ref{design-Xi})--(\ref
{tomography2}). We denote by $\mathbb{P}_{2,n,\brho}$ the
joint distribution of $(\bM_{k}, N_k)$, $k=1,\ldots, n$, for the
quantum state tomography model.
Both $\mathbb{P}_{1,n,\brho}$ and $\mathbb{P}_{2,n,\brho}$ are
probability distributions
on measurable space $(\cB^n \times\mathbb{R}^n, \cF_\cB^n \times
\cF_\mathbb{R}^n)$, where $\cF_\cB$ consists of all subsets of $\cB$.
% and $\cF_\mathbb{R}$ is the Borel $\sigma$-field on $\mathbb{R}$.
%Denote by $\mathbb{P}_{1,n,\brho}$ the joint distributions of
%observations $(Y_k, \bX_k)$
%from the trace regression model with $\varepsilon_k \sim N(0, \operatorname{tr}(\bX_k

Denote by $\Theta$ a class of semi-positive Hermitian matrices with
unit trace. %which will be specifically given in Condition ((C3)) below.
For trace regression and quantum state tomography, we define two
statistical models
%
%e17 #&#
%
\begin{equation}
\label{experiment}\quad \cP_{1n} = \bigl\{(\cX_1,
\cG_1, \mathbb{P}_{1,n,\brho}), \brho\in\Theta\bigr\},\qquad
\cP_{2n} = \bigl\{(\cX_2, \cG_2,
\mathbb{P}_{2,n,\brho}), \brho\in\Theta\bigr\},
\end{equation}
where measurable spaces $(\cX_i, \cG_i)$, $i=1,2$, are either $(\cB
^n \times\mathbb{R}^n, \cF_\cB^n \times\cF_\mathbb{R}^n)$ for
the random design case or $(\mathbb{R}^n, \cF_\mathbb{R}^n)$ for the
fixed design case.
% and $\cF_2$ be two statistical models, respectively, corresponding to
%$\mathbb{P}_{1,\brho}$ and $\mathbb{P}_{2,\brho}$, $ \brho\in\Theta$.
Models $\cP_{1n}$ and $\cP_{2n}$ are called statistical experiments in
the Le Cam paradigm. We use Le Cam's deficiency distance between
$\cP_{1n}$ and $\cP_{2n}$ to compare the two models. Let $\cA$ be a
measurable action space, $L$: $\Theta\times\cA\rightarrow[0, \infty)$ a
loss function, and $\|L\| = \sup\{ L(\brho, \mathbf{a})\dvtx
\brho\in\Theta, \mathbf{a} \in\cA\}$. For model $\cP_{in}$, $i=1, 2$,
denote by $\chi_i$ a decision procedure and $R_i(\chi_i, L, \brho)$ the
risk from using procedure $\chi_i$ when $L$ is the loss function and
$\brho$ is the true value of the parameter. We define deficiency
distance $\Delta(\cP_{1n}, \cP_{2n})$ between $\cP_{1n}$ and $\cP_{2n}$
as the maximum of $\delta(\cP_{1n}, \cP_{2n})$ and $\delta(\cP_{2n},
\cP_{1n})$, where
\[
\delta(\cP_{1n}, \cP_{2n}) = \inf_{\chi_1}
\sup_{\chi_2} \sup_{\brho\in\Theta} \sup
_{L: \|L\|=1} \bigl|R_1(\chi_1, L, \brho) -
R_2(\chi_2, L, \brho)\bigr|
\]
is referred to as the deficiency of $\cP_{1n}$ with respect to $\cP_{2n}$.
If $\Delta(\cP_{1n}, \cP_{2n}) \leq\epsilon$, then every decision
procedure in one of the two experiments $\cP_{1n}$ and $\cP_{2n}$ has
a corresponding
procedure in another experiment that comes within $\epsilon$ of
achieving the same risk for any bounded loss.
%Two experiments $\cP_{1n}$ and $\cP_{2n}$ are called equivalent if $
%in problem $\cP_{1n}$ has a
%corresponding procedure $\xi_2$ in problem $\cP_{2n}$ with the same
%risk, uniformly over $\brho\in\Theta$ and all bounded loss $L$, and
%vice versa.
Two sequences of statistical experiments $\cP_{1n}$ and $\cP_{2n}$
are called asymptotically equivalent if
$\Delta(\cP_{1n},\cP_{2n}) \rightarrow0$,
as $n \rightarrow\infty$. For two asymptotic equivalent experiments
$\cP_{1n}$ and $\cP_{2n}$, any sequence of procedures
$\chi_{1n}$ in model $\cP_{1n}$ has a corresponding sequence of
procedures $\chi_{2n}$ in model $\cP_{2n}$ with risk differences
tending to zero uniformly over $\brho\in\Theta$ and all loss $L$
with $\|L\| = 1$,
and the procedures $\chi_{1n}$ and $\chi_{2n}$ are called
asymptotically equivalent.
See \citet{LeC86}, \citet{LeCYan00} and \citet{Wan02}.

To establish the asymptotic equivalence of trace regression and quantum
state tomography, we need to lay down technical conditions and
make some synchronization arrangement between observables in quantum
state tomography and matrix input variables in trace regression.
%match with matrix variables $\bX_k$ in trace regression.
%
\begin{longlist}
\item[(C1)] Assume that $\cB=\{\bB_1,\ldots, \bB_p\}$, and each
$\bB_j$ is a Hermitian matrix with at most $\kappa$
distinct eigenvalues, where $\kappa$ is a fixed integer.
%With the spectral decomposition (\ref{basis-diagonal}) for each $\bB_j
%consistently estimated from trace regression or quantum state
%tomography.
Matrix input variables $\bX_k$ in trace regression and observables
$\bM_k$ in quantum state tomography
are taken from~$\cB$. For the fixed design case, we assume $p=n$, and
$\bX_k=\bM_k=\bB_k$, $k=1,\ldots, n$.
For the random design case, $\bX_k$ and $\bM_k$ are independently
sampled from $\cB$ according to
distributions $\Pi(j)$ and $\Xi(j)$, respectively, and assume that as
$n, p \rightarrow\infty$, $n \gamma_p \rightarrow0$,
where
%
%e18 #&#
%
\begin{equation}
\label{Pi} \gamma_p = \max_{1 \leq j \leq p} \biggl[ \biggl
\llvert1 - \frac{\Pi
(j)}{\Xi(j)} \biggr\rrvert+ \biggl\llvert1 - \frac{\Xi(j)}{\Pi(j)}
\biggr\rrvert\biggr].
\end{equation}
%
%most $\kappa$ distinct eigenvalues and
%form an orthonormal basis in $\mathbb{C}^{d\times d}$.
%
\item[(C2)] Suppose that two models $\cP_{1n}$ and $\cP_{2n}$ are
identifiable. For trace regression, we assume that $(\bX_1,
\varepsilon_1),\ldots, (\bX_n, \varepsilon_n)$ are independent,
and given $\bX_k$, $\varepsilon_k$ follows a normal
distribution with mean zero and variance
%
%e19 #&#
%
\begin{equation}
\label{C2-var} \operatorname{Var}(\varepsilon_k|\bX_k) =
\frac{1}{m} \bigl\{\operatorname{tr}\bigl(\bX_k^2 \brho\bigr) -
\bigl[\operatorname{tr}(\bX_k \brho)\bigr]^2\bigr\}.
\end{equation}
%
%the statistical inference problem for $\brho$ is identifiable for
%trace regression based on
%$(\bX_1,Y_1),\ldots, (\bX_n,Y_n)$ or quantum state tomography based
%on $(\bM_1,N_1),\ldots, (\bM_n,N_n)$.
%can be expressed linearly by $\{\bB_j, j=1,\ldots, p\}$ so that it
%can be consistently estimated from trace regression or quantum state
%tomography.
%
\item[(C3)]
%Denote by $\mathbb{H}^{d\times d}$ the class of all $d$ by $d$
%semi-positive Hermitian matrices with unit trace, and $\Theta$ a
%subset of $\mathbb{H}^{d\times d}$.
For $\bB_j \in\cB$ with spectral decomposition (\ref
{basis-diagonal}), $j=1,\ldots, p$, let
%
%e20 #&#
%
\begin{equation}
\label{cI} % p_{ja}(\brho) = \operatorname{tr}(\bQ_{ja} \brho), a=1,\ldots, r_j,
\cI_j(\brho) = \bigl\{a\dvtx0 < \operatorname{tr}(
\bQ_{ja} \brho) <1, 1 \leq a \leq r_j \bigr\}.
\end{equation}
%
%0 < p_{ja}(\brho) <1, 1 \leq a \leq\kappa, j=1,\ldots, p \}. \]
Let $c_0$ and $c_1$ be two fixed constants with $0<c_0 \leq c_1 <1$.
Assume for $\brho\in\Theta$,
%
%e21 #&#
%
\begin{equation}
\label{theta} c_0 \leq\min_{a \in\cI_j(\brho)} \operatorname{tr}(
\bQ_{ja} \brho) \leq\max_{a \in\cI_j(\brho)} \operatorname{tr}(\bQ_{ja}
\brho) \leq c_1,\qquad j=1,\ldots, p. % \Theta=\left\{\brho\in\mathbb{H}^{d
%%c_0 \leq\min_{(j,a) \in\cI} p_{ja} \leq\max_{(j,a) \in\cI} p_{ja}
% c_0 \leq\min_{a \in\cI_j(\brho)} p_{ja}(\brho) \leq\max_{a \in
% \right\},
\end{equation}
%
%and $|\cI_j(\brho)|$ denotes the cardinal number of $\cI_j(\brho)$.
%and $s$ is an integer to specify sparsity of $\brho$.
%the subset of $\mathbb{C}^{d\times d}$ consisting of all $d$ by $d$
%Hermitian matrices.
\end{longlist}

%re1 #&#
%
\begin{remark}\label{Rem1}
Condition (C1) synchronizes matrices used as matrix
input variables in trace regression and as observables in
quantum state tomography so that we can compare the two models. The
synchronization is needed for applying matrix completion
methods to quantum state tomography [\citet{Groetal10}]. The
finiteness assumption on $\kappa$ is due to the practical
consideration. Observables in quantum state tomography and matrix input
variables in trace regression are often of
large size. Mathematically the numbers of their distinct eigenvalues
could grow with the size, however, in practice matrices with a few
distinct eigenvalues are usually chosen as observables to perform
measurements in quantum state tomography and as
matrix input variables to mask the entries of $\brho$ in matrix
completion [\citet{CanRec09}, \citet{Gro11},
\citet{Groetal10}, \citet{Kol11}, \citet{KolLouTsy11}, \citet{NieChu00}, \citet{Rec11}, \citet{RohTsy11}].
%As we have discussed in Section \ref{sec2}, $p$ is the number of
%matrix input
%variables in trace regression and the number of observables in quantum
%state tomography, and in general we need $p$ to be
%large enough (at least $d+1$) so that it is possible to estimate all
%entries of a density matrix. The main purpose of this paper is
%to study asymptotic equivalence of quantum state tomography and trace
%regression but not how good or bad are density matrix estimators,
%%means that each statistical procedure for one model has a
%corresponding equal-performance statistical procedure for another
%model.
%so Condition (C1) does not impose any specific requirement on $p$ and
%$d$.
%Multinomial distributions are used to describe the distributions of
%$N_k$ in quantum state tomography, and multivariate
%normal distributions are employed to bridge the multinomial
%distribution with trace regression.
Condition (C2) is to match the variance of $N_k$ in quantum state
tomography with the variance of random error $\varepsilon_k$
in trace regression in order to obtain the asymptotic equivalence,
since $N_k$ and $Y_k$ always have the same mean.
%As we have discussed in Section \ref{sec2.2}, to recover all entries
%of $\brho$
%we require both models identifiable in (C2).
%%on $\bB_j$ as observables in quantum state tomography as well as
%matrix input variables in trace regression. Condition (C2) imposes
%%%such a requirement.
%The identifiability imposes some restrictions on $\bB_j$ as well as
%$(n,d)$.
%For example, consider the case that $p=d^2$, $\bB_1,\ldots, \bB_p$
%form an orthonormal basis, and each $\bB_j$ has exact $\kappa$
%distinct eigenvalues. In this case, for both trace regression and
%quantum state tomography we could possibly recover $(\kappa-1) n$
%entries of $\brho$ based on quantum measurements on $n$ different
%observables $\bB_j$. For $\brho$ with $d^2-1$ free entries, we may
%%require $(\kappa-1) n \geq d^2-1$, that is, with finite $\kappa$, $n$
%is at least of order $d^2$.
Regarding condition (C3),
from (\ref{measurement-U})--(\ref{measurement-RU}) and (\ref
{tomography0})--(\ref{tomography5}) we may see that each $N_k$ is
determined by the counts of random variables $R_{k\ell}$ taking
eigenvalues $\lambda_{ja}$, and the counts jointly follow
a multinomial distribution with parameters of $m$ trials and cell
probabilities $\operatorname{tr}(\bQ_{ja} \brho)$, $a=1,\ldots, r_j$.
Condition (C3) is to ensure that the multinomial distributions (with
uniform perturbations) can be well approximated by
multivariate normal distributions %with small enough approximation
%errors
so that we can calculate the Hellinger distance between the
distributions of $N_k$ (with uniform perturbations) in quantum state
tomography and the
distributions of $\varepsilon_k$ in trace regression and thus
establish the asymptotic equivalence of quantum state tomography and
trace regression.
Index $\cI_j(\brho)$ in (\ref{cI}) is to exclude all the cases with
$\operatorname{tr}(\bQ_{ja} \brho)=0$ or $\operatorname{tr}(\bQ_{ja} \brho)=1$,
under which measurement results on $\bB_j$ are certain, either never
yielding measurement results $\lambda_{ja}$ or always
yielding results $\lambda_{ja}$, and their contributions to $N_k$ are
deterministic and can be completely separated out from $N_k$.
See further details in Remark \ref{Rem4} below and the proofs of
Theorems \ref{theo1} and
\ref{theo2} in Section~\ref{proofs}.
\end{remark}

The following theorem provides bounds on deficiency distance $\Delta
(\cP_{1n}, \cP_{2n})$
and establishes the asymptotic equivalence of trace regression and
quantum state tomography
under the fixed or random designs.
%
%th1 #&#
%
\begin{theorem}\label{theo1}
Assume that conditions \textup{(C1)--(C3)} are satisfied.
\begin{longlist}[(a)]
\item[(a)] For the random design case, we have
%
%e22 #&#
%
\begin{equation}
\label{randomPP} %\Delta(\mathbb{P}_{1,n,\brho}, \mathbb{P}_{2,n,
\Delta(\cP_{1n},
\cP_{2n}) \leq n \gamma_p + C \biggl( \frac{n
\zeta_p }{m}
\biggr)^{1/2}, %&&\qquad \leq\sqrt{n p} \max_{j=1}^p |\Pi(j) - \Xi(j) |^{1/2}
%% \left\{ |\Pi(j) - \Xi(j) |^{1/2}, 1 \leq j \leq p \right\}
% + C_\kappa(n/m)^{1/2},
\end{equation}
where $C$ is a generic constant %free of $(m,n,p,d)$,
depending only on $(\kappa, c_0, c_1)$, integer $\kappa$ and
constants $(c_0, c_1)$ are, respectively,
specified in conditions \textup{(C1)} and \textup{(C3)}, $\gamma_p$ is defined in (\ref
{Pi}), and $\zeta_p$ is given by %defined below is bounded by $1$,
%
%e23 #&#
%
\begin{equation}
\label{zeta-p} \zeta_p = \max_{\brho\in\Theta} \Biggl\{ \sum
_{j=1}^p \Pi(j) 1\bigl(\bigl|\cI_j(
\brho)\bigr|\geq2\bigr), \sum_{j=1}^p \Xi(j) 1
\bigl(\bigl|\cI_j(\brho)\bigr|\geq2\bigr) \Biggr\} \leq1.
\end{equation}
In particular, if $\Pi(j) = \Xi(j)=1/p$ for $j=1,\ldots, p$, then
%
%e24 #&#
%
\begin{equation}
\label{randomPP1} \Delta(\cP_{1n}, \cP_{2n}) \leq C \biggl(
\frac{n \zeta_p }{m} \biggr)^{1/2},
\end{equation}
where now $\zeta_p$ can be simplified as
%
%e25 #&#
%
\begin{equation}
\label{zeta-p1} \zeta_p = \max_{\brho\in\Theta} \Biggl\{
\frac{1}{p} \sum_{j=1}^p 1\bigl(\bigl|
\cI_j(\brho)\bigr|\geq2\bigr) \Biggr\} \leq1.
\end{equation}
\item[(b)] For the fixed design case, we have
%
%e26 #&#
%
\begin{equation}
\label{fixPP} \Delta(\cP_{1n}, \cP_{2n}) \leq C \biggl(
\frac{ n \zeta_p}{m} \biggr)^{1/2},
\end{equation}
where $C$ is the same as in \textup{(a)}, and $\zeta_p$ is given by (\ref{zeta-p1}).
\end{longlist}
\end{theorem}

%re2 #&#
%
\begin{remark}\label{Rem2}
Theorem \ref{theo1} establishes bounds on the deficiency
distance between trace regression and quantum state tomography. If
the deficiency distance bounds in (\ref{randomPP}), (\ref{randomPP1})
and (\ref{fixPP}) go to zero, trace regression and quantum
state tomography are asymptotically equivalent under the corresponding cases.
$\zeta_p$ defined in (\ref{zeta-p}) and (\ref{zeta-p1}) has an
intuitive interpretation as follows.
Proposition \ref{prop1} shows that each observable corresponds to a
multinomial distribution in quantum state tomography.
Of the $p$ multinomial distributions in quantum state tomography,
$\zeta_p$~is the maximum of the average fraction of
the nondegenerate multinomial distributions (i.e., with at least two
cells). As we discussed in Remark \ref{Rem1}, the multinomial
distributions have
cell probabilities $\operatorname{tr}(\bQ_{ja} \brho)$, $a=1,\ldots, r_j$. Since
for each $\bB_j$, $\operatorname{tr}(\bQ_{ja} \brho)$
is the trace of the density matrix $\brho$ restricted to the
corresponding eigen-space, and $\sum_{a=1}^{r_j} \operatorname{tr}(\bQ_{ja} \brho)
= \operatorname{tr}(\brho) =1$, thus if $|\cI_j(\brho)|\geq2$, $\brho$ cannot
live on any single eigen-space corresponding to one eigenvalue of $\bB
_j$; otherwise %under which
measurement results on $\bB_j$ are certain, and the corresponding
multinomial and normal distributions are reduced to the same
degenerate distribution and hence are always equivalent.\vadjust{\goodbreak} Therefore, to
bound the deficiency distance between quantum state
tomography and trace regression we need to consider only the
nondegenerate multinomial distributions, and thus $\zeta_p$ appears
in all the deficiency distance bounds. Since $\zeta_p$ is always
bounded by $1$, from Theorem \ref{theo1} we have that if $n/m
\rightarrow0$, the
two models are asymptotically equivalent. As we will see in Sections~\ref{section-sparse} and \ref{section-low-rank},
depending on density
matrix class $\Theta$ as well as the matrix set $\cB$, $\zeta_p$ may
or may not go to zero, and we will show that if it approaches
to zero, we may have asymptotic equivalence in ultra-high dimensions
where $d$ may be comparable to or exceed $m$.
\end{remark}

%Since $\zeta_p$ is always bounded by $1$, from Theorem \ref{theo1} we
%have that
%if $n/m \rightarrow0$,
%%further condition on the sampling distributions $\Pi(j)$ and $\Xi(j)$
%to make $n\gamma_p \rightarrow0$ for the random design case,
%the two models are asymptotically equivalent regardless of matrix size
%$d$ of $\brho$.
%%and the asymptotic equivalence holds regardless of $p$.
%In other words,
%%trace regression and quantum state tomography can be asymptotically
%equivalent regardless the matrix size of $\brho$, and
%the asymptotic equivalence of trace regression and quantum state
%tomography can be free of dimensionality of the matrix parameter.
%This is the first asymptotic equivalence result in ultra-high
%dimensions where dimension $p$ is allowed to be any number: it may
%exceed sample sizes $m$ and $n$, or even grow exponentially faster
%than the sample sizes.

%For example, consider $n$ independent
%binomial distributions $\operatorname{Bin}(m, p_{j})$, $j=1,\ldots, n$, and $n$
%independent normal distributions with the same means and variances
%as the binomial distributions. The squared deficiency distance between
%binomial distributions and the normal distributions is of order
%$\sum_{j=1}^n 1/[m p_j (1-p_j)]$, which is of order $n/m$, as
%Condition (C3) implies $p_j (1-p_j) \geq c_0 (1-c_1)>0$. %for $p_j \in
%(0,1)$.
%%A special case is that for $j \in\cI$, $p_j \sim1/r$ for some fixed
%$r$, then $p_j (1-p_j) \sim(r-1)/r^2$, and
%%$\sum_{j=1}^n 1/[m p_j (1-p_j)] \sim r^2 n/[(r-1) m]$.

%re3 #&#
%
\begin{remark}\label{Rem3}
The asymptotic equivalence results indicate that we
may apply matrix completion methods to quantum
state tomography by substituting $(\bM_k, N_k)$ from quantum state
tomography for $(\bX_k, Y_k)$ from trace regression.
For example, suppose that $\cB$ is an orthonormal basis and $\brho$
has %representation (\ref{sparse})
an expansion $\brho=\sum_j \alpha_j \bB_j$ with
$\alpha_j = \operatorname{tr}(\brho\bB_j)$. For trace regression, we may estimate
$\alpha_j$ by
the average of those $Y_k$ with corresponding $\bX_k=\bB_j$.
Replacing $(\bX_k, Y_k)$
from trace regression by $(\bM_k, N_k)$ from quantum state tomography
we construct an estimator of $\alpha_j$
by taking the average of those $N_k$ with corresponding $\bM_k=\bB
_j$. In fact, the resulting estimator based on $N_k$ can be naturally
derived from quantum state tomography. From (\ref{measurement}), (\ref
{tomography2}) and (\ref{tomography5}), we have $\alpha_j = \operatorname{tr}(\brho
\bB_j) = E(R)$,
where $R$ is the outcome of measuring $\bB_j$, and hence
it is natural to estimate $\alpha_j$ by the average of quantum
measurements $R_{k\ell}$ with corresponding $\bM_k=\bB_j$.
%See more details in Section \ref{section-Pauli} for density matrix
%estimation for Pauli measurements.
As statistical procedures and fast algorithms are available for trace
regression, these statistical methods and computational
techniques can be easily used to implement quantum state tomography
based on the summarized measurements [\citet{Groetal10} and
\citet{Kol11}].
\end{remark}

%s4 #&#
\section{Fine scale trace regression}\label{sec4}
\label{section-fine}

In Section~\ref{section-P-equivalence} for quantum state tomography we
define $\mathbb{P}_{2,n,\brho}$ and $\cP_{2 n}$ in (\ref {experiment})
based on the average measurements $N_k$, and the asymptotic equivalence
results show that trace regression matches quantum state tomography
with the summarized measurements $(\bM_k, N_k)$, $k=1,\ldots, n$. We
may use individual measurements $R_{k1},\ldots, R_{km}$ instead of
their averages $N_k$ [see (\ref{tomography0})--(\ref{tomography5}) for
their definitions and relationships], and replace
$\mathbb{P}_{2,n,\brho}$ in (\ref {experiment}) by the joint
distribution, $\mathbb{Q}_{2,n,\brho}$, of $(\bM_k, R_{k1},\ldots,
R_{km})$, $k=1,\ldots, n$, for the random design case [or
$(R_{k1},\ldots, R_{km})$, $k=1,\ldots, n$, for the fixed design case]
to define a new statistical experiment for quantum state tomography
with the individual measurements,
%
%e27 #&#
%
\begin{equation}
\label{experimentQ2} %\cQ_{1n} = \{(\cX_1, \cG_1, \mathbb{Q}_{1,n,
\cQ_{2n} = \bigl\{(
\cX_2, \cG_2, \mathbb{Q}_{2,n,\brho}), \brho\in\Theta
\bigr\},
\end{equation}
where measurable space $(\cX_2, \cG_2)$ is either $(\cB^n \times
\mathbb{R}^{m n}, \cF_\cB^n \times\cF_\mathbb{R}^{m n})$ for the
random design case
or $(\mathbb{R}^{m n}, \cF_\mathbb{R}^{m n})$ for the fixed design case.

In general, $\cP_{1n}$ and $\cQ_{2n}$ may not be asymptotically
equivalent. As individual measurements $R_{k1},\ldots, R_{km}$ may
contain more
information than their average $N_k$, $\cQ_{2n}$ may be more
informative than $\cP_{2n}$, and hence $\delta(\cQ_{2n}, \cP
_{2n})=0$ but
$\delta(\cP_{2n}, \cQ_{2n})$ may be bounded away from zero. As a
consequence, we may have $\delta(\cQ_{2n}, \cP_{1n})$ goes to zero but
$\delta(\cP_{1n}, \cQ_{2n})$ and $\Delta(\cP_{1n}, \cQ_{2n})$ are
bounded away from zero. For the special case of $\kappa=2$ where all
$\bB_j$
have at most two distinct eigenvalues such as Pauli matrices in (\ref
{Pauli-basis}), $N_k$ are sufficient statistics for the distribution of
$(R_{k1}, R_{k2})$, and
hence $\cP_{2n}$ and $\cQ_{2n}$ are equivalent, that is, $\Delta(\cP
_{2n}, \cQ_{2n})=0$, $\Delta(\cP_{1n}, \cP_{2n}) = \Delta(\cP
_{1n}, \cQ_{2n})$,
and $\cP_{1n}$ and $\cQ_{2n}$ can still be asymptotically equivalent.
In summary, generally trace regression can be asymptotically equivalent to
quantum state tomography with summarized measurements but not with
individual measurements. In fact, the individual measurements
$(R_{k1},\ldots, R_{km})$, $k=1,\ldots, n$, from quantum state
tomography contain information
about $\operatorname{tr}(\bQ_{ja} \brho)$, $a=1,\ldots, r_j$, while observations
$Y_k$, $k=1,\ldots, n$, from trace regression have information only
about $\operatorname{tr}(\bB_j \brho)$. From (\ref{basis-diagonal}) we get $\operatorname{tr}(\bB
_j \brho) = \sum_{a=1}^{r_j} \lambda_{ja} \operatorname{tr}(\bQ_{ja} \brho)$, so
the individual measurements $(R_{k1},\ldots, R_{km})$ from quantum
state tomography may be more informative
than observations $Y_k$ from trace regression for statistical inference
of $\brho$.
To match quantum state tomography with individual measurements, we may
introduce a fine scale trace regression
model and treat trace regression (\ref{trace-regression}) as a coarse
scale model aggregated from the fine scale model as follows.
Suppose that matrix input variable $\bX_k$ has the following spectral
decomposition:
%
%e28 #&#
%
\begin{equation}
\label{trace-regression2} \bX_k = \sum
_{a=1}^{r^X_k} \lambda^X_{ka}
\bQ^X_{ka},
\end{equation}
where $\lambda^X_{ka}$ are $r^X_k$ real distinct eigenvalues of $\bX
_k$, and $\bQ^X_{ka}$ are the projections onto the eigen-spaces
corresponding to $\lambda^X_{ka}$.
The fine scale trace regression model assumes that observed random
pairs $(\bQ^X_{ka}, y_{ka})$ obey
%
%e29 #&#
%
\begin{equation}
\label{trace-regression1} % V_{ka} = m \operatorname{tr}(\bQ^X_{ka} \brho) + z_{ka},
%k=1,\ldots, n,
% a=1,\ldots, r^X_k,
y_{ka} = \operatorname{tr}\bigl(
\bQ^X_{ka} \brho\bigr) + z_{ka},\qquad k=1,\ldots, n,
a=1,\ldots, r^X_k,
\end{equation}
where $z_{ka}$ are random errors with mean zero.
%$(z_{k1},\ldots, z_{k r^X_k})$ for different $k$ are independent, and
%for each $k$, $(z_{k1},\ldots, z_{k r^X_k})^\prime$ is a multivariate
%normal random vector with mean zero and
%$\operatorname{Var}(z_{ka})= m \operatorname{tr}(\bQ^X_{ka} \brho) [1-\operatorname{tr}(\bQ^X_{ka} \brho) ]$
%and $\operatorname{Cov}(z_{ka}, z_{kb})= - m \operatorname{tr}(\bQ^X_{ka} \brho) \operatorname{tr}(\bQ^X_{kb}

Models (\ref{trace-regression}) and (\ref{trace-regression1}) are
trace regression at two different scales and connected through (\ref
{trace-regression2}) and
the following aggregation relations:
%
%e30 #&#
%
\begin{equation}
\label{trace-regression0}\qquad % Y_k = \frac{1}{m} \sum_{a=1}^{r^X_k}
Y_k = \sum
_{a=1}^{r^X_k} \lambda^X_{ka}
y_{ka},\qquad \varepsilon_k = \sum_{a=1}^{r^X_k}
\lambda^X_{ka} z_{ka},\qquad \operatorname{tr}(\bX_k
\brho) = \sum_{a=1}^{r^X_k}
\lambda^X_{ka} \operatorname{tr}\bigl(\bQ^X_{ka}
\brho\bigr).
\end{equation}
The fine scale trace regression model specified by (\ref
{trace-regression1}) matches quantum state tomography with the individual
measurements $(\bM_k, R_{k1},\ldots, R_{km})$, $k=1,\ldots, n$.
Indeed, as (\ref{trace-regression2}) indicates a one to one
correspondence between $\bX_k$ and $\{\lambda^X_{ka}, \bQ^X_{ka}$,
$a=1,\ldots, r^X_k\}$, we replace $Y_k$ by
$(y_{k1},\ldots, y_{k r^X_k})$ and $\mathbb{P}_{1,n,\brho}$ in
(\ref{experiment}) by the joint distribution,
$\mathbb{Q}_{1,n,\brho}$, of $(\bX_k, y_{k1},\ldots, y_{k
r^X_k})$, $k=1,\ldots, n$, for the random design case [or
$(y_{k1},\ldots, y_{k r^X_k})$, $k=1,\ldots, n$, for the fixed
design case], and define the statistical experiment for fine
scale trace regression (\ref{trace-regression1}) as follows:
%
%e31 #&#
%
\begin{equation}
\label{experimentQ1} \cQ_{1n} = \bigl\{(\cX_1,
\cG_1, \mathbb{Q}_{1,n,\brho}), \brho\in\Theta\bigr\},
\end{equation}
where measurable space $(\cX_1, \cG_1)$ is either $(\cB^n \times
\mathbb{R}^{m n}, \cF_\cB^n \times\cF_\mathbb{R}^{m n})$ for the
random design case
or $(\mathbb{R}^{m n}, \cF_\mathbb{R}^{m n})$ for the fixed design case.

%The proof of Theorem \ref{theo1} has in fact derived bounds on $\Delta(
%trace regression model
%(\ref{trace-regression1}) and quantum state tomography with individual
%measurements $(\bM_k, R_{k1},\ldots, R_{km})$, $k=1,\ldots, n$.
%has two scales but trace regression has only one scale. ?? Trace
%regression does not have any scale match with raw data $R_{jm}$?
To study the asymptotic equivalence of fine scale trace regression and
quantum state tomography with individual measurements, we need to
replace condition (C2)
by a new condition for fine scale trace regression:
\begin{longlist}[(C2$^*$)]
\item[(C2$^*$)] Suppose that two models $\cQ_{1n}$ and $\cQ_{2n}$ are
identifiable.
For fine scale trace regression (\ref{trace-regression1}), random
errors $(z_{k1},\ldots, z_{k r^X_k})$, $k=1,\ldots, n$, are
independent, and %for each $k$,
given $\bX_k$, $(z_{k1},\ldots, z_{k r^X_k})^\prime$ is a
multivariate normal random vector with mean zero and for $a, b =
1,\ldots, r_k^X$, $a \neq b$,
%
%e32 #&#
%
\begin{eqnarray}
\label{C2*-var}
\operatorname{Var}(z_{ka}|\bX_k)&=& \frac{1}{m} \operatorname{tr}
\bigl(\bQ^X_{ka} \brho\bigr) \bigl[1-\operatorname{tr}\bigl(\bQ
^X_{ka} \brho\bigr) \bigr],\nonumber\\[-8pt]\\[-8pt]
\operatorname{Cov}(z_{ka},
z_{kb}|\bX_k)&=& - \frac{1}{m} \operatorname{tr}\bigl(
\bQ^X_{ka} \brho\bigr) \operatorname{tr}\bigl(\bQ^X_{kb}
\brho\bigr).\nonumber
\end{eqnarray}
\end{longlist}
We provide bounds on $\Delta(\cQ_{1n}, \cQ_{2n})$ and establish the
asymptotic equivalence of $\cQ_{1n}$ and $\cQ_{2n}$ in the following theorem.

%
%th2 #&#
%
\begin{theorem}\label{theo2}
Assume that conditions \textup{(C1)}, \textup{(C2$^*$)} and
\textup{(C3)} are satisfied.
\begin{longlist}[(b)]
\item[(a)] For the random design case, we have
%
%e33 #&#
%
\begin{equation}
\label{randomQQ} \Delta(\cQ_{1n}, \cQ_{2n}) \leq n
\gamma_p + %\kappa
C \biggl( \frac{n \zeta_p }{m}
\biggr)^{1/2},
\end{equation}
where as in Theorem \ref{theo1}, $C$ is a generic constant depending
only on
$(\kappa, c_0, c_1)$, integer $\kappa$ and
constants $(c_0, c_1)$ are, respectively, specified in conditions
\textup{(C1)}
and~\textup{(C3)}, and $\gamma_p$ and $\zeta_p$ are given
by (\ref{Pi}) and (\ref{zeta-p}), respectively. In particular, if
$\Pi(j) = \Xi(j)=1/p$ for $j=1,\ldots, p$, then
%
%e34 #&#
%
\begin{equation}
\label{randomQQ1} \Delta(\cQ_{1n}, \cQ_{2n}) \leq C \biggl(
\frac{n \zeta_p }{m} \biggr)^{1/2},
\end{equation}
where $\zeta_p$ is given by (\ref{zeta-p1}).
\item[(b)] For the fixed design case, we have
%
%e35 #&#
%
\begin{equation}
\label{fixQQ} \Delta(\cQ_{1n}, \cQ_{2n}) \leq C \biggl(
\frac{ n \zeta_p}{m} \biggr)^{1/2},
\end{equation}
where $C$ is the same as in \textup{(a)}, and $\zeta_p$ is given by (\ref{zeta-p1}).
\end{longlist}
\end{theorem}

%re4 #&#
%
\begin{remark}\label{Rem4}
For quantum state tomography we regard summarized
measurements and individual measurements
as quantum measurements at coarse and fine scales, respectively. Then
Theorems \ref{theo1} and \ref{theo2} show that quantum state
tomography and trace regression are asymptotically equivalent at both
coarse and fine scales. Moreover, as
measurements at the coarse scale are aggregated from measurements at
the fine scale for both quantum state
tomography and trace regression, their asymptotic equivalence at the
coarse scale is a consequence of their
asymptotic equivalence at the fine scale. Specifically, the deficiency
distance bounds in (\ref{randomQQ})--(\ref{fixQQ}) of Theorem
\ref{theo2}
are derived essentially from the deficiency distance between $n$
independent multinomial distributions in quantum state tomography and
their corresponding multivariate normal distributions in fine scale
trace regression,
and the deficiency distance bounds in (\ref{randomPP}), (\ref
{randomPP1}) and (\ref{fixPP}) of Theorem \ref{theo1}
are the consequences of corresponding bounds in Theorem \ref{theo2}.
%As we have discussed in Remark \ref{Rem1} and may see from the proofs
%of
%Theorems \ref{theo1} and \ref{theo2} in Section \ref{proofs},
%(\ref{tomography0}) and (\ref{trace-regression0}) show that $N_k$ and $
%$(R_{k1},\ldots, R_{km})$ and $m (z_{k1},\ldots, z_{k r^X_k})$,
%respectively.
Fine scale trace regression (\ref{trace-regression1}) and condition
(C2$^*$) indicate that for each $k$,
$(y_{k1},\ldots, y_{k r^X_k})$ follows a multivariate normal
distribution. From (\ref{measurement-U}) and
(\ref{tomography1})--(\ref{tomography5})
we see that given $\bM_k$, $(R_{k1},\ldots, R_{km})$ is jointly
determined by the counts of $R_{k1},\ldots, R_{km}$
taking the eigenvalues of $\bM_k$, and the counts jointly follow a
multinomial distribution, with mean and covariance matching
with those of $m (y_{k1},\ldots, y_{k r^X_k})$. To prove Theorems \ref{theo1}
and \ref{theo2}, we need to derive
the Hellinger distances of the multivariate normal distributions and
their corresponding multinomial distributions with uniform perturbations.
%and then establish bounds on the deficiency distance between quantum
%state tomography and trace regression.
\citet{Car02} has established a bound on deficiency distance between a
multinomial distribution and its corresponding
multivariate normal distribution through the total variation distance
between the multivariate normal distribution and
the multinomial distribution with uniform perturbation. The main
purpose of the multinomial deficiency bound in
\citet{Car02} is the asymptotic equivalence study for density
estimation. Consequently, the multinomial distribution in
\citet{Car02} is allowed to have a large number of cells, with bounded
cell probability ratios, and his proof techniques
are geared up for managing such a multinomial distribution under total
variation distance.
Since quantum state tomography involves many independent multinomial
distributions all with a small number of cells,
Carter's result is not directly applicable for proving Theorems \ref
{theo1} and
\ref{theo2}, nor his approach suitable for the current model setting.
%Carter's result cannot be directly applied to prove Theorems
%key proof techniques like conditional pairing and associated Markov
%kernel
%mapping arguments are inappropriate or unsuitable for our setting.
%derive the needed Hellinger distance results.
To show Theorems \ref{theo1} and \ref{theo2}, we deal with $n$
independent multinomial
distributions in quantum state tomography by
deriving the Hellinger distances between the perturbed multinomial
distributions and the corresponding multivariate normal distributions,
and then we establish bounds on the deficiency distance between quantum
state tomography and trace regression at the fine scale. Moreover, from
(\ref{measurement-RU}), (\ref{tomography0}) and (\ref
{trace-regression0}) we derive $N_k$ from the counts of individual measurements
$R_{k1},\ldots, R_{km}$ for quantum state tomography and $Y_k$ from
fine scale observations $y_{ka}$ for trace regression by the
same aggregation relationship, and %condition (C2$^*$) implies condition C2,
(\ref{C2*-var}) implies (\ref{C2-var}),
so bounds on $\Delta(\cP_{1n}, \cP_{2n})$ can be obtained from those
on $\Delta(\cQ_{1n}, \cQ_{2n})$. Thus, Theorem~\ref{theo1} may be
viewed as a
consequence
of Theorem~\ref{theo2}. For more details see the proofs of Theorems \ref
{theo1} and \ref{theo2} in
Section~\ref{proofs}.
\end{remark}

%s5 #&#
\section{Sparse density matrices}\label{sec5}
\label{section-sparse}
Since all deficiency distance bounds in Theorems \ref{theo1} and \ref
{theo2} depend on
$\zeta_p$, we further investigate $\zeta_p$ for
two special classes of density matrices: sparse density matrices in
this section and low rank density matrices in Section~\ref{section-low-rank}.

%
%co1 #&#
%
\begin{cor} \label{cor-sparse}\label{cor1}
Denote by $\Theta_s$ a collection of density matrices with at most $s$
nonzero entries,
where $s$ is an integer.
%deterministic function of $p$ that grows very slowly in $p$ with an
%example of $s = c$ and $s = c \log p$ for some constant $c$.
Assume that $\cB$ is selected as basis (\ref{Hermitian-basis}), and
$\Pi(j)=\Xi(j)=1/p$. Then
\[
\zeta_p = \max_{\brho\in\Theta_s} \Biggl\{ \frac{1}{p}
\sum_{j=1}^p 1\bigl(\bigl|\cI_j(
\brho)\bigr|\geq2\bigr) \Biggr\} \leq\frac{s_d}{d},
\]
where $s_d$ is the maximum number of nonzero diagonal entries of $\brho
$ over $\Theta_s$.
Furthermore, if conditions \textup{(C1), (C2), (C2$^*$)} and \textup{(C3)} are satisfied, we have
\[
\Delta(\cP_{1n}, \cP_{2n}) \leq C \biggl( \frac{ n s_d}{m d}
\biggr)^{1/2},\qquad \Delta(\cQ_{1n}, \cQ_{2n}) \leq C
\biggl( \frac{ n s_d}{m d} \biggr)^{1/2},
\]
where $C$ is the same generic constant as in Theorems \ref{theo1} and
\ref{theo2}.
\end{cor}

%re5 #&#
%
\begin{remark}\label{Rem5}
Since $p=d^2$, $s_d \leq s$, and the deficiency
distance bounds in Corollary \ref{cor-sparse} are of order
$[n s_d/(m d)]^{1/2}$, if $s_d/d$ goes to zero as $d \rightarrow
\infty$, we may have that as $m,n,d \rightarrow\infty$,
$n s_d/(m d) \rightarrow0$ and hence the asymptotic equivalence of
quantum state tomography and trace regression, while
$n/m$ may not necessarily go to zero.
%although $n/m$ does not go to zero. %is very large or even goes to
%infinity.
Thus, even though sparsity is not required in the asymptotic
equivalence of quantum state tomography and trace regression, Corollary
\ref{cor-sparse} shows that with the sparsity the asymptotic
equivalence is much
easier to achieve. For example, consider the case that $s_d$ is
bounded, and $n$ is of order $d^2$
(suggested by the bounded $\kappa$ and the identifiability discussion
at the end of Section~\ref{sec2.3}).
In this case the deficiency distance bounds in Corollary \ref
{cor-sparse} are of order $(d/m)^{1/2}$, and we obtain the asymptotic
equivalence
of quantum state tomography and trace regression, if $d/m \rightarrow
0$ with an example $d=O(m/\log m)$.
\end{remark}

We illustrate below that the sparse density matrices studied in
Corollary~\ref{cor-sparse} have a sparse representation under basis
(\ref{Hermitian-basis}).
In general, assume that $\cB$ is an orthogonal basis for complex
Hermitian matrices. %(\ref{Hermitian-basis}).
Then every density matrix $\brho$ has a representation under the basis
$\cB$,
%
%e36 #&#
%
\begin{equation}
\label{sparse} \brho= \sum_{j=1}^p
\alpha_j \bB_j,
\end{equation}
where $\alpha_j$ are coefficients. We say a density matrix $\brho$ is
$s$-sparse under the basis~$\cB$, if
%$\cB= \{\bB_1,\ldots, \bB_p\}$ is Pauli basis (\ref{Pauli-basis}),
the representation (\ref{sparse}) of $\brho$ under the basis $\cB$
has at most $s$ nonzero coefficients~$\alpha_j$.
The sparsity definition via representation (\ref{sparse}) is in line
with the vector sparsity concept through orthogonal expansion in
compressed sensing.
It is easy to see that a density matrix $\brho$ with at most $s$
nonzero entries is the same as that $\brho$ is $s$-sparse
%of nonzero coefficients $\alpha_j$ in the representation (
under basis (\ref{Hermitian-basis}).
% $s$ is a deterministic function of $p$ that grows very slowly in $p$
%with an example of $s = C$ and $s = C \log p$ for some constant $C$.
%For example, if $\Theta_s$ is a class of density matrices with at most
%$s$ nonzero elements,
%and we take $\cB$ as basis (\ref{Hermitian-basis}), then any $\brho$
%in $\Theta_s$ has an representation under
%basis (\ref{Hermitian-basis}) with at most $s$ nonzero coefficients.
However, a $s$-sparse matrix under the Pauli basis (\ref{Pauli-basis})
may have more than $s$ nonzero entries. In fact,
%while a $s$-sparse density matrix under basis (\ref{Hermitian-basis})
%has at most $s$ nonzero entries,
it %a $s$-sparse density matrix under Pauli basis (\ref{Pauli-basis})
may have up to $s d$ nonzero entries. The following corollary exhibits
the different behavior of $\zeta_p$ for
sparse density matrices under the Pauli basis.

%
%co2 #&#
%
\begin{cor}\label{cor-sparse-Pauli}\label{cor2}
Denote by $\Theta^p_{s}$ the class of all density matrices that are
$s$-sparse under the Pauli basis,
where $s$ is an integer.
% that in expansion (\ref{sparse}), of those $j \neq1$ with nonzero
%coefficients $\alpha_j$, their corresponding $\bB_j = \bsigma_{\ell_1}
%contains no more than $\gamma$ number of $\bsigma_0$, where $\gamma$
%is a fixed integer free of $p$.
Assume that $\cB$ is selected as the Pauli basis (\ref{Pauli-basis}),
and $\Pi(j)=\Xi(j)=1/p$. Then
\[
1 \geq\zeta_p = \max_{\brho\in\Theta^p_{s}} \Biggl\{
\frac{1}{p} \sum_{j=1}^p 1\bigl(\bigl|
\cI_j(\brho)\bigr|\geq2\bigr) \Biggr\} \geq1 - \frac
{1}{p}.
\]
Furthermore, if conditions \textup{(C1), (C2), (C2$^*$)} and
\textup{(C3)} are satisfied, we have
\[
\Delta(\cP_{1n}, \cP_{2n}) \leq C \biggl( \frac{ n}{m}
\biggr)^{1/2},\qquad \Delta(\cQ_{1n}, \cQ_{2n}) \leq C
\biggl( \frac{ n}{m} \biggr)^{1/2},
\]
where $C$ is the same generic constant as in Theorems \ref{theo1} and
\ref{theo2}.
\end{cor}
%
%re6 #&#
%
\begin{remark}\label{Rem6}
Corollary \ref{cor1} shows that for sparse matrices under
basis (\ref{Hermitian-basis}), %and low rank (perhaps dense) matrices
%we may have that
as $d \rightarrow\infty$, if $s_p/d \rightarrow0$, $\zeta_p$ goes
to zero, and hence the sparsity enables us to establish the asymptotic
equivalence of quantum state
tomography and trace regression under weaker conditions on $m$ and $n$.
However, Corollary \ref{cor2} demonstrates that $\zeta_p$ does
not go to zero for sparse matrices under the Pauli basis. Corollary \ref{cor1}
indicates that for a density matrix with $s$ nonzero entries, in order
to have small $s_p/d$, we must make its
nonzero diagonal entries as less as possible. The Pauli basis is the
worst in a sense that a sparse matrix under the Pauli basis has at
least $d$
nonzero entries, and the Pauli basis tends to put many nonzero entries
on the diagonal. From Corollaries \ref{cor1} and \ref{cor2} we see that
$\zeta_p$ depends on sparsity of the density matrix class, but more
importantly it is determined by how the sparsity is specified by $\cB$.
%defines observables and matrix variables in the two models as well as
\end{remark}

%s6 #&#
\section{Low rank density matrices}\label{sec6}
\label{section-low-rank}

Consider the case of low rank density matrices.
%Let $\Theta_r$ consist of density matrices with rank at most $r$ and
%satisfying (\ref{theta}).
Assume density matrix $\brho$ has rank at most $r$, where $r \ll d$.
%%$r$ is far smaller than $d$.
Then $\brho$ has at most $r$ nonzero eigenvalues, and thus its
positive eigenvalues are sparse.
%Matrix completion is to recover the unknown low rank matrix of $p=d^2$
%entries based on $n$ observations.
%%so the sample size could be required of order $p$.
%It has been shown that for both noisy and noiseless matrix completion
%problems, to recover an unknown matrix of rank at
%most $r$ we usually need only $n_*$ observations, where $n_*$ is of
%order $d r \log d $ instead of potential order $d^2$ (see Cand\'es
%and Plan (2009a, b),
%Cand\'es and Recht (2008), \citet{Gro11}, Koltchinskii, Lounici and
%Tsybakovet (2011), \citet{Liu}, and \citet{RohTsy11}). The
%matrix completion result and Theorems \ref{theo1} and \ref{theo2}
%together indicate that
%if $n \geq n_*$ and $n \zeta_p /m \rightarrow0$, trace regression and
%quantum state tomography are asymptotic equivalent,
%and quantum state tomography can reconstruct the unknown density
%matrix of rank up to $r$. %We study $\zeta_p$ for low rank density
%matrices.
%The following corollary indicates that for special low rank matrices,
%we need very weak conditions on $(m,n)$ to achieve asymptotic
%equivalence.
The following corollary derives the behavior of $\zeta_p$ for low rank
density matrices and the Pauli basis.

%
%co3 #&#
%
\begin{cor} \label{cor-lowrank-Pauli}\label{cor3}
Denote by $\Theta_{r}$ the collection of all density matrices $\brho$
with rank up to $r \ll d$.
Assume that\vadjust{\goodbreak} $\cB$ is the Pauli basis (\ref{Pauli-basis}), and $\Pi
(j)=\Xi(j)=1/p$. Then
\[
1 \geq\zeta_p = \max_{\brho\in\Theta_r} \Biggl\{
\frac{1}{p} \sum_{j=1}^p 1\bigl(\bigl|
\cI_j(\brho)\bigr|\geq2\bigr) \Biggr\} \geq1 - \frac
{1}{p}.
\]
Furthermore, if conditions \textup{(C1), (C2), (C2$^*$)} and
\textup{(C3)} are satisfied, we have
\[
\Delta(\cP_{1n}, \cP_{2n}) \leq C \biggl( \frac{ n}{m}
\biggr)^{1/2},\qquad \Delta(\cQ_{1n}, \cQ_{2n}) \leq C
\biggl( \frac{ n}{m} \biggr)^{1/2},
\]
where $C$ is the same generic constant as in Theorems \ref{theo1} and
\ref{theo2}.
\end{cor}
We construct a low rank density matrix class and matrix set for which
$\zeta_p$ goes to zero in the following corollary.
%
%co4 #&#
%
\begin{cor} \label{cor-lowrank}\label{cor4}
Suppose that $\bg_1,\ldots, \bg_d$ form an orthonormal basis in
$\mathbb{R}^d$, and
%$\cB$ is the matrix set we build from $\bg_1,\ldots, \bg_d$ the same
%way as we build basis (\ref{Hermitian-basis}) from $\be_1,\ldots,
%
\begin{eqnarray*}
\cB&=& \biggl\{ \bg_{\ell} \bg_{\ell}^\prime,
\frac{1}{\sqrt{2}} \bigl( \bg_{\ell_1} \bg_{\ell_2}^\prime+
\bg_{\ell_2} \bg_{\ell
_1}^\prime\bigr), \frac{\sqrt{-1}}{\sqrt{2}}
\bigl( \bg_{\ell_2} \bg_{\ell_1}^\prime- \bg_{\ell_1}
\bg_{\ell
_2}^\prime\bigr),\\
&&\hspace*{111.6pt} \ell, \ell_1,
\ell_2=1,\ldots, d, \ell_1 <\ell_2 \biggr
\}.
\end{eqnarray*}
Assume that $\gamma\ll d$ and $r \ll d$ are integers. Denote by
$\Theta_{r\gamma}$ a collection of density matrices $\brho$ with the form
%If $\brho$ is rank $r$, we have representation
%
%e37 #&#
%
\begin{equation}
\label{rank-r} \brho= \sum_{j=1}^r
\xi_j U_j U_j^\dagger,
\end{equation}
where $\xi_j \geq0$, $\xi_1+\cdots+\xi_r=1$, and $U_j$ are unit
vectors in $\mathbb{C}^{d}$ whose real and imaginary parts
are linear combinations of $\bg_{\ell_1},\ldots, \bg_{\ell_k}$,
$1\leq\ell_1,\ldots, \ell_k \leq d$ and $1 \leq k \leq\gamma$.
%at most $\gamma$ basis vectors selected from $\bg_1,\ldots, \bg_d$,
Assume $\Pi(j)=\Xi(j)=1/p$. Then
\[
\zeta_p = \max_{\brho\in\Theta_{r\gamma}} \Biggl\{ \frac{1}{p}
\sum_{j=1}^p 1\bigl(\bigl|\cI_j(
\brho)\bigr|\geq2\bigr) \Biggr\} \leq\frac{2 r
\gamma(4 \gamma+1)}{p}.
\]
Furthermore, if conditions \textup{(C1), (C2), (C2$^*$)} and
\textup{(C3)} are satisfied, we have
\[
\Delta(\cP_{1n}, \cP_{2n}) \leq C \biggl( \frac{ n r \gamma
^2}{m p}
\biggr)^{1/2},\qquad \Delta(\cQ_{1n}, \cQ_{2n}) \leq C
\biggl( \frac{ n r \gamma
^2}{m p} \biggr)^{1/2},
\]
where $C$ is the same generic constant as in Theorems \ref{theo1} and
\ref{theo2}.
\end{cor}
%
%re7 #&#
%
\begin{remark}\label{Rem7}
It is known that a density matrix of rank up to $r$
has representation (\ref{rank-r}), and matrix $\brho$ with
representation (\ref{rank-r}) has rank at most $r$. Corollary \ref
{cor-lowrank-Pauli} shows that for the class of density
matrices with rank at most $r$, $\zeta_p$ does not go to zero under
the Pauli basis.
% and the deficiency distance bounds are of order $(n/m)^{1/2}$.
%To make $\zeta_p$ go to zero, we need to take a special subclass of
%low rank density matrices.
Corollary \ref{cor-lowrank} constructs a basis $\cB$ and a subclass
of low rank density matrices, for which $\zeta_p$ can
go to zero, and the deficiency distance bounds are of order $[n r
\gamma^2/(m p)]^{1/2}$. Since $r, \gamma\ll d$ and
$p=d^2$, $r \gamma^2/ p$ may go to zero very fast as $d \rightarrow
\infty$.
%we may have that for fixed $m$,
As $m,n,d \rightarrow\infty$, if $n r \gamma^2/(m p) \rightarrow
0$, we obtain the asymptotic equivalence of quantum state tomography
and trace regression. %, in spite of small $m$.
%Low rank is not required to establish the asymptotic equivalence of
%quantum state tomography and trace regression, but Corollary
%the asymptotic equivalence under the low rank density matrix case
%needs much weaker conditions on $(m,n)$.
For example, consider the case that $r$ and $\gamma$ are bounded, and
$n$ is of order $d^2$
(suggested by the bounded $\kappa$ and the identifiability discussion
at the end of Section~\ref{sec2.3}).
In this case the deficiency distance bounds in Corollary \ref
{cor-lowrank} are of order $m^{-1/2}$, and we conclude that if $m
\rightarrow\infty$, the two models are asymptotically equivalent for
any $(n,d)$ compatible with the model identifiability
condition. A particular example is that $n=d^2$ and $d$ grows
exponentially faster than $m$.
%regardless of matrix size $d$.
\end{remark}

%re8 #&#
%
\begin{remark}\label{Rem8}
The low rank condition $r \ll d$ on a density
matrix indicates that it has a relatively small number of positive
eigenvalues, that is, its positive eigenvalues are sparse. We may also
explain the condition on the eigenvectors $U_j$ in (\ref{rank-r}) via
sparsity as follows. Since $\{\bg_1,\ldots, \bg_d\}$ is an
orthonormal basis in $\mathbb{R}^d$, the real part, $\operatorname{Re}(U_j)$, and
imaginary part, $\operatorname{Im}(U_j)$, of $U_j$ have the following expansions under
the basis:
%
%e38 #&#
%
\begin{equation}
\label{rank-r-1} \operatorname{Re}(U_j) = \sum_{\ell=1}^d
\alpha^j_{1\ell} \bg_\ell,\qquad \operatorname{Im}(U_j) =
\sum_{\ell=1}^d \alpha^j_{2\ell}
\bg_\ell,
\end{equation}
where $\alpha^j_{1\ell}$ and $\alpha^j_{2 \ell}$ are coefficients. Then
a low rank density matrix with representation (\ref{rank-r}) belongs to
$\Theta_{r\gamma}$, if for $j=1,\ldots, r$, $\{\ell, \alpha^j_{1\ell}
\neq0\}$ and $\{\ell, \alpha^j_{2\ell} \neq0\}$ have cardinality at
most $\gamma $, that is, there are at most $\gamma$ nonzero
coefficients in the expansions (\ref{rank-r-1}). As $\gamma\ll d$, the
eigenvectors $U_j$ have sparse representations. Thus, the subclass
$\Theta_{r\gamma}$ of density matrices imposes some sparsity conditions
on not only the eigenvalues but also the eigenvectors of its members.
In fact, \citet{WitTibHas09} indicates that we need some sparsity
on both eigenvalues and eigenvectors for estimating large matrices. An
important class of quantum states are pure states, which correspond to
density \mbox{matrices} of rank one. In order to have a pure state in
$\Theta_{r\gamma}$, its eigenvector $U_1$ corresponding to
eigenvalue~$1$ must be a liner combination of at most $\gamma$ basis
vectors $\bg _\ell$. Such a requirement can be met for a large class of
pure states through the selection of proper $\gamma$ and suitable bases
in $\mathbb{R}^d$. It is interesting to see that matrices themselves in
$\Theta_{r\gamma
}$ of Corollary \ref{cor4} may not be sparse. %in fact they can be very
%dense.
For example, taking $\bg_1,\ldots, \bg_d$ as the Haar basis in
$\mathbb{R}^d$ [see \citet{Vid99}],
%that $\bg_1=(1,1,\ldots,1)^\prime$ and $\bg_2=(1,\ldots, 1, -1,
%,\ldots, -1)^\prime$. It is easy to see
we obtain that rank one matrix
$\brho= (1,1,\ldots,1)^\prime(1,1,\ldots,1)/d$ and rank two matrix
$\brho= 3 (1,1,\ldots,1)^\prime(1,1,\ldots,1)/(4 d) +
(1,\ldots, 1, -1,\ldots, -1)^\prime(1,\break\ldots, 1, -1,\ldots,
-1)/(4 d)$, which are inside $\Theta_{r\gamma}$ for $(r,\gamma
)=(1,1)$ and $(r,\gamma)=(2,2)$, respectively, but not sparse.
\end{remark}

%re9 #&#
%
\begin{remark}\label{Rem9}
From Corollaries \ref{cor1}--\ref{cor4}, we see that whether $\zeta_p$
goes to zero or
not is largely dictated by $\cB$ used in the two models.
As we discussed in Remarks \ref{Rem5} and \ref{Rem7}, for certain
classes of sparse or
low rank density matrices, $\zeta_p$ goes to zero,
and we can achieve the asymptotic equivalence of quantum state
tomography and trace regression when $d$ is comparable to or
exceeds $m$. In particular for a special subclass of low rank density
matrices we can obtain the asymptotic equivalence even
when $d$ grows exponentially faster than $m$. We should emphasize that
the claimed asymptotic equivalences
in the ultra high dimension setting are under some sparse circumstances
for which $\zeta_p$ goes to zero, that is, of the
$p$ multinomial distributions in the quantum state tomography model, a
relatively small number of multinomial distributions
are nondegenerate, and similarly, the trace regression model as the
approximating normal experiment consists of the same
small number of corresponding nondegenerate normal distributions. In
other words, the asymptotic equivalence in ultra high
dimensions may be interpreted as the approximation of a sparse quantum
state tomography model by a sparse Gaussian
trace regression model. This is the first asymptotic equivalence result
in ultra high dimensions. It leads us to speculate
that sparse Gaussian experiments may play an important role in the
study of asymptotic equivalence in the ultra high dimension setting.
%for ultra high dimensional statistical models.
\end{remark}

%s7 #&#
\section{Proofs}\label{sec7}
\label{proofs}

%s7.1 #&#
\subsection{Basic facts and technical lemmas}\label{sec7.1}
We need some basic results about the Markov kernel method which are
often used to bound $\delta(\cP_{2n}, \cP_{1n})$ and prove
asymptotic equivalence of $\cP_{1n}$ and $\cP_{2n}$
%and total variation to prove Theorem \ref{theo1} %and the
%Kullback-Leibler
%divergence for handling the deficiency
[see \citet{LeC86} and \citet{LeCYan00}].
%The Markov kernel method is often used to bound $\delta(\cP_{2n},
%One effective approach is the Markov kernel method.
A Markov kernel $K(\omega, A)$ is defined for $\omega\in\cX_2$ and
$A \in\cG_1$
%$\omega\in\cB^n \times\mathbb{R}^n$ and $A \in\cF_\cB^n \times\cF_
such that for a given $\omega\in\cX_2$, $K(\omega, \cdot)$ is a
probability measure on the $\sigma$-field $\cG_1$, and for a fixed $A
\in\cG_1$, $K(\cdot, A)$ is a measurable function on $\cX_2$. The
Markov kernel maps any $\mathbb{P}_{2,n,\brho} \in\cP_{2n}$ into
another probability measure $[K (\mathbb{P}_{2,n,\brho})](A) = \int
K(\omega, A) \mathbb {P}_{2,n,\brho}(d \omega) \in\cP_{1n}$. We have
the following result:
%
%e39 #&#
%
\begin{equation}
\label{kernel} \delta(\cP_{2n}, \cP_{1n}) \leq\inf
_{K} \sup_{\brho\in\Theta} \bigl\| \mathbb{P}_{1,n,\brho}
- K (\mathbb{P}_{2,n,\brho}) \bigr\|_{\mathrm{TV}},
\end{equation}
where the infimum is over all Markov kernels, and $\|\cdot\|_{\mathrm{TV}}$ is
the total variation norm.

We often use the Hellinger distance %$H(\cdot, \cdot)$ and the
%Kullback-Leibler divergence $D(\cdot, \cdot)$
to bound total variation norm and handle product probability measures.
For two probability measures $P$ and $Q$ on a common measurable space,
we define
the Hellinger distance
%Suppose that random variables $(U_1, U_2)$ and $(V_1, V_2)$ have joint
%distributions $P$ and $Q$, respectively. Then
%
%e40 #&#
%
\begin{equation}
\label{Hellinger0} H^2(P, Q) = \int\biggl\llvert\sqrt{
\frac{dP}{d\mu}} - \sqrt{\frac
{dQ}{d\mu}} \biggr\rrvert^2 \,d\mu,
\end{equation}
where $\mu$ is any measure that dominates $P$ and $Q$, and if $P$ and
$Q$ are equivalent,
%
%e41 #&#
%
\begin{equation}
\label{Hellinger1} H^2(P, Q) = 2 - 2 E_P \biggl[\sqrt{
\frac{dQ}{dP}} \biggr], %
%D(P, Q) = E_P \left[ \log\frac{dP}{dQ} \right],
\end{equation}
where $E_P$ denotes expectation under $P$. We have
%
%e42 #&#
%
\begin{equation}
\label{TV-Hellinger-KL} \| P - Q \|_{\mathrm{TV}} \leq H(P, Q), %\leq D(P, Q),
\end{equation}
%
%Suppose that probability measures $P$ and $Q$ are on the same
%measurable space.
and for any event $A$,
%&& \| P - Q \|_{\mathrm{TV}}^2 \leq H^2(P, Q) %\nonumber\\ &&
%&&\qquad \leq2 \int_{A^c} \left[\frac{dP}{d\nu} + \frac{dQ}{d\nu} \right] \,d
%&&\qquad = 2 [P(A^c) + Q(A^c)] + \int_{A} \left|\sqrt{\frac{dP}{d\nu}} -
%We may bound the Hellinger distance by the likelihood ratio $dP/dQ$
%used in the Kullback-Leibler divergence,
%
%e43 #&#
%
\begin{eqnarray}
\label{TV-KL1} % \| P - Q \|_{\mathrm{TV}}^2 \leq
H^2(P, Q) &\leq&2 - 2 E_P
\biggl[ 1_A \sqrt{\frac{dQ}{dP}} \biggr] = 2 P
\bigl(A^c\bigr) + 2 E_P \biggl[ 1_A
\biggl(1-\sqrt{\frac{dQ}{dP}} \biggr) \biggr]
\nonumber\\[-8pt]\\[-8pt]
&\leq&2 P\bigl(A^c\bigr) + E_P \biggl[
1_A \log\frac{dP}{dQ} \biggr],\nonumber
\end{eqnarray}
where the last inequality is from the fact that $x - 1 \geq\log x$ for
any $x>0$.
%In particular, if $P = \prod_{k=1}^n P_k$ and $Q = \prod_{k=1}^n Q_k$,
% H^2(P, Q) \leq\sum_{k=1}^n H^2(P_k, Q_k) \leq2
% \sum_{k=1}^n P(A^c_k) + \sum_{k=1}^n E_P \left[ 1_{A_k} \log

\citet{Car02} has established an asymptotic equivalence of a
multinomial distribution and its corresponding multivariate normal
distribution through bounding the total variation distance between the
multivariate normal distribution and the multinomial distribution
with uniform perturbation. The approach in \citet{Car02} is to break
dependence in the multinomial
distribution and create independence by successively conditioning on
pairs and thus establish a bound on the total variation distance of
the perturbed multinomial distribution and the multivariate normal distribution.
\citet{Car02} works for the multinomial distribution with a large
number of cells,
while quantum state tomography involves many independent multinomial
distributions all with a small number of cells.
To handle the many small independent multinomial distributions for
quantum state tomography and prove Theorems \ref{theo1} and \ref
{theo2}, we need
to derive the Hellinger distances between the perturbed multinomial
distributions and
multivariate normal distributions instead of total variation distance.
Carter's approach is geared up for total variation distance and the
result cannot be directly used to prove Theorems \ref{theo1} and \ref
{theo2}. %derive the
%needed Hellinger distance results.
Our approach to proving Lemma \ref{lemMulti-norm} below is to directly
decompose a multinomial distribution as products of conditional distributions
and then establish a bound on the Hellinger distance between the
perturbed multinomial distribution and its corresponding multivariate
normal distribution.

Denote by $C$ a generic constant whose value may change from appearance
to appearance. The value of $C$ may depends on
fixed constants $(\kappa, c_0, c_1)$ given by conditions (C1) and (C3) but
is free of $(m,n,d,p)$ and individual~$\brho$.

First, we describe a known result between binomial and normal
distributions %whose proof can be found in
[see Carter (\citeyear{Car02}), B2 of the Appendix].
%
%le1 #&#
%
\begin{lem} \label{lemBin-norm}
Suppose that $P$ is a binomial distribution $\operatorname{Bin}(m, \theta)$ with
$\theta\in(0,1)$, and $Q$ is a normal distribution with mean $m
\theta$ and variance $m \theta(1-\theta)$.
Let $P^{*}$ be the convolution distribution of $P$ and an independent
uniform distribution on $(-1/2, 1/2)$. Then
\[
P^*\bigl(A^c\bigr) \leq\exp\bigl(- C m^{1/3}\bigr),\qquad
E_{P^*} \biggl[ 1(A) \log\frac
{dP^*}{dQ} \biggr] \leq
\frac{C}{m \theta(1-\theta)},
\]
where %$C$ is a generic constant,
$A = \{ | U - m \theta| \leq m [\theta(1-\theta)]^{2/3}
\}$, and random variable $U$ has the distribution $P$.
\end{lem}

We give bounds on the Hellinger distances between the perturbed
multinomial distributions and their corresponding multivariate normal
distributions in next two lemmas whose proofs are collected in the
\hyperref[app]{Appendix}.

%
%le2 #&#
%
\begin{lem} \label{lemMulti-norm}
Suppose that $P$ is a multinomial distribution ${\cal M}(m, \theta
_1,\ldots, \theta_{r})$, where $r \geq2$ is a fixed integer,
\[
\theta_1 + \cdots+ \theta_{r}=1,\qquad c_0 \leq
\min(\theta_1,\ldots, \theta_{r}) \leq\max(
\theta_1,\ldots, \theta_{r}) \leq c_1
\]
and $0< c_0 \leq c_1<1$ are two fixed constants. Denote by $Q$ the
multivariate normal distribution whose mean and covariance
are the same as $P$. Let $P^{*}$ be the convolution of the distribution
$P$ and the distribution of $(\psi_1,\ldots, \psi_r)$, where
$\psi_1,\ldots, \psi_{r-1}$ are independent and follow a uniform
distribution on $(-1/2, 1/2)$, and $\psi_r = - \psi_1 - \cdots- \psi
_{r-1}$. Then
\[
%H(P^*, Q) \leq\frac{C}{\sqrt{m} }.
H\bigl(P^*, Q\bigr) \leq r^2 \exp\bigl(- C
m^{1/3}\bigr) + \frac{C r}{\sqrt{m}}.
\]
\end{lem}

%
%le3 #&#
%
\begin{lem} \label{lemMany-Multi-norm}
Suppose that for $k=1,\ldots, n$,
$P_k$ is a multinomial distribution $\cM(m, \theta_{k1},\ldots,
\theta_{k \nu_k})$, where $\nu_k \leq\kappa$, $\kappa$ is a
fixed integer, $\theta_{k1} + \cdots+ \theta_{k \nu_k}=1$, and for
constants $c_0$ and $c_1$,
\[
0< c_0 \leq\min(\theta_{k1},\ldots,
\theta_{k \nu_k}) \leq\max(\theta_{k1},\ldots,
\theta_{k \nu_k}) \leq c_1 <1.
\]
Denote by $Q_k$ the multivariate normal distribution whose mean and
covariance are the same as
$P_k$. If $\nu_k \geq2$, following the same way as in Lemma \ref
{lemMulti-norm} we define $P_k^{*}$ as the convolution of $P_k$
and an independent uniform distribution on $(-1/2, 1/2)$, and if $\nu
_k \leq1$ let $P_k^*=P_k$. Assume that $P_k, P_k^*, Q_k$
for different $k$ are independent, and define product probability measures
\[
P = \prod_{k=1}^n P_k,\qquad P^* =
\prod_{k=1}^n P_k^*,\qquad Q =
\prod_{k=1}^n Q_k.
\]
Then we have
\[
H^2\bigl(P^*, Q\bigr) \leq\frac{C \kappa^2}{m} \sum
_{k=1}^n 1(\nu_k\geq2).
\]
\end{lem}

We need the following lemma on total variation distance of two joint
distributions whose proof is in the \hyperref[app]{Appendix}.
%
%le4 #&#
%
\begin{lem} \label{lemTV}
Suppose that $U_1$ and $V_1$ are discrete random variables, and
random variables $(U_1, U_2)$ and $(V_1, V_2)$ have joint distributions
$F$ and $G$, respectively.
%If $U_1$ and $U_2$ are independent with respective marginal
%distributions $F_1$ and $F_2$, and $V_1$ and $V_2$ are independent with
%respective marginal distributions $G_1$ and $G_2$, then
%% D(P_1 \times P_2, Q_1 \times Q_2) = D(P_1, Q_1) + D(P_2, Q_2).
% \| F_1 \times F_2 - G_1 \times G_2 \|_{\mathrm{TV}} \leq\| F_1 - G_1 \|_{\mathrm{TV}}
%+ \| F_2 - G_2 \|_{\mathrm{TV}}.
%In general,
Let $F(u_1, u_2)=F_1(u_1) \times F_{2|1}(u_2|u_1)$ and $G(v_1,
v_2)=G_1(v_1) \times G_{2|1}(v_2|v_1)$, where $F_1$ and $G_1$ are the
respective marginal distributions
of $U_1$ and~$V_1$, and $F_{2|1}$ and $G_{2|1}$ are the conditional
distributions of $U_2$ given $U_1$ and
$V_2$ given $V_1$, respectively. Then %from Carter (2007, section 3.1)
%
%e44 #&#
%
\begin{eqnarray}
\label{TV2} \| F - G \|_{\mathrm{TV}} &\leq&\max_x \biggl
\llvert1 - \frac{P(U_1=x)}{P(V_1=x)} \biggr\rrvert\nonumber\\[-8pt]\\[-8pt]
&&{}+ E_{F_1} \bigl[
\bigl\|F_{2|1}(\cdot|U_1) - G_{2|1}(
\cdot|V_1) \bigr\|_{\mathrm{TV}} |U_1=V_1
\bigr],\nonumber % D(P, Q) = D(P_1, Q_1) + E_{P_1} [D(P_{2|U_1},
%Q_{2|V_1})|U_1=V_1],
\end{eqnarray}
where $E_{F_1}$ denotes expectation under $F_1$, $\|F_{2|1}(\cdot|U_1)
- G_{2|1}(\cdot|V_1)\|_{\mathrm{TV}}$ denotes the total variation norm of the
difference of the two conditional distributions $F_{2|1}$ and
$G_{2|1}$, and the value of the second term on the right-hand side of
(\ref{TV2}) is clearly specified as follows:
\begin{eqnarray*}
&&
E_{F_1} \bigl[\bigl\|F_{2|1}(\cdot|U_1) -
G_{2|1}(\cdot|V_1) \bigr\|_{\mathrm{TV}} |U_1=V_1
\bigr] \\
&&\qquad= \sum_{x} \bigl\|F_{2|1}(\cdot|x) -
G_{2|1}(\cdot|x) \bigr\|_{\mathrm{TV}} P(U_1=x).
\end{eqnarray*}
\end{lem}

%Now we proceed to prove Theorem \ref{theo1}.
%s7.2 #&#
\subsection{\texorpdfstring{Proofs of Theorems \protect\ref{theo1} and \protect\ref{theo2}}
{Proofs of Theorems 1 and 2}}\label{sec7.2}

\mbox{}

\begin{pf*}{Proof of Theorem \ref{theo1}}
Denote by $\mathbb{P}_{1,n,\brho}^k$ the distribution of $(\bX_k, Y_k)$
and $\mathbb{P}_{2,n,\brho}^k$ the distribution of $(\bM_k, N_k)$,
$k=1,\ldots,n$. For different $k$, $(\bX_k, Y_k)$ from trace regression
are independent, and $(\bM_k, N_k)$ from quantum state tomography are
independent, so $\mathbb{P}_{1,n,\brho}^k$ and
$\mathbb{P}_{2,n,\brho}^k$ for different $k$ are independent, and
%
%e45 #&#
%
\begin{equation}
\label{thm-P1P2} \mathbb{P}_{1,n,\brho} = \prod_{k=1}^n
\mathbb{P}_{1,n,\brho}^k,\qquad \mathbb{P}_{2,n,\brho} = \prod
_{k=1}^n \mathbb{P}_{2,n,\brho}^k,
\end{equation}
where $\mathbb{P}_{1,n,\brho}$ and $\mathbb{P}_{2,n,\brho}$ are
given in (\ref{experiment}).

Suppose that $\bM_{k}$ has $\nu_k$ different eigenvalues, and let
$U_{ka} = \sum_{\ell=1}^m 1(R_{k\ell}=\lambda_{ka})$, $a=1,\ldots,
\nu_k$, and $\bU_k=(U_{k1},\ldots, U_{k\nu_k})^\prime$. Denote by
$\mathbb{Q}_{2,n,\brho}^k$ the distribution of $(\bM_k, \bU_k)$. If
$\nu_k \geq2$, we let $\mathbb{Q}_{2,n,\brho}^{k*}$ be the distribution
of $(\bM_k, \bU^*_{k})$, where $\bU^*_k=(U_{k1}^*,\ldots,
U_{k\nu_k}^*)^\prime$, $U_{ka}^*$ is equal to $U_{ka}$ plus an
independent uniform random variable on $(-1/2, 1/2)$, $a=1,\ldots,
\nu_k-1$ and $U_{k\nu_k}^*=m - U_{k1}^*-\cdots-U_{k, \nu_k-1}^*$. Note
that $\mathbb{P}_{2,n,\brho}^k$ is the distribution of $(\bM_k, N_k)$,
and
%
%e46 #&#
%
\begin{equation}
\label{NRU} N_k=(R_{k1} + \cdots+ R_{km})/m=(
\lambda_{k1} U_{k1} + \cdots+ \lambda_{k\nu_k}
U_{k\nu_k})/m.
\end{equation}
Analog to the expression (\ref{NRU}) of $N_k$ in terms of $\bU
_k=(U_{k1},\ldots, U_{km})^\prime$, we define
%
%e47 #&#
%
\begin{equation}
\label{N*U} N^*_k=\bigl(\lambda_{k1} U^*_{k1} +
\cdots+ \lambda_{k\nu_k} U^*_{k\nu_k}\bigr)/m,
\end{equation}
and denote by $\mathbb{P}_{2,n,\brho}^{k*}$ the distribution of $(\bM
_k, N^*_{k})$.
If $\nu_k \leq1$, let $\mathbb{Q}_{2,n,\brho}^{k*}=\mathbb
{Q}_{2,n,\brho}^{k}$ and $\mathbb{P}_{2,n,\brho}^{k*}=\mathbb
{P}_{2,n,\brho}^{k}$.
As $\mathbb{Q}_{2,n,\brho}^k$, $\mathbb{Q}_{2,n,\brho}^{k*}$, and
$\mathbb{P}_{2,n,\brho}^{k*}$ for different $k$ are independent, define
their product probability measures
%
%e48 #&#
%
\begin{equation}
\label{thm-Q2Q*2P*2}\qquad \mathbb{Q}_{2,n,\brho} = \prod
_{k=1}^n \mathbb{Q}_{2,n,\brho}^{k},\qquad
\mathbb{Q}_{2,n,\brho}^{*} = \prod_{k=1}^n
\mathbb{Q}_{2,n,\brho}^{k*},\qquad \mathbb{P}_{2,n,\brho}^* = \prod
_{k=1}^n \mathbb{P}_{2,n,\brho}^{k*}.
\end{equation}
Note that, since $\bU_k$ and $(R_{k1},\ldots, R_{km})$ have a one to
one correspondence, and the two statistical experiments formed by the
distribution of $(\bM_k, \bU_k)$ and the distribution of $(\bM_k,
R_{k1},\ldots, R_{km})$ have zero deficiency distance,
without confusion we abuse the notation $\mathbb{Q}_{2,n,\brho}$
by using it here for the joint distribution of $(\bM_k, \bU_k)$,
$k=1,\ldots, n$, as well as in (\ref{experimentQ2}) for the joint
distribution of
$(\bM_k, R_{k1},\ldots, R_{km})$, $k=1,\ldots, n$.

Given $\bM_k=\bB_{j_k}$, let $\nu_k=r_{j_k}$, and
$\bU_k=(U_{k1},\ldots, U_{k r_{j_k}})^\prime$ follows a multinomial
distribution $\cM(m, \operatorname{tr}(\bQ_{j_k 1} \brho),\ldots,\operatorname{tr}(\bQ_{j_k
r_{j_k}} \brho) )$,
where $r_j$ and $\bQ_{j a}$ are defined in (\ref{basis-diagonal}), and
\begin{eqnarray}
E(U_{ka}|\bM_k=\bB_{j_k}) &=& m \operatorname{tr}(
\bQ_{j_k a} \brho),\nonumber\\
\operatorname{Var}(U_{ka} |\bM_k=
\bB_{j_k}) &=& m \operatorname{tr}(\bQ_{j_k a} \brho) \bigl[1-\operatorname{tr}(\bQ
_{j_k a} \brho)\bigr],
\nonumber\\
\operatorname{Cov}(U_{ka},U_{kb} | \bM_k=
\bB_{j_k}) &=& - m \operatorname{tr}(\bQ_{j_k a} \brho) \operatorname{tr}(\bQ_{j_k b}
\brho),\nonumber\\
&&\eqntext{a \neq b, a, b =1,\ldots, r_{j_k}.}
\end{eqnarray}
Then
\begin{eqnarray*}
E(N_k|\bM_k=\bB_{j_k}) &=& \sum
_{a=1}^{r_{j_k}} \lambda_{j_k a} \operatorname{tr}(
\bQ_{j_k a} \brho) = \operatorname{tr}(\bB_{j_k} \brho)=\operatorname{tr}(\bM_{k}
\brho),
\\
\operatorname{Var}(N_k|\bM_k=\bB_{j_k}) &=&
\frac{1}{m}\sum_{a=1}^{r_{j_k}}
\lambda_{j_k a}^2 \operatorname{tr}(\bQ_{j_k a} \brho) \bigl[ 1 - \operatorname{tr}(
\bQ_{j_k a} \brho)\bigr]
\\
&&{} - \frac{2}{m}\sum_{a=1}^{r_{j_k}} \sum
_{b=a+1}^{r_{j_k}} \lambda_{j_k a}
\lambda_{j_k b} \operatorname{tr}(\bQ_{j_k a} \brho) \operatorname{tr}(\bQ_{j_k
b} \brho)
\\
&=& \frac{1}{m} \bigl\{\operatorname{tr}\bigl(\bB_{j_k}^2 \brho
\bigr) - \bigl[\operatorname{tr}(\bB_{j_k} \brho)\bigr]^2\bigr\}\\
&=&
\frac{1}{m} \bigl\{\operatorname{tr}\bigl(\bM_k^2 \brho\bigr) -
\bigl[\operatorname{tr}(\bM_k \brho)\bigr]^2\bigr\}.
\end{eqnarray*}
From (\ref{trace-regression2}) and (\ref{trace-regression1}), we have
that given $\bX_k=\bB_{j_k}$, $r^X_k=r_{j_k}$,
and multivariate normal random vector\vadjust{\goodbreak} $\bV_k=(V_{k1},\ldots, V_{k
r_{j_k}})^\prime= m (y_{k1},\ldots, y_{k r_{j_k}})^\prime$ has
conditional mean and conditional covariance matching those of $\bU
_k=(U_{k1},\ldots, U_{k r_{j_k}})^\prime$.
With $\bX_k=\bB_{j_k}$ we may rewrite (\ref{trace-regression1}) and
(\ref{trace-regression0}) as follows:
%
%e49 #&#
%
\begin{eqnarray}
\label{epsilonV} V_{ka} &=& m \operatorname{tr}(\bQ_{j_k a} \brho) + m
z_{ka},\qquad a = 1,\ldots, r_{j_k},\nonumber\\[-8pt]\\[-8pt]
Y_k &=& %\sum_{a=1}^{r_{j_k}}
\frac{1}{m} \sum_{a=1}^{r_{j_k}}
\lambda_{ka} V_{ka},\qquad \varepsilon_k= \sum
_{a=1}^{r_{j_k}} \lambda_{ka} z_{ka}.\nonumber
\end{eqnarray}
%
%where $z_{ka}$ is defined in (\ref{trace-regression1}).
Denote by $\mathbb{Q}_{1,n,\brho}^k$ the distribution of $(\bX_k,
\bV_{k})$. Then $\mathbb{Q}_{1,n,\brho}^k$ for different $k$ are
independent, and
%
%e50 #&#
%
\begin{equation}
\label{thm-Q1} \mathbb{Q}_{1,n,\brho} = \prod_{k=1}^n
\mathbb{Q}_{1,n,\brho}^{k},
\end{equation}
where $\mathbb{Q}_{1,n,\brho}$ is the joint distribution of $(\bX_k,
V_{k1},\ldots, V_{k r^X_k})$, $k=1,\ldots, n$.
Note that, since $\bV_k=(V_{k1},\ldots, V_{k r_{j_k}})^\prime= m
(y_{k1},\ldots, y_{k r_{j_k}})^\prime$, and the two statistical
experiments formed by
the distribution of $(\bX_k, V_{k1},\ldots, V_{k r_{j_k}})$ and the
distribution of $(\bX_k, y_{k1},\ldots, y_{k r_{j_k}})$ have zero
deficiency distance,
without confusion we abuse the notation $\mathbb{Q}_{1,n,\brho}$ by
using it here for the joint distribution of
$(\bX_k, V_{k1},\ldots, V_{k r^X_k})$, $k=1,\ldots, n$, as well as
in (\ref{experimentQ1}) for the joint distribution of
$(\bX_k, y_{k1},\ldots, y_{k r^X_k})$, $k=1,\ldots, n$.

Conditional on $\bM_k=\bB_{j_k}$, for $k=1,\ldots, n$, if $|\cI
_{j_k}(\brho)| \leq1$, $\mathbb{Q}_{1,n,\brho}^{k}$ and $\mathbb
{Q}_{2,n,\brho}^{k}$ are the same degenerate distribution;
if $|\cI_{j_k}(\brho)| \geq2$, $\mathbb{Q}_{2,n,\brho}^{k}$ is a
multinomial distribution with $\mathbb{Q}_{2,n,\brho}^{k*}$
its uniform perturbation, and $\mathbb{Q}_{1,n,\brho}^{k}$ is a
multivariate normal distribution with mean and covariance matching
those of $\mathbb{Q}_{2,n,\brho}^{k}$.
%Following Carter (\citeyear{Car02}, Theorem 1 and its proof) we have that
%conditional on $\bX_k=\bM_k=\bB_{j_k}$,
% \| \mathbb{Q}_{1,n,\brho}^{k} - \mathbb{Q}_{2,n,\brho}^{k*} \|_{\mathrm{TV}}
%where $D(\cdot, \cdot)$ is the Kullback{\chr"96}Leibler divergence.
Thus applying Lemma \ref{lemMany-Multi-norm}, we obtain that
given $(\bX_1,\ldots, \bX_n)=(\bM_1,\ldots, \bM_n)=(\bB_{j_1},\ldots,
\bB_{j_n})$,
%
%e51 #&#
%
\begin{eqnarray}
\label{QQ}\quad %&& \| \mathbb{Q}_{1,n,\brho} - \mathbb{Q}_{2,n,\brho}^{*}
%%&& D(\mathbb{Q}_{1,n,\brho}, \mathbb{Q}_{2,n,\brho}^{*}) \leq
\bigl\| \mathbb{Q}_{1,n,\brho} - \mathbb{Q}_{2,n,\brho}^{*}
\bigr\|_{\mathrm{TV}}^2 &\leq& H^2\bigl(\mathbb{Q}_{1,n,\brho},
\mathbb{Q}_{2,n,\brho}^{*} \bigr) \leq\frac{C \kappa^2}{m} \sum
_{k=1}^n 1\bigl(\bigl|\cI_{j_k}(\brho)\bigr|
\geq2\bigr), %\leq\frac{s C_\kappa}{\sqrt{m}},
\end{eqnarray}
where the first inequality is due to (\ref{TV-Hellinger-KL}).
As (\ref{N*U}) and (\ref{epsilonV}) imply that $N_k^*$ and $Y_k$ %$
are the same weighted averages of
components of $\bU^*_k$ and $\bV_k$, respectively, $\mathbb
{P}_{1,n,\brho}$ and $\mathbb{P}_{2,n,\brho}^*$ are the same respective
marginal probability measures of $\mathbb{Q}_{1,n,\brho}$ and
$\mathbb{Q}_{2,n,\brho}^*$.
Hence, conditional on $(\bX_1,\ldots, \bX_n)=(\bM_1,\ldots, \bM_n)$,
%
%e52 #&#
%
\begin{equation}
\label{PQ} %D(\mathbb{P}_{1,n,\brho}, \mathbb{P}_{2,n,\brho}^*) \leq D(
\bigl\| \mathbb{P}_{1,n,\brho} -
\mathbb{P}_{2,n,\brho}^* \bigr\|_{\mathrm{TV}} \leq\bigl\| \mathbb{Q}_{1,n,\brho} -
\mathbb{Q}_{2,n,\brho}^* \bigr\|_{\mathrm{TV}}.
\end{equation}
With $\bX_k$ and $\bM_k$ are sampled from $\cB$ according to
distributions $\Pi$ and $\Xi$, respectively, we have
%
%e53 #&#
%
\begin{eqnarray}
\label{PP}
&&\bigl\|\mathbb{P}_{1,n,\brho} - \mathbb{P}_{2,n,\brho}^{*}
\bigr\|_{\mathrm{TV}} \nonumber\\
&&\qquad\leq\max_{1 \leq j \leq p} \biggl\llvert1 -
\frac{\Pi^n(j)}{\Xi^n(j)} \biggr\rrvert
\nonumber
\\
&&\qquad\quad{} + E_\Pi\bigl(E_\Pi\bigl[ \bigl\|\mathbb
{P}_{1,n,\brho}-\mathbb{P}_{2,n,\brho}^{*} \bigr\|_{\mathrm{TV}} |
\bX_{1}=\bM_1,\ldots, \bX_n=
\bM_n \bigr] \bigr)
\\
%&& D(\mathbb{P}_{1,n,\brho}, \mathbb{P}_{2,n,\brho}^{*}) = D(\Pi^n,
&&\qquad\leq n \max_{1 \leq j \leq p}\biggl
\llvert1 - \frac{\Pi(j)}{\Xi(j)} \biggr\rrvert\nonumber\\
&&\qquad\quad{}+ E_\Pi\bigl(
E_\Pi\bigl[ \bigl\| \mathbb{Q}_{1,n,\brho} - \mathbb{Q}_{2,n,\brho}^{*}
\bigr\|_{\mathrm{TV}} | \bX_{1}=\bM_1,\ldots,
\bX_n=\bM_n \bigr] \bigr)
\nonumber
\\
&&\qquad\leq n \gamma_p + \frac{C \kappa}{\sqrt{m}} E_\Pi\Biggl(
\Biggl[ \sum_{k=1}^n 1\bigl(\bigl|
\cI_{j_k}(\brho)\bigr|\geq2\bigr) \Biggr]^{1/2} \Biggr)
\nonumber
\\
&&\qquad\leq n \gamma_p + \frac{C \kappa}{\sqrt{m}} \Biggl( \sum
_{k=1}^n E_\Pi\bigl[ 1\bigl(\bigl|
\cI_{j_k}(\brho)\bigr|\geq2\bigr) \bigr] \Biggr)^{1/2}
\nonumber
\\
&&\qquad\leq n \gamma_p + \frac{C \kappa}{\sqrt{m}} \Biggl( \sum
_{k=1}^n \sum_{j=1}^p
\Pi(j) 1\bigl(\bigl|\cI_{j}(\brho)\bigr|\geq2\bigr) \Biggr)^{1/2}
\nonumber
\\
&&\qquad= n \gamma_p + \frac{C \kappa}{\sqrt{m}} \Biggl( n \sum
_{j=1}^p \Pi(j) 1\bigl(\bigl|\cI_{j}(\brho)\bigr|
\geq2\bigr) \Biggr)^{1/2}
\nonumber
\\
&&\qquad\leq n \gamma_p + C \kappa\biggl(\frac{n \zeta_p}{m}
\biggr)^{1/2},\nonumber
\end{eqnarray}
where the first three inequalities are, respectively, from Lemma \ref
{lemTV}, (\ref{PQ}) and (\ref{QQ}),
the fourth inequality is applying H\"older's inequality, and the fifth
inequality is due the fact that $\bX_k$ and $\bM_k$ are the
i.i.d. sample from $\cB$. Combining (\ref{kernel}) and (\ref{PP}),
we obtain
%
%e54 #&#
%
\begin{eqnarray}
\label{deltaP21} \delta(\cP_{2n}, \cP_{1n}) &\leq&\inf
_{K} \sup_{\brho\in\Theta
} \bigl\|\mathbb{P}_{1,n,\brho}
- K (\mathbb{P}_{2,n,\brho})\bigr\|_{\mathrm{TV}}
\nonumber
\\
&\leq&\sup_{\brho\in\Theta} \bigl\| \mathbb{P}_{1,n,\brho} - \mathbb
{P}_{2,n,\brho}^* \bigr\|_{\mathrm{TV}}
\\
&\leq& n \gamma_p + C \kappa\biggl(\frac{n \zeta_p}{m}
\biggr)^{1/2}.\nonumber
\end{eqnarray}
To bound $\delta(\cP_{1n}, \cP_{2n})$, we %follow Carter (\citeyear{Car02},
%section 5) by using
employ a round-off procedure to invert the uniform perturbation used to obtain
$\mathbb{Q}_{2,n,\brho}^{*}$ and $\mathbb{P}_{2,n,\brho}^{*}$ in
(\ref{thm-Q2Q*2P*2}) [also see Carter (\citeyear{Car02}), Section~5].
Specifically let $\bV_{k}^*=(V_{k1}^*,\ldots, V_{k\nu_k}^*)^\prime$,
where $V_{ka}^*$ is a random vector obtained by rounding $V_{ka}$ off
to the nearest integer, $a=1,\ldots, \nu_k-1$, and
$V_{k \nu_k}^* = m - V_{k1}^* - \cdots- V_{k, \nu_k-1}^*$. %applying
%the round-off procedure to $\bV_{ka}$
Denote by $\mathbb{Q}_{1,n,\brho}^{k*}$ the distribution of
$(\bX_k, \bV_{k}^*)$ and $\mathbb{P}^{k*}_{1,n,\brho}$ the
distribution of $(\bX_k, (\lambda_{k1} V_{k1}^* + \cdots+ \lambda
_{k \nu_k} V_{k \nu_k}^*)/m)$,
and let
%
%e55 #&#
%
\begin{equation}
\label{thm-Q*1P*1} \mathbb{Q}^*_{1,n,\brho} = \prod
_{k=1}^n \mathbb{Q}_{1,n,\brho
}^{k*},\qquad
\mathbb{P}_{1,n,\brho}^* = \prod_{k=1}^n
\mathbb{P}_{1,n,\brho}^{k*}.
\end{equation}
It is easy to see that for any integer-valued random variable $W$,
\[
\mbox{round-off of } \bigl[W + \operatorname{uniform}(-1/2,1/2) \bigr] = W,
\]
and thus the round-off procedure inverts the uniform perturbation
procedure. Denote by $K_0$ and $K_1$ the uniform perturbation
and the round-off procedure, respectively. Then from (\ref
{thm-Q2Q*2P*2}), (\ref{thm-Q1}) and (\ref{thm-Q*1P*1}) we have
%
%e56 #&#
%
\begin{eqnarray}
\label{thm-K1K0} K_1(\mathbb{Q}_{1,n,\brho}) &=&
\mathbb{Q}^*_{1,n,\brho},\qquad K_0(\mathbb{Q}_{2,n,\brho}) =
\mathbb{Q}_{2,n,\brho}^*,\nonumber\\[-8pt]\\[-8pt]
K_1\bigl[K_0(\mathbb{Q}_{2,n,\brho})\bigr] &=& K_1\bigl[\mathbb{Q}_{2,n,\brho}^*
\bigr]= \mathbb{Q}_{2,n,\brho}.\nonumber
\end{eqnarray}
%
%Similar to Section 5 of \citet{Car02}
From (\ref{thm-K1K0}), we show that conditional on $(\bX_1,\ldots,
\bX_n)=(\bM_1,\ldots, \bM_n)$,
%
%e57 #&#
%
\begin{eqnarray}
\label{Q*1Q2} \bigl\| \mathbb{Q}^*_{1,n,\brho} - \mathbb{Q}_{2,n,\brho}
\bigr\|_{\mathrm{TV}} &=& \bigl\| K_1(\mathbb{Q}_{1,n,\brho}) -
K_1\bigl[K_0(\mathbb{Q}_{2,n,\brho})\bigr]
\bigr\|_{\mathrm{TV}} \nonumber\\
&=& \bigl\| K_1\bigl[ \mathbb{Q}_{1,n,\brho} -
K_0(\mathbb{Q}_{2,n,\brho})\bigr] \bigr\| _{\mathrm{TV}}
\nonumber\\[-8pt]\\[-8pt]
&\leq&\bigl\| \mathbb{Q}_{1,n,\brho} - K_0(\mathbb{Q}_{2,n,\brho})
\bigr\| _{\mathrm{TV}} \nonumber\\
&=& \bigl\| \mathbb{Q}_{1,n,\brho} - \mathbb{Q}_{2,n,\brho}^*
\bigr\|_{\mathrm{TV}},\nonumber
\end{eqnarray}
which is bounded by (\ref{QQ}).
Using the same arguments for showing (\ref{PQ}) and (\ref{PP}) we
derive from (\ref{QQ}) and (\ref{Q*1Q2}) the following result:
%
%e58 #&#
%
\begin{eqnarray}
\label{PPreverse}
&&
\bigl\|\mathbb{P}_{1,n,\brho}^* - \mathbb{P}_{2,n,\brho}
\bigr\|_{\mathrm{TV}} \nonumber\\
&&\qquad\leq n \max_{1\leq j \leq p} \biggl\llvert1 -
\frac{\Xi(j)}{\Pi(j)} \biggr\rrvert+ \frac{C \kappa}{\sqrt{m}} \Biggl(
n \sum
_{j=1}^p \Xi(j) 1\bigl(\bigl|\cI_{j}(\brho)\bigr|
\geq2\bigr) \Biggr)^{1/2}
\\
&&\qquad\leq n \delta_p + C \kappa\biggl(\frac{n \zeta_p}{m}
\biggr)^{1/2},\nonumber
\end{eqnarray}
and applying (\ref{kernel}) we conclude
%
%e59 #&#
%
\begin{eqnarray}
\label{deltaP12} \delta(\cP_{1n}, \cP_{2n}) &\leq& \inf
_{K} \sup_{\brho\in\Theta
} \bigl\| K (\mathbb{P}_{1,n,\brho})
- \mathbb{P}_{2,n,\brho}\bigr\|_{\mathrm{TV}}
\nonumber
\\
&\leq&\sup_{\brho\in\Theta} \bigl\| \mathbb{P}_{1,n,\brho}^* -
\mathbb{P}_{2,n,\brho} \bigr\|_{\mathrm{TV}}
\\
&\leq& n \delta_p + C \kappa\biggl(\frac{n \zeta_p}{m}
\biggr)^{1/2}.\nonumber %\leq n \max_j \left| 1 - \frac{\Xi(j)}{\Pi(j)} \right|
%+ \frac{C
\end{eqnarray}
Collecting together the deficiency bounds in (\ref{deltaP21}) and
(\ref{deltaP12}) we establish (\ref{randomPP}) to bound the
deficiency distance
$\Delta(\cP_{1n}, \cP_{2n})$ for the random design case. For the
special case of $\Pi(j)=\Xi(j)=1/p$, $\gamma_p=0$ and
\begin{eqnarray*}
\zeta_p&=&\max\Biggl\{ \sum_{j=1}^p
\Pi(j) 1\bigl(\bigl|\cI_{j}(\brho)\bigr|\geq2\bigr), \sum
_{j=1}^p \Xi(j) 1\bigl(\bigl|\cI_{j}(\brho)\bigr|
\geq2\bigr) \Biggr\} \\
&=& \frac
{1}{p} \sum_{j=1}^p
1\bigl(\bigl|\cI_{j}(\brho)\bigr|\geq2\bigr).
\end{eqnarray*}
The result (\ref{randomPP1}) follows.
\end{pf*}

For the fixed design case, the arguments for proving (\ref{fixPP}) are
the same except for now we simply combine (\ref{QQ}), (\ref{PQ}) and
(\ref{Q*1Q2}) but no need for (\ref{PP}) and
(\ref{PPreverse}).

\begin{pf*}{Proof of Theorem \ref{theo2}}
The proof of Theorem \ref{theo1} has essentially
established Theorem \ref{theo2}. All we need is to modify the arguments
as follows. As in the derivation of (\ref{PP}) we apply Lemma \ref
{lemTV} directly to $\mathbb{Q}_{1,n,\brho}$ and $\mathbb
{Q}_{2,n,\brho}^{*}$ and use (\ref{QQ}) to get
\begin{eqnarray*}
&&
\bigl\|\mathbb{Q}_{1,n,\brho} - \mathbb{Q}_{2,n,\brho}^{*}
\bigr\|_{\mathrm{TV}} \\
&&\qquad\leq\max_{1 \leq j \leq p} \biggl\llvert1 -
\frac{\Pi^n(j)}{\Xi^n(j)} \biggr\rrvert
\\
&&\qquad\quad{} + E_\Pi\bigl(E_\Pi\bigl[ \bigl\|\mathbb
{Q}_{1,n,\brho}-\mathbb{Q}_{2,n,\brho}^{*} \bigr\|_{\mathrm{TV}} |
\bX_{1}=\bM_1,\ldots, \bX_n=
\bM_n \bigr] \bigr)
\\
%&&\leq n \max_{1 \leq j \leq p}\left| 1 - \frac{\Pi(j)}{\Xi(j)}
% \bX_{1}=\bM_1,\ldots, \bX_n=\bM_n \right]\right) \nonumber\\
&&\qquad\leq n \gamma_p +
C \kappa\biggl(\frac{n \zeta_p}{m} \biggr)^{1/2},
\end{eqnarray*}
and then we obtain, instead of (\ref{deltaP21}), the following result:
%
%e60 #&#
%
\begin{eqnarray}
\label{deltaP21-thm2} \delta(\cQ_{2n}, \cQ_{1n}) &\leq&\inf
_{K} \sup_{\brho\in\Theta
} \bigl\| \mathbb{Q}_{1,n,\brho}
- K (\mathbb{Q}_{2,n,\brho}) \bigr\|_{\mathrm{TV}}
\nonumber
\\
&\leq&\sup_{\brho\in\Theta} \bigl\| \mathbb{Q}_{1,n,\brho} - \mathbb
{Q}_{2,n,\brho}^* \bigr\|_{\mathrm{TV}}
\\
&\leq& n \gamma_p + C \kappa\biggl(\frac{n \zeta_p}{m}
\biggr)^{1/2}.\nonumber
\end{eqnarray}
As in the derivation of (\ref{PPreverse}), we apply Lemma \ref{lemTV}
to $\mathbb{Q}_{1,n,\brho}^{*}$ and $\mathbb{Q}_{2,n,\brho}$ and
use (\ref{QQ}) and (\ref{Q*1Q2}) to get
\begin{eqnarray*}
\bigl\|\mathbb{Q}_{1,n,\brho}^* - \mathbb{Q}_{2,n,\brho} \bigr\|_{\mathrm{TV}}
&\leq& n \max_{1\leq j \leq p} \biggl\llvert1 - \frac{\Xi(j)}{\Pi(j)}
\biggr
\rrvert\\
&&{}+ \frac{C \kappa}{\sqrt{m}} \Biggl( n \sum_{j=1}^p
\Xi(j) 1\bigl(\bigl|\cI_{j}(\brho)\bigr|\geq2\bigr) \Biggr)^{1/2}
\\
&\leq& n \delta_p + C \kappa\biggl(\frac{n \zeta_p}{m}
\biggr)^{1/2},
\end{eqnarray*}
and then we obtain, instead of (\ref{deltaP12}), the following result:
%
%e61 #&#
%
\begin{eqnarray}
\label{deltaP12-thm2} \delta(\cQ_{1n}, \cQ_{2n}) &\leq&\inf
_{K} \sup_{\brho\in\Theta
} \bigl\| K (\mathbb{Q}_{1,n,\brho})
- \mathbb{Q}_{2,n,\brho}\bigr\|_{\mathrm{TV}}
\nonumber
\\
&\leq&\sup_{\brho\in\Theta} \bigl\| \mathbb{Q}_{1,n,\brho}^* -
\mathbb{Q}_{2,n,\brho} \bigr\|_{\mathrm{TV}}
\\
&\leq& n \delta_p + C \kappa\biggl(\frac{n \zeta_p}{m}
\biggr)^{1/2}.\nonumber
\end{eqnarray}
Putting together the deficiency bounds in (\ref{deltaP21-thm2}) and
(\ref{deltaP12-thm2}) we establish (\ref{randomQQ}) to bound the
deficiency distance
$\Delta(\cQ_{1n}, \cQ_{2n})$ for the random design case.
\end{pf*}

%s7.3 #&#
\subsection{Proofs of corollaries}\label{sec7.3}

To prove corollaries, from Theorems \ref{theo1} and \ref{theo2} we need
to show the given
bounds on $\zeta_p$ and then substitute them into (\ref{randomPP1})
and (\ref{randomQQ1}).
Below we will derive $\zeta_p$ for each case.

\begin{pf*}{Proof of Corollary \ref{cor-sparse}}
We first analyze the eigen-structures of basis matrices given by (\ref
{Hermitian-basis}). For diagonal basis matrix $\bB_j$ with $1$ on
$(\ell, \ell)$ entry and $0$ elsewhere, its eigenvalues
are $1$ and $0$. Corresponding to eigenvalue $1$, the eigenvector is
$\be_{\ell}$, and corresponding to
eigenvalue $0$, the eigen-space is the orthogonal complement of $\operatorname{span}\{
\be_\ell\}$. %$\be_1,\ldots, \be_{\ell-1}, \be_{\ell-1},\ldots,
Denote by $\bQ_{j0}$ and $\bQ_{j1}$ the projections on the
eigen-spaces corresponding to eigenvalues $0$ and $1$,
respectively.\vspace*{2pt}

For real symmetric nondiagonal $\bB_j$ with $1/\sqrt{2}$ on $(\ell
_1,\ell_2)$ and $(\ell_2, \ell_1)$ entries and $0$ elsewhere, the
eigenvalues are
$1$, $-1$ and $0$. Corresponding to eigenvalues $\pm1$, the
eigenvectors are $(\be_{\ell_1} \pm\be_{\ell_2})/\sqrt{2}$,
respectively, and corresponding to
eigenvalue $0$, the eigen-space is the orthogonal complement of $\operatorname{span}\{
\be_{\ell_1} \pm\be_{\ell_2} \}$.
%$\be_1,\ldots, \be_{\ell_1-1}, \be_{\ell_1 + 1},\ldots, \be_{
Denote by $\bQ_{j0}$, $\bQ_{j1}$ and $\bQ_{j,-1}$ the projections on
the eigen-spaces corresponding to eigenvalues $0$, $1$ and $-1$, respectively.

For imaginary Hermitian $\bB_j$ with $-\sqrt{-1}/\sqrt{2}$ on $(\ell
_1,\ell_2)$ entry, $\sqrt{-1}/\sqrt{2}$ on $(\ell_2, \ell_1)$
entry and $0$ elsewhere,
the eigenvalues are $1$, $-1$ and $0$. Corresponding to eigenvalues
$\pm1$, the eigenvector are $(\be_{\ell_1} \pm\sqrt{-1} \be_{\ell
_2})/\sqrt{2}$,
respectively, and corresponding to eigenvalue $0$, the eigen-space is
the orthogonal complement of $\operatorname{span}\{ \be_{\ell_1} \pm\sqrt{-1} \be
_{\ell_2} \}$.
%$ (\be_1,\ldots, \be_{\ell_1-1}, \be_{\ell_1 + 1},\ldots, \be_{
Denote by $\bQ_{j0}$, $\bQ_{j1}$ and $\bQ_{j,-1}$ the projections on
the eigen-spaces corresponding to eigenvalues
$0$, $1$ and $-1$, respectively.

For diagonal $\bB_j$ with $1$ on $(\ell,\ell)$ entry, it is a
binomial case,
\[
\operatorname{tr}(\brho\bQ_{j0} ) = 1 - \operatorname{tr}(\brho\bQ_{j1}),\qquad
\operatorname{tr}(\brho\bQ
_{j1}) = \be_{\ell}^\prime\brho\be_{\ell} =
\rho_{\ell\ell}
\]
and
\[
\bigl|\cI_j(\brho)\bigr| = 2\cdot1\bigl(0< \operatorname{tr}(\brho\bQ_{j1}) <1
\bigr) + 1\bigl(\operatorname{tr}(\brho\bQ_{j1} )=1\bigr) + 1\bigl(\operatorname{tr}(\brho
\bQ_{j1})=0\bigr).
\]
In order to have $|\cI_j(\brho)| \geq2$, we need $\operatorname{tr}(\brho\bQ
_{j1}) = \rho_{\ell\ell} \in(0,1)$.
Since $\brho$ has at most $s_d$ nonzero diagonal entries, among all
the $d$ diagonal matrices $\bB_j$ there are at most $s_d$ of diagonal
matrices $\bB_j$
for which it is possible to have $\operatorname{tr}(\brho\bQ_{j 1}) \in(0,1)$ and
thus $|\cI_j(\brho)|\geq2$.

For nondiagonal $\bB_j$, it is a trinomial case,
\[
\operatorname{tr}(\brho\bQ_{j0}) = 1 - \operatorname{tr}(\brho\bQ_{j1}) - \operatorname{tr}(\brho
\bQ_{j,-1}),
\]
and $\operatorname{tr}(\brho\bQ_{j\pm1})$ depend on whether $\bB_j$ is real or complex.

For real symmetric nondiagonal $\bB_j$ with $1/\sqrt{2}$ on $(\ell
_1,\ell_2)$ and $(\ell_2,\ell_1)$ entries,
\begin{eqnarray*}
\operatorname{tr}(\brho\bQ_{j\pm1}) &=& (\be_{\ell_1} \pm\be_{\ell_2})^\prime
\brho(\be_{\ell_1} \pm\be_{\ell_2})/2
\\
&=& (\rho_{\ell_1 \ell_1} + \rho_{\ell_2 \ell_2} \pm\rho_{\ell
_1 \ell_2} \pm
\rho_{\ell_2 \ell_1})/2
\\
&=& \frac{1}{2} (1, \pm1) \pmatrix{\rho_{\ell_1 \ell_1} &
\rho_{\ell_1 \ell_2}
\cr
\rho_{\ell_2 \ell_1} & \rho_{\ell_2 \ell_2}} \pmatrix{ 1
\cr
\pm1};
\end{eqnarray*}
and for imaginary Hermitian nondiagonal $\bB_j$ with $-\sqrt{-1}/\sqrt
{2}$ on $(\ell_1,\ell_2)$ entry and $\sqrt{-1}/\sqrt{2}$
on $(\ell_2,\ell_1)$ entry,
\begin{eqnarray*}
\operatorname{tr}(\brho\bQ_{j\pm1}) &=& (\be_{\ell_1} \pm\sqrt{-1}
\be_{\ell
_2})^\dagger\brho(\be_{\ell_1} \pm\sqrt{-1}
\be_{\ell_2})/2
\\
&=& (\rho_{\ell_1 \ell_1} + \rho_{\ell_2 \ell_2} \pm\sqrt{-1}
\rho_{\ell_1 \ell_2} \mp\sqrt{-1} \rho_{\ell_2 \ell_1})/2
\\
&=& \frac{1}{2} (1, \mp\sqrt{-1}) \pmatrix{\rho_{\ell_1 \ell_1} &
\rho_{\ell_1 \ell_2}
\cr
\rho_{\ell_2 \ell_1} & \rho_{\ell_2 \ell_2}} \pmatrix{1
\cr
\pm\sqrt{-1}}.
\end{eqnarray*}
As $\brho$ is semi-positive with trace $1$, matrix
\[
\pmatrix{\rho_{\ell_1 \ell_1} & \rho_{\ell_1 \ell_2}
\cr
\rho_{\ell_2 \ell_1} &
\rho_{\ell_2 \ell_2}}
\]
must be semi-positive with trace no more than $1$. Of $\rho_{\ell_1
\ell_1}$ and $\rho_{\ell_2 \ell_2}$, if one of them is zero, the
semi-positiveness implies
$\rho_{\ell_1 \ell_2}=\rho_{\ell_2 \ell_1}=0$. Thus, the $2$ by
$2$ matrix has four scenarios: %$\rho_{\ell_1 \ell_1} \neq0$ and $
%$\rho_{\ell_1 \ell_1} \neq0$ and $\rho_{\ell_2 \ell_2}=0$, $\rho_{
%
\[
\pmatrix{\rho_{\ell_1 \ell_1} & \rho_{\ell_1 \ell_2}
\cr
\rho_{\ell_2 \ell_1} &
\rho_{\ell_2 \ell_2}} \quad\mbox{or}\quad \pmatrix{ \rho_{\ell_1 \ell_1} & 0
\cr
0 & 0}
\quad\mbox{or}\quad \pmatrix{0 & 0
\cr
0 & \rho_{\ell_2 \ell_2}} \quad\mbox{or}\quad \pmatrix{0 & 0
\cr
0 & 0 }.
\]
For the last three scenarios under both real symmetric and imaginary
Hermitian cases, we obtain
\[
\operatorname{tr}(\brho\bQ_{j1}) = \operatorname{tr}(\brho\bQ_{j,-1}) =
\rho_{\ell_1 \ell_1}/2 \quad\mbox{or}\quad \rho_{\ell_2 \ell_2}/2 \quad\mbox{or}\quad 0.
\]
For both real symmetric and imaginary Hermitian cases, in order to have
\mbox{$|\cI_j(\brho)|\geq2$} possible, at lease one of $\rho_{\ell_1 \ell
_1}$ and $\rho_{\ell_2 \ell_2}$
needs to\vspace*{1pt} be nonzero. Since $\brho$ has at most $s_d$ nonzero
diagonal entries,\vspace*{1pt}
among $(d^2-d)/2$ real symmetric nondiagonal matrices $\bB_j$ [or
$(d^2-d)/2$ imaginary Hermitian nondiagonal matrices $\bB_j$],
there are at most $d s_d - s_d (s_d +1)/2$ of real symmetric
nondiagonal $\bB_j$ (or imaginary Hermitian nondiagonal matrices
$\bB_j$)
for which it is possible to have $\operatorname{tr}(\brho\bQ_{j1}) \in(0,1)$ or
$\operatorname{tr}(\brho\bQ_{j, -1}) \in(0,1)$ and thus $|\cI_j(\brho)|\geq2$.

Finally, for $\brho\in\Theta_s$, putting together the results on the
number of $\bB_j$ for which it is possible to have $|\cI_j(\brho
)|\geq2$ in the diagonal, real symmetric and imaginary Hermitian
cases, we conclude
\[
\sum_{j=1}^p 1\bigl(\bigl|\cI_j(
\brho)\bigr|\geq2\bigr) \leq d s_d - s_d (s_d+1) +
s_d \leq d s_d
\]
and
\[
\zeta_p = \max_{\brho\in\Theta_s} \Biggl\{\frac{1}{p}
\sum_{j=1}^p 1\bigl(\bigl|\cI_j(
\brho)\bigr|\geq2\bigr) \Biggr\} \leq\frac{s_d}{d}. % =
\]
\upqed
\end{pf*}

\begin{pf*}{Proof of Corollary \ref{cor-sparse-Pauli}}
The Pauli basis (\ref{Pauli-basis}) has $p=d^2$ matrices with $d=2^b$.
We identify index $j =1,\ldots, p$ with
$(\ell_1, \ell_2,\ldots, \ell_b) \in\{0, 1, 2, 3\}^b$, $j=1$
corresponds to $\ell_1= \cdots= \ell_b=0$, and $\bB_1=\bI_d$.
In two dimensions, Pauli matrices satisfy $\operatorname{tr}(\sigma_0)=2$, and
$\operatorname{tr}(\sigma_1)=\operatorname{tr}(\sigma_2)=\operatorname{tr}(\sigma_3)=0$.
Consider $\bB_j = \bsigma_{\ell_1} \otimes\bsigma_{\ell_2}
\otimes\cdots\otimes\bsigma_{\ell_b}$.
$\operatorname{tr}(\bB_{j})=\operatorname{tr}(\bsigma_{\ell_1}) \operatorname{tr}(\bsigma_{\ell_2}) \cdots
\operatorname{tr}(\bsigma_{\ell_b})$;
$\operatorname{tr}(\bB_{1})=d$; for $j \neq1$ [or $(\ell_1,\ldots, \ell_b) \neq
(0,\ldots, 0)$], $\operatorname{tr}(\bB_{j})=0$ and
$\bB_j$ has eigenvalues $\pm1$. Denote by $\bQ_{j \pm}$ the
projections onto the eigen-spaces corresponding to eigenvalues $\pm1$,
respectively.
Then for $j \neq1$,
\begin{eqnarray*}
\bB_{j} &=& \bQ_{j +} - \bQ_{j -},\qquad
\bB_j^2=\bQ_{j +} + \bQ_{j -}=
\bI_d,\qquad \bB_{j} \bQ_{j \pm} = \pm
\bQ_{j \pm}^2= \pm\bQ_{j \pm},
\\
0 &=& \operatorname{tr}(\bB_{j}) = \operatorname{tr}(\bQ_{j +}) - \operatorname{tr}(
\bQ_{j -}),\qquad d = \operatorname{tr}(\bI_d) = \operatorname{tr}(\bQ_{j +}) + \operatorname{tr}(
\bQ_{j -}),
\end{eqnarray*}
and solving the equations we get
%
%e62 #&#
%
\begin{equation}
\label{tr-BQ} \operatorname{tr}(\bQ_{j \pm}) = d/2,\qquad \operatorname{tr}(\bB_{j}
\bQ_{j \pm}) = \pm \operatorname{tr}(\bQ_{j \pm}) = \pm d/2,\qquad j \neq1.
\end{equation}
For $j \neq j^\prime$, $\bB_{j}$ and $\bB_{j^\prime}$ are orthogonal,
\[
0=\operatorname{tr}(\bB_{j^\prime} \bB_{j}) = \operatorname{tr}(\bB_{j^\prime}
\bQ_{j+}) - \operatorname{tr}(\bB_{j^\prime} \bQ_{j -})
\]
and further if $j, j^\prime\neq1$,
\begin{eqnarray*}
\bB_{j^\prime} \bQ_{j+} + \bB_{j^\prime}
\bQ_{j -} &=& \bB_{j^\prime} (\bQ_{j+} +
\bQ_{j -}) = \bB_{j^\prime},\\
\operatorname{tr}(\bB_{j^\prime}
\bQ_{j+}) + \operatorname{tr}(\bB_{j^\prime} \bQ_{j -}) &=& \operatorname{tr}(
\bB_{j^\prime})=0,
\end{eqnarray*}
which imply
%
%e63 #&#
%
\begin{equation}
\label{tr-BQ1} \operatorname{tr}(\bB_{j^\prime} \bQ_{j \pm}) = 0,\qquad j \neq
j^\prime, j, j^\prime\neq1.
\end{equation}

For any density matrix $\brho$ with representation (\ref{sparse})
under the Pauli basis (\ref{Pauli-basis}), we have
$1= \operatorname{tr}(\brho)= \alpha_1 \operatorname{tr}(\bB_1) = d \alpha_1$ and hence $\alpha
_1= 1/d$.
Consider special density matrices $\brho\in\Theta_{s}$ with expression
%
%e64 #&#
%
\begin{equation}
\label{sparse1} \brho= \frac{1}{d} \bI_d + \frac{\beta}{d}
\bB_{j^*},
\end{equation}
where $\beta$ is a real number with $|\beta|<1$, and index $j^* \neq
1$.\eject % will be determined later.
%has at most $s$ terms in representation (\ref{sparse}).

To check if $|\cI_j(\brho)| \geq2$, we need to evaluate $\operatorname{tr}(\brho
\bQ_{j \pm})$ for $\brho$ given by (\ref{sparse1}), $j=1,\ldots, p$.

For $j=1$, $\bB_1=\bQ_{1 +}=\bI_d$, and since $\operatorname{tr}(\bB_{j^*})=0$, we have
%
%e65 #&#
%
\begin{equation}
\label{tr-Qrho1} \operatorname{tr}( \brho\bQ_{1 +}) = \frac{1}{d} \operatorname{tr}(
\bI_d) + \frac{\beta}{d} \operatorname{tr}(\bB_{j^*}) = 1.
\end{equation}
For $j=j^*$, from (\ref{tr-BQ}) we have $\operatorname{tr}(\bQ_{j^* \pm})=d/2$ and
$\operatorname{tr}(\bB_{j^*} \bQ_{j^* \pm}) = \pm d/2$, and thus
%
%e66 #&#
%
\begin{equation}
\label{tr-Qrho2} \operatorname{tr}( \brho\bQ_{j^* \pm}) = \frac{1}{d} \operatorname{tr}(
\bQ_{j^* \pm}) + \frac
{\beta}{d} \operatorname{tr}(\bB_{j^*}
\bQ_{j^* \pm} ) = \frac{1 \pm\beta}{2} \in(0,1).
\end{equation}
For $j \neq j^*$ or $1$ [i.e., $(\ell_1,\ldots, \ell_b) \neq(\ell
_1^*,\ldots, \ell_b^*)$
or $(0,\ldots, 0)$], from (\ref{tr-BQ1}) we have $\operatorname{tr}(\bB_{j^*} \bQ
_{j \pm}) = 0$, and thus
%
%e67 #&#
%
\begin{equation}
\label{tr-Qrho3} \operatorname{tr}(\brho\bQ_{j \pm}) = \frac{1}{d} \operatorname{tr}(
\bQ_{j \pm}) + \frac{\beta
}{d} \operatorname{tr}(\bB_{j^*}
\bQ_{j \pm}) = \frac{1}{d} \operatorname{tr}(\bQ_{j \pm}) =
\frac{1}{2}.
\end{equation}
Equations (\ref{tr-Qrho1})--(\ref{tr-Qrho3}) immediately show that for $\brho$
given by (\ref{sparse1}) and $j \neq1$,
$\operatorname{tr}( \brho\bQ_{j \pm}) \in[(1-|\beta|)/2, (1+|\beta|)/2]$, $|\cI
_j(\brho)|=2$, and
\[
\sum_{j=1}^p 1\bigl(\bigl|\cI_j(
\brho)\bigr|\geq2\bigr) = p - 1,
\]
which implies
\[
\max_{\brho\in\Theta^p_s} \Biggl\{ \frac{1}{p} \sum
_{j=1}^p 1\bigl(\bigl|\cI_j(\brho)\bigr|\geq2
\bigr) \Biggr\} \geq1 - \frac{1}{p}.
\]
\upqed
\end{pf*}

\begin{pf*}{Proof of Corollary \ref{cor-lowrank-Pauli}}
We use the notation and facts about the Pauli basis
(\ref{Pauli-basis}) in the proof of Corollary \ref{cor-sparse-Pauli}:
$p=d^2$, $d=2^b$, and we identify index $j =1,\ldots, p$ with $(\ell
_1, \ell_2,\ldots, \ell_b) \in\{0, 1, 2, 3\}^b$. Consider $\bB_j =
\bsigma_{\ell_1} \otimes\bsigma_{\ell_2}
\otimes\cdots\otimes\bsigma_{\ell_b}$. For $j=1$ [or $\ell_1= \cdots=
\ell_b=0$], $\bB_1=\bI_d$, and for $j \neq1$ [or $(\ell_1,\ldots,
\ell_b) \neq(0,\ldots, 0)$], $\bB_j$ has eigenvalues $\pm1$, $\bQ_{j
\pm}$ are the projections onto the eigen-spaces corresponding to
eigenvalues $\pm1$, respectively, $\bB_{j} = \bQ_{j +} - \bQ_{j -}$,
and $\bI_d=\bQ_{j +} + \bQ_{j -}$.

Let
\[
\be= \sqrt{2/7} \bigl[(\sqrt{3}/2, 1/2)^\prime+ (\sqrt{3}/2, \sqrt
{-1}/2)^\prime\bigr]
%[(\sqrt{3}/2, 1/2)^\prime/2 + (\sqrt{3}/2, \sqrt{-1}/2)^\prime/2
%]/sqrt{7/8}
= (\sqrt{6/7}, \sqrt{1/14} + \sqrt{-1/14})^\prime.
\]
Then for $\ell=0,1,2,3$, $\varpi_\ell=\be^\dagger\bsigma_{\ell} \be$ is
equal to $1$, $ 2\sqrt{3}/7$, $2 \sqrt{3}/7$ and $5/7$, respectively.
Let $U = \be^{\otimes b}$ and $\brho= U U^\dagger$. Then $\brho$ is a
rank one density matrix, and
\begin{eqnarray*}
\operatorname{tr}(\brho\bQ_{j +}) + \operatorname{tr}(\brho\bQ_{j -}) &=& \operatorname{tr}(\brho) = 1,
\\
\operatorname{tr}(\brho\bQ_{j +}) - \operatorname{tr}(\brho\bQ_{j -}) &=& \operatorname{tr}(\brho
\bB_j) %=
= U^\dagger\bB_j U
=
\bigl(\be^\dagger\bsigma_{\ell_1} \be\bigr) \times\cdots\times
\bigl(\be^\dagger\bsigma_{\ell_b} \be\bigr) \\
&=& \varpi_{\ell_1}
\cdots\varpi_{\ell_b}.
\end{eqnarray*}
Solving the two equations we obtain $\operatorname{tr}(\brho\bQ_{j \pm}) = (1 \pm
\varpi_{\ell_1} \cdots\varpi_{\ell_b})/2$.\eject

For $j \neq1$ [or $(\ell_1,\ldots, \ell_b) \neq(0,\ldots, 0)$],
$(\varpi_{\ell_1},\ldots, \varpi_{\ell_b}) \neq(1,\ldots, 1)$, and
$ 0 \leq\varpi_{\ell_1} \cdots\varpi_{\ell_b} \leq5/7$, and thus
$\operatorname{tr}(\brho\bQ_{j +}) \geq1/2$ and $\operatorname{tr}(\brho\bQ_{j -}) \geq1/7$,
which immediately shows that for the given rank one density matrix
$\brho$ and $j \neq1$, $|\cI_j(\brho)|=2$, and
\[
\sum_{j=1}^p 1\bigl(\bigl|\cI_j(
\brho)\bigr|\geq2\bigr) = p - 1,
\]
which implies
\[
\max_{\brho\in\Theta_r} \Biggl\{ \frac{1}{p} \sum
_{j=1}^p 1\bigl(\bigl|\cI_j(\brho)\bigr|\geq2
\bigr) \Biggr\} \geq1 - \frac{1}{p}.
\]
\upqed
\end{pf*}

\begin{pf*}{Proof of Corollary \ref{cor-lowrank}}
Since under $\bg_1,\ldots, \bg_d$, basis matrices $\bB_j$ defined in
the corollary have
the same behavior as matrix basis (\ref{Hermitian-basis}) under $\be
_1,\ldots, \be_d$,
from the proof of Corollary \ref{cor1} on the eigen-structures of
matrix basis
(\ref{Hermitian-basis}) we see that under
$\bg_1,\ldots, \bg_d$, $\bB_j$ has possible eigenvalues $0$ and
$1$ for diagonal $\bB_j$ and eigenvalues $0$, $1$ and $-1$ for
nondiagonal $\bB_j$. For the diagonal case,
corresponding to eigenvalue $1$, the eigenvector is $\bg_\ell$; for
the real symmetric nondiagonal case, corresponding to eigenvalues $\pm1$,
the eigenvectors are $(\bg_{\ell_1} \pm\bg_{\ell_2})/\sqrt{2}$,
respectively; and for the complex Hermitian nondiagonal case,
corresponding to eigenvalue $\pm1$,
the eigenvectors are $ (\bg_{\ell_1} \pm\sqrt{-1} \bg_{\ell
_2})/\sqrt{2}$, respectively. Denote by $\bQ_{j0}$, $\bQ_{j1}$ and
$\bQ_{j,-1}$ the projections on the eigen-spaces
corresponding to eigenvalues $0$, $1$ and $-1$, respectively.

For diagonal $\bB_j$ with $j$ corresponding to $(\ell, \ell)$, it is
a binomial case,
\[
\operatorname{tr}(\brho\bQ_{j0}) = 1 - \operatorname{tr}(\brho\bQ_{j1}),
\qquad\operatorname{tr}(\brho\bQ
_{j1}) = \bg_{\ell}^\prime\brho\bg_{\ell} =
\sum_{a=1}^r \xi_a
\bigl|U_a^\dagger\bg_{\ell}\bigr|^2
\]
and
\begin{eqnarray*}
\bigl|\cI_j(\brho)\bigr| &=& 2\cdot1\bigl(0< \operatorname{tr}(\brho\bQ_{j1}) <1
\bigr) + 1\bigl(\operatorname{tr}(\brho\bQ_{j1})=1\bigr) \\
&&{}+ 1\bigl(\operatorname{tr}(\brho
\bQ_{j1})=0\bigr).
\end{eqnarray*}
In order to have $|\cI_j(\brho)|\geq2$ possible, we need $\operatorname{tr}(\brho
\bQ_{j1}) \in(0,1)$.
Since $\brho$ is generated by at most $r$ vectors $U_a$, and for each
$U_a$ there are at most $2 \gamma$ of $\bg_\ell$ with $U_a^\dagger
\bg_{\ell} \neq0$,
among all the $d$ diagonal matrices $\bB_j$ there are at most $ 2 r
\gamma$ of diagonal matrices $\bB_j$ for which
it is possible to have $\operatorname{tr}(\brho\bQ_{j1}) \in(0,1)$ and thus $|\cI
_j(\brho)|\geq2$.

For nondiagonal $\bB_j$, it is a trinomial case,
\[
\operatorname{tr}(\brho\bQ_{j0}) = 1 - \operatorname{tr}(\brho\bQ_{j1}) - \operatorname{tr}(\brho
\bQ_{j,-1}),
\]
and $\operatorname{tr}(\brho\bQ_{j\pm1})$ depend on whether $\bB_j$ is real or complex.

For real symmetric nondiagonal $\bB_j$ with $j$ corresponding to
$(\ell_1, \ell_2)$,
\[
\operatorname{tr}(\brho\bQ_{j\pm1}) = (\bg_{\ell_1} \pm\bg_{\ell_2})^\prime
\brho(\bg_{\ell_1} \pm\bg_{\ell_2})/2 = \sum
_{a=1}^r \xi_a \bigl|U_a^\dagger(
\bg_{\ell_1} \pm\bg_{\ell_2})\bigr|^2/2;
\]
and for imaginary Hermitian nondiagonal $\bB_j$ with $j$
corresponding to $(\ell_1, \ell_2)$,
\begin{eqnarray*}
\operatorname{tr}(\brho\bQ_{j\pm1}) &=& (\bg_{\ell_1} \pm\sqrt{-1}
\bg_{\ell
_2})^\dagger\brho(\bg_{\ell_1} \pm\sqrt{-1}
\bg_{\ell_2})/2 \\
&=& \sum_{a=1}^r
\xi_a \bigl|U_a^\dagger(\bg_{\ell_1} \pm
\sqrt{-1} \bg_{\ell_2})\bigr|^2/2.
\end{eqnarray*}
In order to have $|\cI_j(\brho)|\geq2$ possible, we need $\operatorname{tr}(\brho
\bQ_{j1}) \in(0,1)$ or $\operatorname{tr}(\brho\bQ_{j -1}) \in(0,1)$.
Since $\brho$ is generated by at most $r$ vectors $U_a$, and for each
$U_a$ there are at most $2 \gamma$ number of $\bg_\ell$ with
$U_a^\dagger\bg_{\ell} \neq0$,
among $(d^2-d)/2$ real symmetric nondiagonal matrices $\bB_j$ [or
$(d^2-d)/2$ imaginary Hermitian nondiagonal matrices $\bB_j$],
there are at most $4 r \gamma^2$ of real symmetric nondiagonal
$\bB_j$ (or imaginary Hermitian nondiagonal matrices $\bB_j$)
for which it is possible to have $\operatorname{tr}(\brho\bQ_{j 1}) \in(0,1)$ or
$\operatorname{tr}(\brho\bQ_{j, -1}) \in(0,1)$ and thus $|\cI_j(\brho)|\geq2$.

Finally, for $\brho\in\Theta_{r\gamma}$, combining the results on
the number of $\bB_j$ for which it is possible to have $|\cI_j(\brho
)|\geq2$ in the diagonal, real symmetric and imaginary Hermitian
cases, we conclude
\[
\sum_{j=1}^p 1\bigl(\bigl|\cI_j(
\brho)\bigr|\geq2\bigr) \leq8 r \gamma^2 + 2 r \gamma,
\]
%
%4 r \gamma(d-1) + 2 r \gamma\leq4 r \gamma d, \]
and
\[
\zeta_p = \max_{\brho\in\Theta_{r\gamma}} \Biggl\{\frac{1}{p}
\sum_{j=1}^p 1\bigl(\bigl|\cI_j(
\brho)\bigr|\geq2\bigr) \Biggr\} \leq\frac{ 2 r
\gamma(4 \gamma+1)}{p}.
\]
\upqed
\end{pf*}

%denote by $\cR(\brho)$ the range of $\brho$ and $\cR^\bot(\brho)$ the
%orthogonal complement of $\cR(\brho)$. Then $\mathbb{C}^{d\times d} =
%Since the range of each $\brho\in\Theta_r$ is generated by at most
%$r$ number of $\bg_1,\ldots, \bg_d$, among all eigen-vectors of $
%the number of the eigen vectors falling inside $\cR(\rho)$ is at most
%$r$, and the number of the eigen vectors falling inside $\cR^\bot(
%there are at most $r$ number of $\bB_j$ whose eigen vectors
%corresponding to eigenvalues $\pm1$ have nonzero projections on $\cR(
%we have $|\cI_j(\brho)| \geq2$ only for $r$ number of $\bB_j$, which
%implies
%&& %\max_{\brho\in\Theta_r} \left\{ \sum_{j=1}^p 1(|\cI_j(\brho)|
% \sum_{j=1}^p 1(|\cI_j(\brho)|\geq2) \leq r, \\
%&& \zeta_p = \max_{\brho\in\Theta} \left\{\frac{1}{p} \sum_{j=1}^p
%1(|\cI_j(\brho)|\geq2) \right\} \leq\frac{r}{p}.

%sA #&#
%
\begin{appendix}\label{app}
\section*{Appendix: Proofs of Lemmas 2--4}\label{sec8}

%Suppose that $P$ is a multinomial distribution ${\cal M}(m, \theta_1,
%,\ldots, \theta_{r})$, where $r \geq2$ is a fixed integer,
%,\ldots, \theta_{r}) \leq\max(\theta_1,\ldots, \theta_{r}) \leq c_1,
%and $0< c_0 \leq c_1<1$ are two fixed constants. Denote by $Q$ the
%multivariate normal distribution whose mean and covariance
%are the same as $P$. Let $P_k^{*}$ be the convolution of distribution
%$P_k$ and the distribution of $(\psi_1,\ldots, \psi_r)$, where
%$\psi_1,\ldots, \psi_{r-1}$ are independent and follow a uniform
%distribution on $(-1/2, 1/2)$, and $\psi_r = - \psi_1 - \cdots-

\begin{pf*}{Proof of Lemma \ref{lemMulti-norm}}
For $r=2$, it is the binomial case, and the lemma is a consequence of
(\ref{TV-KL1}) and Lemma \ref{lemBin-norm}.

For $r=3$, write $(U_1, U_2, U_3) \sim P$ and $(V_1, V_2, V_3) \sim Q$.
Add independent uniforms on $(-1/2, 1/2)$ to $U_1$ and $U_2$, denote
the resulting random variables by $U_1^*$ and $U_2^*$, respectively,
and let $U_3^*= m - U_1^* - U_2^*$. Then
$(U_1^*, U_2^*, U_3^*) \sim P^*$. Note that
$U_1 + U_2 + U_3=U_1^* + U_2^* + U_3^* = V_1 + V_2 + V_3=m$, and $U_1$
and $U_2$ are equal to the round-offs, $[U_1^*]$ and $[U_2^*]$,
of $U_1^*$ and $U_2^*$, respectively, here round-off $[x]$ means
rounding $x$ off to the nearest integer.

For trinomial random variable $(U_1, U_2, U_3) \sim{\cal M}(m, \theta
_1, \theta_2, \theta_3)$, we have
$U_1 \sim \operatorname{Bin}(m, \beta_1)=P_1$, the conditional distribution of $U_2$
given $U_1$: $U_2|U_1 \sim \operatorname{Bin}(m - U_1, \beta_2)=P_2$, and
$U_3 = m - U_1 - U_2$, where $\beta_1=\theta_1$, $\beta_2 = \theta
_2/(\theta_2+\theta_3)$, $\beta_3 = \theta_3/(\theta_2+\theta_3)$.
Since $\theta_j$ are between $c_0$ and $c_1$, $\beta_2$ and $\beta
_3$ are between $c_0/(c_0+c_1)$ and $c_1/(c_0+c_1)$.
We have decomposition $P=P_1 P_2$.

Denote by $P_1^*$ the distribution of $U_1^*$ and $P_2^*$ the
conditional distribution of $U_2^*$ given $U_1^*$. Then $P_1^*$ is the
convolution
of $P_1$ and an independent uniform distribution on $(-1/2,1/2)$. Since
the added uniforms are independent of $U_j$, and
$U_j$ is the round-off of $U_j^*$, the conditional distribution of
$U_2^*$ given $U_1^*$ is equal to the conditional distribution of
$U_2^*$ given
$U_1=[U_1^*]$, which in turn is equal to the convolution of $P_2$ and
an independent uniform distribution on $(-1/2,1/2)$. We have
decomposition $P^*=P^*_1 P^*_2$.

For trivariate normal random variable $(V_1, V_2, V_3) \sim Q$, we have
$V_1 \sim N(m \beta_1, m \beta_1 (1-\beta_1))=Q_1$,
the conditional distribution of $V_2$ given $V_1$: $V_2|V_1 \sim
N((m-V_1)\beta_2, m (1-\beta_1) \beta_2 \beta_3)=Q_2$, and
$V_3=m - V_1 - V_2$. We have decomposition $Q=Q_1 Q_2$.

As there is a difference in conditional variance between $P_2$ and
$Q_2$, we define $V_2^\prime\sim
Q_2^\prime=N((m-V_1)\beta_2, (m-V_1) \beta_2 \beta_3)$ to match the
conditional variance of~$P_2$, and $V_3^\prime= m - V_1 - V_2^\prime$.
Simple direct calculations show that given $V_1$,
%
%eA.1 #&#
%
\begin{equation}
\label{HellingerQQprime} H^2\bigl(Q_2,
Q_2^\prime\bigr) \leq\frac{3}{2} \biggl( 1 -
\frac{m-V_1}{m
(1-\beta_1)} \biggr)^2.
\end{equation}
Note that $P^*=P_1^* P_2^*$ and $Q=Q_1 Q_2$ are probability measures on
$\{(x_1, x_2,\break x_3)\dvtx x_1 + x_2 + x_3=m\}$.
Define probability measures $Q_1 Q_2^\prime$ and $P_1^* Q_2^\prime$
on $\{(x_1, x_2, x_3)\dvtx x_1 + x_2 + x_3=m\}$, where
$Q_1$ and $P_1^*$ are their respective marginal distributions of the
first component, and $Q_2^\prime$ is their conditional
distribution of the second component given the first component. We use
$Q_1 Q_2^\prime$ and $P_1^* Q_2^\prime$ to bridge between
$P^*=P_1^* P_2^*$ and $Q=Q_1 Q_2$. Applying triangle inequality we obtain
%
%eA.2 #&#
%
\begin{eqnarray}
\label{P*Q}\quad H\bigl(P^*, Q\bigr) &\leq& H\bigl(P^*, Q_1
Q_2^\prime\bigr) + H\bigl(Q_1
Q_2^\prime, Q\bigr)
\nonumber\\
&\leq& H\bigl(P_1^* P_2^*, P_1^*Q_2^\prime
\bigr) + H\bigl(P_1^* Q_2^\prime,
Q_1 Q_2^\prime\bigr) \\
&&{}+ H\bigl(Q_1
Q_2^\prime, Q_1 Q_2\bigr).\nonumber
\end{eqnarray}
Using (\ref{Hellinger0}), (\ref{TV-KL1}), Lemma \ref{lemBin-norm}
and (\ref{HellingerQQprime}) we evaluate the Hellinger distances on
the right-hand side of (\ref{P*Q}) as follows:
%
%eA.3 #&#
%
\begin{eqnarray}\label{QprimeQ}
H^2\bigl(Q_1 Q_2^\prime,
Q_1 Q_2\bigr) &=& \int\biggl\llvert\sqrt{
\frac
{dQ_1}{dx_1} \frac{dQ_2}{dx_2}} - \sqrt{\frac{dQ_1}{dx_1}
\frac
{dQ_2^\prime}{dx_2} } \biggr\rrvert^2 \,dx_1
\,dx_2
\nonumber
\\
&=& \int dQ_1 \int\biggl\llvert\sqrt{\frac{dQ_2}{dx_2}} - \sqrt
{\frac
{dQ_2^\prime}{dx_2} } \biggr\rrvert^2 \,dx_2
\nonumber
\\
&=& E_{Q_1} \bigl[ H^2\bigl(Q_2,
Q_2^\prime\bigr)\bigr]
\\
&\leq& E_{Q_1} \biggl[\frac{3}{2} \biggl(1 -
\frac{m-V_1}{m
(1-\beta_1)} \biggr)^2 \biggr]
\nonumber
\\
& %= \frac{3}{2} \frac{m \beta_1 (1-\beta_1)}{ m^2 (1-\beta_1)^2}
=& \frac{3\beta_1}{2m(1-\beta_1)} \leq\frac{3 \theta_1 }{2 m
(\theta_2+\theta_3)} \leq
\frac{C}{m},\nonumber
\end{eqnarray}
where (\ref{HellingerQQprime}) is used to bound $H^2(Q_2, Q_2^\prime
)$ and obtain the first inequality
%
%eA.4 #&#
%
\begin{eqnarray}\label{P*Q1}
H^2\bigl(P_1^* Q_2^\prime,
Q_1 Q_2^\prime\bigr) &=& \int\biggl\llvert\sqrt
{\frac
{dP_1^*}{dx_1}} - \sqrt{\frac{dQ_1}{dx_1} } \biggr\rrvert
^2 \,dx_1 \int dQ_2^\prime
\nonumber
\\
&=& \int\biggl\llvert\sqrt{\frac{dP_1^*}{dx_1}} - \sqrt{
\frac
{dQ_1}{dx_1} } \biggr\rrvert^2 \,dx_1
= H^2\bigl(P_1^*, Q_1\bigr)
\\
&\leq&\exp
\bigl(- C m^{1/3}\bigr) + \frac{C}{m \theta_1
(1-\theta_1)} \leq\frac{C}{m},\nonumber
\end{eqnarray}
where Lemma \ref{lemBin-norm} and (\ref{TV-KL1}) are used to bound
$H^2(P_1^*, Q_1)$ and obtain the first inequality
%
%eA.5 #&#
%
\begin{eqnarray}\label{P*Q2}
H^2\bigl(P_1^*P_2^*, P_1^*
Q_2^\prime\bigr) &=& \int dP_1^* \int\biggl\llvert
\sqrt{\frac{dP_2^*}{dx_2}} - \sqrt{\frac{dQ_2^\prime}{dx_2} } \biggr
\rrvert^2 \,dx_2
\nonumber
\\
&=& E_{P_1^*} \bigl[ H^2\bigl(P_2^*,
Q_2^\prime\bigr)\bigr]
\nonumber\\[-8pt]\\[-8pt]
&\leq& 2 - 2 E_{P^*_1} \biggl\{ E_{P^*_2} \biggl[
1_A \sqrt{\frac{dP^*_2}{dQ_2^\prime}} \Big|U_1 \biggr]
\biggr\}
\nonumber
\\
&\leq& 2 P^*\bigl(A^c\bigr) + E_{P^*_1} \biggl\{
1_{A_1} E_{P^*_2} \biggl[1_{A_2} \log
\frac{dP^*_2}{dQ_2^\prime} \Big| U_1 \biggr] \biggr\},\nonumber
\end{eqnarray}
where\vspace*{-1pt} we use (\ref{TV-KL1}) to bound $H^2(P_2^*, Q_2^\prime)$ and
obtain the last two inequalities,
%We define $A$ and bound $P^*(A^c)$ and $E_{P^*} \left[ 1_{A} \log
$A = A_1 \cap A_2$, and
\begin{eqnarray*}
A_1 &=& \bigl\{ |U_1 - m \beta_1| \leq
\bigl[m \beta_1 (1-\beta_1)\bigr]^{2/3} \bigr
\},\\
A_2 &=& \bigl\{ \bigl|U_2 - (m-U_1)
\beta_2\bigr| \leq\bigl[(m-U_1) \beta_2 (1-\beta
_2)\bigr]^{2/3} \bigr\}.
\end{eqnarray*}
%
%A_2 = \{ |U_2 - (m-U_1) \beta_2| \leq[(m-U_1) \beta_2 (1-
We evaluate $P^*(A^c)$ as follows:
%
%eA.6 #&#
%
\begin{eqnarray}
\label{PAr=3}
P^*\bigl(A^c\bigr) &=& P\bigl(A_1^c
\cup\bigl[A_2^c \cap A_1\bigr] \bigr)= P
\bigl(A_1^c\bigr) + P\bigl(A_2^c
\cap A_1\bigr)\nonumber\\
&=& P_1\bigl(A_1^c
\bigr) + E_{P}\bigl[1_{A_1} P\bigl(A_2^c|U_1
\bigr)\bigr]
\nonumber
\\
&\leq& \exp\bigl(- C m^{1/3}\bigr) + E_{P}\bigl[
1_{A_1} \exp\bigl(- C \{m-U_1\}^{1/3} \bigr)
\bigr]
\\
&\leq& \exp\bigl(- C m^{1/3}\bigr) + \exp\bigl(- C \bigl\{m- m
\beta_1 - \bigl[m \beta_1 (1-\beta_1)
\bigr]^{2/3}\bigr\}^{1/3} \bigr) \nonumber\\
&\leq&2 \exp\bigl(- C
m^{1/3}\bigr),\nonumber
\end{eqnarray}
where we utilize\vspace*{1pt} Lemma \ref{lemBin-norm} to derive $P_1(A_1^c)$ and
$P(A_2^c|U_1)$, and bound $m - U_1$ by using the fact that on $A_1$,
$U_1 \leq m \beta_1 + [m \beta_1 (1-\beta_1)]^{2/3}$.
Again we apply Lemma~\ref{lemBin-norm} to bound $E_{P^*_2} [
1_{A_2} \log\frac{dP^*_2}{dQ_2^\prime} | U_1 ]$
and obtain
%
%eA.7 #&#
%
\begin{eqnarray}
\label{EAr=3} %&& E_{P^*} \left[ 1_{A} \log\frac{dP^*}{dQ} \right] =
%E_{P^*} \left[
%1_{A_1 \cap A_2} \log\frac{dP^*_1}{dQ_1} + 1_{A_1 \cap A_2} \log
%&&\qquad \leq E_{P^*} \left[ 1_{A_1} \log\frac{dP^*_1}{dQ_1} \right] +
%E_{P^*} \left\{1_{A_1} E_{P^*}\left[ 1_{A_2} \log
%&& E_{P^*} \left[ 1_{A_1} \log\frac{dP^*_1}{dQ_1} \right] \leq
&&E_{P^*_1} \biggl\{
1_{A_1} E_{P^*_2} \biggl[ 1_{A_2} \log
\frac{dP^*_2}{dQ_2^\prime} \Big| U_1 \biggr] \biggr\}\nonumber \\
&&\qquad\leq E_{P^*_1}
\biggl\{1_{A_1} \frac{C}{ (m-U_1) \beta_2(1-\beta_2)} \biggr\}
\\
&&\qquad\leq\frac{C}{ (m- m \beta_1 - [m \beta_1 (1-\beta_1)]^{2/3} )
\beta_2(1-\beta_2)} \leq\frac{C}{m},\nonumber
\end{eqnarray}
where to bound $1/(m-U_1)$ we use the fact that on $A_1$, $U_1 \leq m
\beta_1 + [m \beta_1 (1-\beta_1)]^{2/3}$.

Substituting (\ref{PAr=3}) and (\ref{EAr=3}) into (\ref{P*Q2}) and
then combining it with (\ref{P*Q})--(\ref{P*Q1})
%(\ref{P*Q}), (\ref{QprimeQ}), (\ref{P*Q1}) and (\ref{P*Q2})
we prove that the lemma is true for $r=3$.

Consider the $r+1$ case.
Write $(U_1,\ldots, U_r, U_{r+1}) \sim P$, $U_1 + \cdots+
U_{r+1}=m$, and decompose
$P=P_1 P_2 \cdots P_{r-1} P_r$, where $U_1 \sim P_1=\operatorname{Bin}(m, \beta_1)$,
$P_j = \operatorname{Bin}(m-T_{j-1}, \beta_j)$ is the conditional distribution of
$U_j$ given $U_1,\ldots, U_{j-1}$,
$T_{j}=U_1 + \cdots+ U_{j}$, $\beta_1=\theta_1$, $\beta_j=\theta
_j/(1-\theta_1 - \cdots-\theta_{j-1})$.
Since $\theta_j$ are between $c_0$ and $c_1$, all $\beta_j$ are
between $c_0/(c_0+ r c_1)$ and $c_1/(c_0+c_1)$
that are bounded away from $0$ and $1$.

Similarly write $(V_1,\ldots, V_r, V_{r+1}) \sim Q$, $V_1 + \cdots+
V_{r+1}=m$, and decompose
$Q=Q_1 Q_2 \cdots Q_{r-1} Q_r$, where $V_1 \sim Q_1=N(m\beta_1, m
\beta_1(1-\beta_1))$, and $Q_j=N((m-S_{j-1})\beta_j,
m (\theta_j +\cdots+\theta_{r+1}) \beta_j(1-\beta_j))$ is the
conditional distribution of $V_j$ given $V_1,\ldots, V_{j-1}$,
where $S_j=V_1 + \cdots+ V_j$.

As there are differences in conditional variance between $P_j$ and
$Q_j$, we handle the differences by introducing
$Q_j^\prime\cdots Q_r^\prime$ as follows. Given $V_1,\ldots,
V_{j-1}$ we define
$(V_j^\prime,\ldots, V_r^\prime, V_{r+1}^\prime) \sim Q_j^\prime
\cdots Q_r^\prime$, where the conditional distribution of $V_\ell
^\prime$
given $V_1,\ldots, V_{j-1}, V_j^\prime,\ldots, V^\prime_{\ell
-1}$ is
$Q_\ell^\prime=N((m-S_{\ell-1}^\prime)\beta_\ell, (m - S_{\ell
-1}^\prime) \beta_\ell(1-\beta_\ell))$ for $\ell= j,\ldots, r$,
$V_{r+1}^\prime=m - V_1 -\cdots- V_{j-1} - V^\prime_j - \cdots-
V_r^\prime$, and
$S_{\ell}^\prime=V_1 + \cdots+ V_{j-1} + V_j^\prime+ \cdots+
V_\ell^\prime$. Then given $V_1,\ldots, V_{j-1}$,
%
%eA.8 #&#
%
\begin{equation}
\label{rHellingerQQprime} H^2\bigl(Q_j,
Q_j^\prime\bigr) \leq\frac{3}{2} \biggl( 1 -
\frac{m-S_{j-1}}{m
(\theta_j +\cdots+\theta_{r+1})} \biggr)^2.
\end{equation}
Add independent uniforms on $(-1/2, 1/2)$ to $U_1,\ldots, U_r$, denote
the resulting corresponding random variables by $U_j^*$, and let
$U_{r+1}^*= m - U_1^* - \cdots- U_r^*$. Then
$(U_1^*,\ldots, U_{r+1}^*) \sim P^*$. Note that
$U_1 + \cdots+ U_{r+1}=U_1^* + \cdots+ U_{r+1}^* = V_1 + \cdots+
V_{r+1}=m$, and $U_j$ is equal to the round-off of $U_j^*$.
Let $P^*=P_1^* P_2^* \cdots P_{r-1}^* P_r^*$, where we denote by
$P_1^*$ the distribution of $U_1^*$ and $P_j^*$ the conditional
distribution of $U_j^*$ given $U_1^*,\ldots, U_{j-1}^*$. Then $P_1^*$
is the convolution\vspace*{2pt} of $P_1$ and an independent uniform
distribution on $(-1/2,1/2)$. Since the added uniforms are independent
of $U_j$, and
$U_j$ is the round-off of $U_j^*$, the conditional distribution of
$U_j^*$ given $U_1^*,\ldots, U_{j-1}^*$ is equal to the
conditional distribution of $U_j^*$ given $U_1=[U_1^*],\ldots,
U_{j-1}=[U^*_{j-1}]$, which in turn is equal to the convolution
of $P_j$ and an independent uniform distribution on $(-1/2,1/2)$.

Note that $P^*=P_1^* \cdots P_r^*$ and $Q=Q_1\cdots Q_r$ are
probability measures on $\{(x_1,\ldots, x_r, x_{r+1})\dvtx x_1 +
\cdots+ x_{r+1}=m\}$. We define probability measures $Q_1 \cdots Q_j
Q_{j+1}^\prime\cdots Q_r^\prime$ and $P^*_1 \cdots P^*_{j-1}
Q_{j}^\prime\cdots Q_r^\prime$ on $\{(x_1,\ldots, x_r, x_{r+1})\dvtx
x_1 + \cdots+ x_{r+1}=m\}$, $j=2,\ldots, r$, and use them to bridge
between $P^*$ and $Q$. Applying triangle inequality, we have
%
%eA.9 #&#
%
\begin{eqnarray}\label{rP*Q1}
H\bigl(P^*, Q\bigr) &\leq& H\bigl(P^*, Q_1\cdots Q_{r-1}
Q_r^\prime\bigr) + H\bigl(Q_1\cdots
Q_{r-1} Q_r^\prime, Q\bigr)
\nonumber
\\
&\leq& H\bigl(P^*, Q_1\cdots Q_{r-2}
Q_{r-1}^\prime Q_r^\prime\bigr) \nonumber\\
&&{}+ H
\bigl(Q_1\cdots Q_{r-2} Q_{r-1}^\prime
Q_r^\prime, Q_1\cdots Q_{r-1}
Q_r^\prime\bigr) \nonumber\\
&&{}+ H\bigl(Q_1\cdots
Q_{r-1} Q_r^\prime, Q\bigr)
\leq \cdots\\
&\leq& H\bigl(P^*, Q_1 Q_2^\prime
\cdots Q_r^\prime\bigr) \nonumber\\
&&{}+ \sum_{j=2}^{r}
H\bigl(Q_1\cdots Q_{j-1} Q_j^\prime
\cdots Q_r^\prime, Q_1\cdots Q_j
Q_{j+1}^\prime\cdots Q_r^\prime\bigr)\nonumber
\end{eqnarray}
and
%
%eA.10 #&#
%
\begin{eqnarray}\label{rP*Q2}
&&
H\bigl(P^*, Q_1 Q_2^\prime\cdots
Q_r^\prime\bigr) \nonumber\\
&&\qquad\leq H\bigl(P^*, P_1^* \cdots
P_{r-1}^* Q_r^\prime\bigr) + H\bigl(P_1^*
\cdots P_{r-1}^* Q_r^\prime, Q_1
Q_2^\prime\cdots Q_r^\prime\bigr)
\nonumber
\\
&&\qquad\leq H\bigl(P^*, P_1^* \cdots P_{r-1}^*
Q_r^\prime\bigr) + H\bigl(P_1^* \cdots
P_{r-1}^* Q_r^\prime, P_1^* \cdots
P_{r-2}^*Q_{r-1}^\prime Q_r^\prime
\bigr)
\\
&&\qquad\quad{} + H\bigl(P_1^* \cdots P_{r-2}^* Q_{r-1}^\prime
Q_r^\prime, Q_1 Q_2^\prime
\cdots Q_r^\prime\bigr)
\nonumber
\\
&&\qquad\leq \cdots\leq\sum_{j=1}^r H
\bigl(P_1^* \cdots P_{j}^* Q_{j+1}^\prime
\cdots Q_r^\prime, P_1^* \cdots
P_{j-1}^*Q_{j}^\prime\cdots Q_r^\prime
\bigr).\nonumber
\end{eqnarray}
Substitute (\ref{rP*Q2}) into (\ref{rP*Q1}) to get
%
%eA.11 #&#
%
\begin{eqnarray}\label{rP*Q3}
H\bigl(P^*, Q\bigr) &\leq& \sum_{j=1}^r H
\bigl(P_1^* \cdots P_{j}^* Q_{j+1}^\prime
\cdots Q_r^\prime, P_1^* \cdots
P_{j-1}^*Q_{j}^\prime\cdots Q_r^\prime
\bigr)
\nonumber\\[-8pt]\\[-8pt]
&&{} + \sum_{j=2}^{r} H
\bigl(Q_1\cdots Q_{j-1} Q_j^\prime\cdots
Q_r^\prime, Q_1\cdots Q_j
Q_{j+1}^\prime\cdots Q_r^\prime\bigr).\nonumber
\end{eqnarray}
Using (\ref{Hellinger0}), (\ref{TV-KL1}), Lemma \ref{lemBin-norm}
and (\ref{rHellingerQQprime}) we evaluate the Hellinger distances on
the right-hand side of (\ref{rP*Q3}) as follows:
%
%eA.12 #&#
%
\begin{eqnarray}\label{rQQprime}
&&
H^2\bigl(Q_1\cdots Q_{j-1}
Q_{j}^\prime\cdots Q_r^\prime,
Q_1\cdots Q_j Q_{j+1}^\prime\cdots
Q_r^\prime\bigr)
\nonumber
\\
&&\qquad = \int dQ_1 \cdots dQ_{j-1} \int\biggl\llvert\sqrt{
\frac{dQ_j}{dx_j}} - \sqrt{\frac{dQ_j^\prime
}{dx_j} } \biggr\rrvert
^2 \,dx_j \int dQ_{j+1}^\prime\cdots
dQ_{r}^\prime
\nonumber
\\
&&\qquad = \int dQ_1 \cdots dQ_{j-1} \int\biggl\llvert\sqrt{
\frac{dQ_j}{dx_j}} - \sqrt{\frac{dQ_j^\prime}{dx_j} } \biggr\rrvert
^2 \,dx_j
\nonumber\\
&&\qquad = E_{Q_1\cdots Q_{j-1}} \bigl[ H^2\bigl(Q_j,
Q_j^\prime\bigr)\bigr]
\\
&&\qquad \leq E_{Q_1\cdots Q_{j-1}} \biggl[\frac{3}{2} \biggl(1 -
\frac
{m-S_{j-1}}{m (\theta_j+\cdots+\theta_{r+1})} \biggr)^2 \biggr]
\nonumber
\\
&&\qquad = %\frac{3}{2} \frac{m (\theta_j+\cdots+\theta_{r+1}) ( 1-\theta_j-
\frac{3 (1-\theta_j-\cdots-\theta_{r+1})}{2 m(\theta_j+\cdots
+\theta_{r+1})} \leq\frac{3 \theta_1}{ 2 m (\theta_r+\theta_{r+1})} \nonumber\\
&&\qquad\leq
\frac{C}{m},\nonumber
%&& H(P^*, Q) \leq H(P^*, Q^\prime) +
\end{eqnarray}
where we use (\ref{rHellingerQQprime}) to bound the Hellinger distance
$H^2(Q_j, Q_j^\prime)$ and obtain the first inequality
%
%eA.13 #&#
%
\begin{eqnarray}\label{rPQprime}\quad
&& H^2\bigl(P_1^* \cdots P_{j}^*
Q_{j+1}^\prime\cdots Q_r^\prime,
P_1^* \cdots P_{j-1}^*Q_{j}^\prime\cdots
Q_r^\prime\bigr)
\nonumber
\\
&&\qquad = \int dP_1^* \cdots dP_{j-1}^* \int\biggl\llvert\sqrt
{\frac
{dP_j^*}{dx_j}} - \sqrt{\frac{dQ_j^\prime}{dx_j} } \biggr\rrvert
^2 \,dx_j \int dQ_{j+1}^\prime\cdots
dQ_{r}^\prime
\nonumber
\\
&&\qquad = \int dP_1^* \cdots dP_{j-1}^* \int\biggl\llvert\sqrt
{\frac
{dP_j^*}{dx_j}} - \sqrt{\frac{dQ_j^\prime}{dx_j} } \biggr\rrvert
^2 \,dx_j
\nonumber
\\
&&\qquad = E_{P_1^*\cdots P_{j-1}^*} \bigl[ H^2\bigl(P_j^*,
Q_j^\prime\bigr)\bigr]
\nonumber
\\
&&\qquad \leq E_{P_1^*\cdots P_{j-1}^*} \biggl( 2 P^*_j\bigl(A_1^c
\cup\cdots\cup A_{j}^c|U_1,\ldots,
U_{j-1}\bigr) \\
&&\qquad\quad\hspace*{48pt}{}+ 1_{A_1 \cdots A_{j-1}} E_{P_j^*} \biggl[
1_{A_j} \log\frac{P_j^*}{Q_j^\prime} \Big|U_1,\ldots,
U_{j-1} \biggr] \biggr)
\nonumber
\\
&&\qquad = 2 P^*\bigl(A_1^c \cup\cdots\cup
A_{j}^c\bigr) \nonumber\\
&&\qquad\quad{}+ E_{P_1^*\cdots P_{j-1}^*} \biggl(
1_{A_1 \cdots A_{j-1}} E_{P_j^*} \biggl[ 1_{A_j} \log
\frac{P_j^*}{Q_j^\prime} \Big|U_1,\ldots, U_{j-1} \biggr] \biggr)
\nonumber
\\
&&\qquad \leq2 P^*\bigl(A_1^c \cup\cdots\cup
A_{j}^c\bigr) + E_{P_1^*\cdots
P_{j-1}^*} \biggl(1_{A_1 \cdots A_{j-1}}
\frac{C}{(m-T_{j-1})\beta
_j(1-\beta_j)} \biggr),\nonumber
\end{eqnarray}
where we use (\ref{TV-KL1}) to bound the Hellinger distance
$H^2(P^*_j, Q_j^\prime)$ and obtain the first inequality, we employ
Lemma \ref{lemBin-norm} to bound $E_{P^*_j} [ 1_{A_j}
\log\frac{dP^*_j}{dQ_j^\prime}|U_1,\ldots, U_{j-1}
]$ and get the last inequality,
and for $\ell=1,\ldots, j$,
\[
A_\ell= \bigl\{ \bigl|U_\ell- (m-U_1-
\cdots-U_{\ell-1}) \beta_\ell\bigr| \leq\bigl[(m-U_1-
\cdots- U_{\ell-1}) \beta_\ell(1-\beta_\ell)
\bigr]^{2/3} \bigr\}.
\]
Note that on $A_{j-1}$, $U_{j-1} \leq(m-T_{j-2}) \beta_{j-1} + [m
\beta_{j-1} (1-\beta_{j-1})]^{2/3}$. Then for $j=1,\ldots, r$ we
have on $A_1 \cdots A_{j-1}$,
%
%eA.14 #&#
%
\begin{eqnarray}
\label{m-T}
&&
m - T_{j-1} \nonumber\\
&&\qquad= m - T_{j-2} - U_{j-1}
\nonumber
\\
&&\qquad\geq (m - T_{j-2}) (1 - \beta_{j-1}) - \bigl[m
\beta_{j-1} (1-\beta_{j-1})\bigr]^{2/3}
\nonumber
\\
&&\qquad\geq (m - T_{j-3}) (1-\beta_{j-2}) (1 -
\beta_{j-1})\nonumber\\
&&\qquad\quad{} - (1-\beta_{j-1}) \bigl[m \beta_{j-2}
(1-\beta_{j-2})\bigr]^{2/3} - \bigl[m \beta_{j-1} (1-
\beta_{j-1})\bigr]^{2/3}
\geq \cdots\nonumber
\\
&&\qquad\geq m (1-\beta_1) \cdots(1-\beta_{j-1})\\
&&\qquad\quad{} -
m^{2/3} \sum_{\ell=1}^{j-1} \bigl[
\beta_{\ell} (1-\beta_{\ell})\bigr]^{2/3} (1-
\beta_\ell) \cdots(1-\beta_{j-1})
\nonumber
\\
&&\qquad\geq C m\nonumber
\end{eqnarray}
and thus
%
%eA.15 #&#
%
\begin{eqnarray}
\label{EAr} E_{P_1^*\cdots P_{j-1}^*} \biggl(1_{A_1 \cdots A_{j-1}}
\frac {C}{(m-T_{j-1})\beta_j(1-\beta_j)}
\biggr) &\leq&\frac{C}{m}. % = E_{P_1^*\cdots P_{j-1}^*} \left( \frac{C
%1_{A_1 \cdots A_{j-1}}
%}{(m-T_{j-2} - U_{j-1})\beta_j(1-\beta_j)} \right) \nonumber\\
%&&\qquad \leq E_{P_1^*\cdots P_{j-1}^*} \left(\frac{C 1_{A_1 \cdots
%A_{j-1}} }{(m-T_{j-2})(1-\beta_{j-1}) - [m \beta_{j-1} (1-
%&&\qquad \leq E_{P_1^*\cdots P_{j-2}^*} \left( \frac{C 1_{A_1 \cdots
%A_{j-2}} }{(m-T_{j-2})(1-\beta_{j-1}) - [m \beta_{j-1} (1-
%
%&&\qquad \leq\frac{C}{\{ m (1-\beta_1)\cdots(1-\beta_{j-1}) - m^{2/3} \sum_{
%%&&\qquad \leq\frac{C}{m (\theta_r + \theta_{r+1}) \beta_j(1-\beta_j)}
%&&\qquad \leq\frac{C}{m}.
\end{eqnarray}
We evaluate $P^*(A_1^c \cup\cdots\cup A_{j}^c)$ as follows:
%
%eA.16 #&#
%
\begin{eqnarray}
\label{PAc} \bigcup_{\ell=1}^j
A_\ell^c &=& \bigcup_{\ell=1}^j
\bigl(A^c_\ell A_{\ell-1}\cdots A_1
\bigr),
\nonumber
\\
P^* \Biggl(\bigcup_{\ell=1}^j
A_\ell^c \Biggr) &=& \sum_{\ell
=1}^j
P^*\bigl(A^c_\ell A_{\ell-1}\cdots
A_1 \bigr)
\nonumber
\\
&=& P^*\bigl(A_1^c\bigr)+ \sum
_{\ell=2}^j E_{P^*}\bigl[1_{A_1 \cdots A_{\ell-1}}
P^*\bigl(A_\ell^c|U_1,\ldots,
U_{\ell-1}\bigr)\bigr]
\\
&\leq& \exp\bigl[ - C m^{1/3} \bigr] +\sum
_{\ell=2}^j E_{P^*} \bigl(1_{A_1 \cdots A_{\ell-1}}
\exp\bigl[ - C (m - T_{\ell
-1})^{1/3} \bigr] \bigr)
\nonumber
\\
&\leq& \sum_{\ell=1}^j \exp\bigl[ - C
m^{1/3} \bigr] \leq j \exp\bigl[ - C m^{1/3} \bigr],\nonumber
\end{eqnarray}
where Lemma \ref{lemBin-norm} is employed to bound $P^*(A^c_1)$ and
$P^*(A^c_\ell|U_1,\ldots, U_{\ell-1})$, and
we use (\ref{m-T}) to bound $m-T_{\ell-1}$.

Plugging (\ref{EAr}) and (\ref{PAc}) into (\ref{rPQprime}) and
combining it together with (\ref{rP*Q3}) and (\ref{rQQprime}), we obtain
\begin{eqnarray*}
H\bigl(P^*, Q\bigr) &\leq&\frac{ C (r-1) }{\sqrt{m}} + \sum_{j=1}^r
\biggl\{2 j \exp\bigl[ - C m^{1/3} \bigr] + \frac{C}{m} \biggr
\}^{1/2} \\
&\leq&\frac{C r }{\sqrt{m}} + r^2 \exp\bigl[ - C
m^{1/3} \bigr],
\end{eqnarray*}
which proves the lemma for the $r+1$ case.
\end{pf*}

%Define
%&& A^c = \bigcup_{j=1}^r A_j^c = \bigcup_{j=1}^r (A^c_j A_{j-1}\cdots
%A_1) \\
%&& P^*(A^c) = \sum_{j=1}^r P^*(A^c_j A_{j-1}\cdots A_1 ) =P^*(A_1^c)+
%U_{j-1})] \\
%&&\qquad \leq\exp\left[ - C m^{1/3} \right] +\sum_{j=2}^r E_{P^*}\left(A_1
%&&\qquad \leq\sum_{j=1}^r \exp\left[ - C m^{1/3} \right] \leq r \exp\left[
%- C m^{1/3} \right].
%$r=4$,
%p_1 (1-p_1)]^{2/3} } \leq\frac{C}{m } \]
%r case
%1_{A_1\cdots A_{r-1}} E \left[ 1_{A_r} \log\frac{dP^*_r}{dQ_r}
%| U_1,\ldots, U_{r-1} \right] \right\}\]

%Suppose that for $k=1,\ldots, n$,
%$P_k$ is a multinomial distribution $\cM(m, \theta_1,\ldots, \theta_{
%$\theta_1 + \cdots+ \theta_{\nu_k}=1$, and for constants $c_0$ and
%$c_1$,
%Denote by $Q_k$ the multivariate normal distribution whose mean and
%covariance are the same as
%$P^k$. If $\nu_k \geq2$, following the same way as in Lemma
%and independent uniform distributions
%on $(-1/2, 1/2)$, and if $\nu_k \leq1$ let $P_k^*=P_k$. Assume that
%$P_k, P_k^*, Q_k$ for different $k$ are independent, and define
%product probability measures
%Then we have

\begin{pf*}{Proof of Lemma \ref{lemMany-Multi-norm}}
Since $P_k, P_k^*, Q_k$ for different $k$ are independent, an
application of the Hellinger distance property for product probability
measures [\citet{LeCYan00}] leads to
\[
H^2\bigl(P^*, Q\bigr) \leq\sum_{k=1}^n
H^2\bigl(P^*_k, Q_k\bigr).
\]
We note that if $\nu_k\leq1$, both $P_k$ and $Q_k$ are point mass at
$m$ and thus $H(P_k, Q_k)=0$. Hence,
\[
H^2\bigl(P^*, Q\bigr) \leq\sum_{k=1}^n
H^2\bigl(P^*_k, Q_k\bigr) 1(
\nu_k\geq2).
\]
Applying Lemma \ref{lemMulti-norm}, we obtain
\[
H^2\bigl(P^*, Q\bigr) \leq\sum_{k=1}^n
\biggl[\kappa^4 \exp\bigl( - C m^{1/3} \bigr) +
\frac{C \kappa^2}{m} \biggr] 1(\nu_k\geq2).
\]
For $m$ exceeding certain integer $m_0$,
\[
\frac{C \kappa^2 }{m} \geq\kappa^4 \exp\bigl( - C m^{1/3}
\bigr)
\]
and hence for $m >m_0$,
\[
H^2\bigl(P^*, Q\bigr) \leq\frac{C \kappa^4}{m} \sum
_{k=1}^n 1(\nu_k\geq2).
\]
For $m \leq m_0$, we may adjust constant $C$ to make the above
inequality still holds for $m \leq m_0$.
\end{pf*}

\begin{pf*}{Proof of Lemma \ref{lemTV}}
%(\ref{TV1}) can be easily shown as follows,
%&& \| F_1 \times F_2 - G_1 \times G_2 \|_{\mathrm{TV}} \leq\| F_1 \times F_2 -
%G_1 \times F_2 \|_{\mathrm{TV}} + \| G_1 \times F_2 - G_1 \times G_2 \|_{\mathrm{TV}} \\
%&&\qquad = \| F_1 - G_1 \|_{\mathrm{TV}} + \| F_2 - G_2 \|_{\mathrm{TV}}.
%For (\ref{TV2}) we have
%
\begin{eqnarray*}
\| F - G \|_{\mathrm{TV}} &=& \bigl\|F_1(x) \times
F_{2|1}(y|x) - G_1(x) \times G_{2|1}(y|x)
\bigr\|_{\mathrm{TV}}
\\
&\leq& \bigl\|F_1(x) \times F_{2|1}(y|x) - F_1(x)
\times G_{2|1}(y|x) \bigr\| _{\mathrm{TV}} \\
&&{}+ \bigl\|F_1(x) \times
G_{2|1}(y|x) - G_1(x) \times G_{2|1}(y|x) \bigr\|
_{\mathrm{TV}}
\\
&=& \bigl\|F_1(x) \bigl[F_{2|1}(y|x) - G_{2|1}(y|x)
\bigr] \bigr\|_{\mathrm{TV}} \\
&&{}+ \bigl\|F_1(x) G(x,y)/G_2(x) - G(x,
y) \bigr\|_{\mathrm{TV}},
\end{eqnarray*}
where
\begin{eqnarray*}
&&\bigl\|F_1(x) \bigl[F_{2|1}(y|x) - G_{2|1}(y|x)
\bigr] \bigr\|_{\mathrm{TV}} \\
&&\qquad= E_{F_1} \bigl[\bigl\| F_{2|1}(
\cdot|U_1) - G_{2|1}(\cdot|V_1)
\bigr\|_{\mathrm{TV}} |U_1=V_1\bigr],
\\
&&\bigl\|F_1(x) G(x,y)/G_2(x) - G(x, y) \bigr\|_{\mathrm{TV}} \\
&&\qquad=
\bigl\|\bigl[F_1(x) / G_1(x) -1\bigr] G(x, y)
\bigr\|_{\mathrm{TV}}
\\
&&\qquad\leq\max_x \biggl\{ \biggl\llvert\frac{P(U_1=x)}{P(V_1=x)} -1
\biggr\rrvert\bigl\|G(x, y) \bigr\|_{\mathrm{TV}} \biggr\} = \max_x
\biggl\llvert\frac
{P(U_1=x)}{P(V_1=x)}-1 \biggr\rrvert.
\end{eqnarray*}
\upqed
\end{pf*}
\end{appendix}


%Yazhen Wang is Professor, Department of Statistics, University of
%Wisconsin-Madison, Madison, WI 53706, USA. Email: yzwang@stat.wisc.edu.
%The authors thank the editor, associate
%editor, and two anonymous referees for stimulating comments and
%suggestions, which led to significant improvements in both substance
%and the presentation of the paper. }

%suskaldyti doi

% imsref loaded by lrinkeviciute, 2013-09-19 10:58:01
%
% imsref loaded by lrinkeviciute, 2013-09-19 12:39:58
%

\printaddresses

\end{document}